\documentclass[12pt,a4paper]{article}
\usepackage{times}
\usepackage[norefpage]{nomencl}
\makenomenclature
\usepackage[utf8]{inputenc}
\usepackage{graphicx}
\usepackage{flushend}
\usepackage{caption}
\usepackage{subcaption}
\usepackage{bbm}
\usepackage{epsfig,epstopdf}
\usepackage{cite}
\usepackage{float}
\usepackage[usenames,dvipsnames]{color}
\usepackage{amssymb,amsmath,amsbsy}
\usepackage{setspace}
\usepackage{bm}
\usepackage{multicol}
\usepackage[margin=25mm]{geometry}
\usepackage{mathtools}
\usepackage{multicol}
\usepackage{soul} 
\usepackage{framed}
\usepackage[colorlinks=true,
            linkcolor=blue,
            citecolor=blue,
             urlcolor=blue,
            filecolor=blue,
          pdfpagemode=None,
          pdfstartview=FitH]{hyperref}
\makeatletter
\newcommand*{\rom}[1]{\expandafter\@slowromancap\romannumeral #1@}
\makeatother
\newcommand\comb[2][n]{\prescript{#1\mkern-0.5mu}{}C_{#2}}

\title{ABC Conjecture: $ABC = 2^m p^n q^r$ with Fermat or Mersenne Primes}

\author{Anupam Saxena\thanks{to $A, B, C, E$ and most importantly, the two $Ks$.}
\\ Indian Institute of Technology Kanpur, INDIA 208016\\ \copyright{Anupam Saxena}}

\begin{document}

\maketitle
\begin{abstract}
\label{label_abstract}
For $p$ and $q$ any two distinct Fermat or Mersenne primes, $m,n,r$ as positive integers and $\mu = \pm 1$ satisfying any diophantine relation, $\mbox{(i)}\; 2^m + \mu = p^nq^r, \mbox{(ii)} \; 2^mp^n  + \mu = q^r \mbox{ or } \mbox{(iii)} \; p^n + \mu q^r = 2^m$, it is shown that the number of triplets $\{A, B, C \}$ with $\gcd(A,B) = 1$ and $C = A + B$, for which their product is of the form $ABC = 2^mp^nq^r$ and which satisfy $C > \mathrm{rad}(ABC)^{1 + \varepsilon}$ for any real $\varepsilon > 0$, is finite.

For the triplet $\{2^{y+1}, 2^{2y}+1, (2^y+1)^2\}$, a solution to (iii) with positive integer $y$ such that $2^y+1$ and $2^{2y}+1$ are primes, $\mathrm{rad}(ABC)^{1 + \varepsilon} > C$ holds for any $\varepsilon > 0$. Furthermore, finiteness of the number of solutions of (iii) when $n$ is even, is demonstrated in \cite{Bugeaud_1997}. All other solutions are enumerated.

\end{abstract}
\noindent\rule{16cm}{1pt}

\section{Introduction}

\noindent The $ABC$ conjecture states that for every $\varepsilon > 0$, there exist finitely many triplets $\{A, B, C\}$ of natural numbers, satisfying $A + B = C, gcd(A,B) = 1$, and $B > A \geq 1$, such that $C > rad(ABC)^{1 + \varepsilon}$ holds, where $rad(ABC)$ is the product of all distinct primes of $ABC$. \\

\noindent One realizes that product $N = ABC$ is even, and can be written in the form $N = 2^mx^ny^r$ where $x$ and $y$ are any odd positive natural numbers, and coprime to each other. $m, n, r$ are positive integer exponents. To satisfy $C = A + B$, six diophantine combinations are possible\footnote{$C = 2^m = x^n y^r + 1 \equiv A + B$;  $C = x^n y^r = 2^m + 1 \equiv A + B$; $ C = x^n =  2^m + y^r $; $C = 2^m =  x^n + y^r $; $C = 2^m x^n =  y^r + 1 $  and $C \equiv y^r = 2^m x^n + 1$.}: 

\begin{eqnarray}
	\label{equality_3_p0}                  \mbox{ (a)   }     2^{m} \pm 1 = x^{n}y^{r};  \;\;\;\;\;\;\;\;\;\; 
	\mbox{ (b)   }    2^{m} = x^{n} \pm y^{r}    \;\;\;\;\;\;\;\;\;\;                    
	\mbox{ (c)   }    2^{m}x^{n} \pm 1 = y^{r};                                                      
\end{eqnarray}
\\

\noindent Each diophantine relation above is addressed respectively in \S\ref{2m_plus_mu_equals_pnqr}-\S\ref{2m_pm_plus_mu_equals_qr} for specific cases when $x\equiv p$ and $y \equiv q$ are Fermat and/or Mersenne primes. Using existing and other elementary approaches, we show for such cases that the ABC conjecture holds and that triplets $\{A, B, C \}$ satisfying $C > rad(ABC)^{1 + \varepsilon}$ are finitely many in number. We further attempt to enumerate such solutions by solving the above diophantine equations. For a solution to Eq. \ref{equality_3_p0} (b) (\S\ref{pn_mu_qr_2m}), that is, $2^{z+1} = (2^{z} + 1)^2 - (2^{2z} + 1)$, $z$ a positive integer such that $2^{z} + 1$ and $2^{2z} + 1$ are primes, $rad(ABC)^{1 + \varepsilon} > C$, that is, $[2 \cdot (2^{z} + 1) \cdot (2^{2z} + 1)]^{1 + \varepsilon} > (2^{z} + 1)^2$ holds for any $\varepsilon > 0$. For another case of Eq. \ref{equality_3_p0} (b) when $n$ is even, number of solutions are shown to be finite\cite{Bugeaud_1997}. \\

\noindent
The $ABC$ conjecture is attributed to Masser \cite{Masser85} and Oesterl\'{e} \cite{Oesterle88}. Many other forms exist, e.g., \cite{Dorian_Goldfeld,wiki_abc}. Significant effort has been devoted into showing the conjecture to be true for special and generic cases, by Mochizuki and others, e.g., \cite{wiki_abc, Mochizuki1, Mochizuki2, Mochizuki3, Mochizuki4, Mochizuki5, Mochizuki6, Mochizuki7}. The conjecture has numerous consequences, and is affiliated strongly with the Fermat's Last Theorem, proved in \cite{Wiles}. Upper bounds on $C = A + B$, are established in \cite{ Stewart_Tijdeman,Stewart_Yu1,Stewart_Yu2}. Waldschmidt \cite{Waldschmidt_2009}, in a comprehensive review, relates Pillai's conjecture \cite{Pillai_SS_1945} (proved in \cite{Stroeker_and_Tijdeman_1982}), Catalan Equation (proved in \cite{Catalan_Mihailescu}), Waring's problem and the ABC conjecture. Da silva et. al. \cite{da_Silva_et_al_2018}, inspired by Lenstra \cite{Lenstra_talk}, re-interpret Gersonides’s theorem as the $ABC$ conjecture being true on the set of harmonic\footnote{a harmonic number is a power of $2$ times a power of $3$.} numbers. They further add to this set, one at a time, a finite set of $ndh$-numbers\footnote{An $ndh$-number is a positive integer that cannot be written as a difference of harmonic numbers.}, the infinite set of primes congruent to $ 41 \; (\bmod \; 48)$, sets of Fermat and Mersenne primes, and show that the $ABC$ conjecture holds on these sets. \\

\noindent Of the many existing, some recent works on the $ABC$ conjecture comprise the following. Lichtman \cite{Lichtman_2025} demonstrates via an elementary and self-contained proof that there are at most $O(N^{\frac{2}{3}})$ many triplets of coprime integers in a cube ($a, b, c$ $\in \{1, ..., N \}^3$) that satisfy the summation constraint and $\mathrm{rad}(abc) < c^{1 - \varepsilon}$. Idowu \cite{Idowu_2025} proposes novel parametric constructions of high-quality abc-triplets using powers of 2 and 3 with modular inversion, and offers algebraic and computational evidence of structured triplets optimized for the ratio $\log c / \log \mathrm{rad}(abc)$. Bright \cite{Bright_2024} shows that there exist infinitely many triples such that $C > \mathrm{rad}(ABC) \exp(6.563 \sqrt{\log C} / \log \log C)$ improving constant from earlier works. Scoones \cite{Scoones_2022} generalises the method of Stewart and Yu to give an improvement on Györy's bound for algebraic integers.

\noindent \textcolor{black}{\rule{16cm}{1pt}}

\section{$2^m + \mu = p^nq^r$}
\label{2m_plus_mu_equals_pnqr}
\noindent Consider $2^{m} + \mu = p^{n}q^{r}$ with $\mu = \pm 1$,  $p,q$ odd primes and $m,n,r$ positive integers. 
$gcd(n,r)$ must be $1$, otherwise, $2^{m} + \mu$ is a perfect power, an impossibility, from Appendix \S \ref{2_k_pm_mu}. \\ Three cases arise: $m$ is composite (\S \ref{Kaali_153_a1}), an odd prime (\S \ref{Kaali_153_a2}), or a power of 2 (\S \ref{2_2_w_plus_mu_pn_qr}). \\
 

\subsection{$2^{av} + \mu = p^nq^r$}
\label{Kaali_153_a1}
Consider $m = av$, composite with $a > 1$ and $v$, an odd prime.
\begin{eqnarray} 2^m + \mu = 2^{av} + \mu^v = (2^a + \mu)R =  (2^a + \mu)\left[1 - 2^a\mu + 2^{2a} -  2^{3a}\mu + ... +  2^{a(v-1)}\right] \nonumber \end{eqnarray}
For $v > 2$ and $a > 1$, $R = \frac{2^{av} + \mu}{2^a + \mu} \geq 2^a + \mu > 1$   since $2^{a(v-1)} \geq 2^a + 2 > 2^a - 2 + \frac{2}{2^a}$. Equality holds  when  $a = \mu = 1$ and $v = 3$. \\

\noindent Per Appendix \S \ref{gcd_of_factors2}, $gcd(2^a + \mu,R) = gcd(2^a + \mu,v) = 1 \mbox{ or } v$. 

\begin{enumerate}
\item Consider $gcd(2^a + \mu, R) = gcd(2^a + \mu, v) = 1$.\\
Since $2^a+\mu \mbox{ and } R$ are both $> 1$ and coprime, and since $(2^a + \mu)R = p^nq^r$, without  
loss of generality, let $2^a + \mu = p^n$ and $R = q^r$. \\
Per Appendix \S \ref{2_k_pm_mu}, $n = 2$ when $a = p = 3 \mbox{ and } \mu = 1$; Else, $n = 1$.

\begin{enumerate}
\item For $n = 2$, $v \neq 3$ since $gcd(2^a + \mu, v) = 1$. $v$ must be $\geq 5$. Then
\begin{eqnarray}
q^r = \frac{2^{3v}  +1}{9} = \left(\frac{2^{v}  +1}{3}\right) \left(\frac{1 - 2^{v}  +2^{2v}}{3}\right) \nonumber
 \end{eqnarray} 
Since $v$ is odd, $3 \mid 2^{v}  +1$. $1 - 2^{v}  + 2^{2v} \equiv 1 - (-1)^v + (-1)^{2v} \equiv  3 \text{ or } 0 \;(\bmod \;3)$ so that $3 \mid 1 - 2^{v}  + 2^{2v}$. \textcolor{black}{Per Appendix \S\ref{gcd_of_factors2}, $gcd(2^{v}  +1, 1 - 2^v + 2^{2v}) = gcd(2^{v}  +1, 3) = 3$ implying $gcd\left[ \left(\frac{2^{v}  +1}{3}\right), \left(\frac{1 - 2^{v}  +2^{2v}}{3}\right)\right] = 1$}. Since $v \geq 5$, $\frac{2^{v}  +1}{3} > 1$ and $\frac{1 - 2^{v}  +2^{2v}}{3}> 1$. $q^r = \frac{2^{3v}  +1}{9}$ does not hold since the $RHS$ is divisible by two or more odd primes. \\

\item For $n = 1$, $2^a + \mu = p$ (prime). $ q^r = R = \frac{(2^a)^v + \mu}{2^a + \mu}$.

\begin{enumerate}
	\item \textcolor{black}{From \cite{Ge_Yimin}, Corollary 6, for a prime number $p$ and an(y) integer $x$, every prime divisor $q$ of $\frac{x^p-1}{x-1} = 1 + x + x^2 + ... + x^{p-1}$ is either $p$ or $\equiv 1  \; (\mod p) $}.
	\item \textcolor{black}{For $v$ prime and $\mu = -1,$ the above holds true for $R = \frac{(2^a)^v -1 }{2^a -1}$. The same also is true for $\frac{(2^a \mu)^v + \mu}{2^a \mu  + \mu}$ and thus $R = \frac{(2^a)^v + 1}{2^a  + 1}$}.
	\item $(2^a)^v + \mu = 2^a[(2^a)^{v-1} - 1] + 2^a + \mu$. $v \mid (2^a)^{v-1} - 1$  per the Fermat's little theorem. $v \nmid 2^a + \mu$ since $gcd(2^a + \mu,v) = 1$. \textcolor{black}{Thus, $q$ must be $\equiv 1 \; (\bmod \; v)$}. 
	\item Let $q$ be a prime of the form $q = 2^d + \beta, \beta = \pm 1$ and $d \geq 1$, an integer.
	\item Further, if $q = 2^d + \beta = \textcolor{black}{\lambda v + 1}$, $\lambda$ a positive even integer, $\beta$ cannot be $1$ since $v \nmid 2^d$.  $\beta$ must be $-1$ such that $q = 2^d - 1$ is prime implying $d$ is either $2$ or an odd prime. \\

\item {\bf Proposition}: For $\mu = -1$, $q^r = (2^d - 1)^r = \frac{(2^a)^v + \mu}{2^a + \mu}$ does not hold. \\
{\bf Proof}:
\begin{enumerate} 
\item Since $p = 2^a - 1$ is prime, $a$ is either $2$ or an odd prime.
\item If $a = 2$, $ q^r = (2^d-1)^r = \frac{(2^a)^v -1}{2^a -1} = \frac{2^{2v}-1}{3} = (\frac{2^v+1}{3})(2^v-1)$ does not hold since $\frac{2^v+1}{3} > 1$ and $gcd(2^v-1,2^v+1) = 1$. $(\frac{2^v+1}{3})(2^v-1)$ is divisible by two or more odd primes. $a$ must be an odd prime.
\item Further, $d \neq 2$ since $2^{av} - 1 \equiv (-1)^{av} - 1 = -2 \; (\bmod \; 3)$. $d$ must be odd prime. 
\item \textcolor{black}{If $2^d-1 \mid 2^{av}-1$, $d$ must divide\footnote{For positive integers $m$ and $n$, $2^m - 1$ and $2^n - 1$ are coprime if and only if $m$ and $n$ are coprime \cite{Thms_Mersenne_primes}. In general, $gcd(a^m-1,a^n-1) = a^{gcd(m,n)} - 1$ for positive integers $a, m, n$.} $av$}. Since $d, a \mbox{ and } v$ are all odd primes, either $d = a$ or $d = v$. $d$ cannot be $a$ since $p \neq q$. $d$ must be $v$ making $2^v-1$ prime. \\
$ q^r = \frac{(2^a)^v -1}{2^a -1}$ may be rewritten as $(2^a -1)(2^v-1)^r = (2^a)^v -1$. 
\item Per Appendix \S \ref{gcd_of_factors2}, $gcd(\frac{ 2^{av} -1}{2^v-1},2^v-1)$ = $1$ or $a$. 

(s1) If $gcd(\frac{ 2^{av} -1}{2^v-1},2^v-1) = 1$, $2^v-1 \nmid \frac{ 2^{av} -1}{2^v-1}$. $r$ cannot be larger than 1. \\
(s2) If $r = 1$, $(2^a-1)(2^v-1) = 2^{av} - 1$, or, $2^{a+v-1} - 2^{a-1} - 2^{v-1} = 2^{av-1} - 1$ is not possible due to opposite parities. \\
(s3) If $gcd(\frac{ 2^{av} -1}{2^v-1},2^v-1) = a$, $2^v - 1$ must be $a$. $q^r = \frac{(2^a)^v -1}{2^a -1}$  becomes $(2^{2^v-1} -1)(2^v-1)^r = 2^{v(2^v-1)} -1$. \\
\textcolor{black}{
(s4) Per Bang's theorem (\cite{Bang1}, \cite{Bang2}; 
Zsigmondy's theorem \cite{Zsigmondy_1892} is its generalization), for positive integers $n, k$ and for $n > 1, n \neq 6$, there exists a prime number (primitive prime divisor) that divides $2^n - 1$ but not $2^k - 1$ for any $k < n$.\\
(s5) $2^{2^v-1} -1 \mid 2^k - 1$ for $k = 2^v-1 < n = v(2^v-1)$. \\
(s6) $2^{v} -1 \mid 2^k - 1$ for $k = v < n = v(2^v-1)$. \\
(s7) $(2^{2^v-1} -1)(2^v-1)^r = 2^{v(2^v-1)} -1$ is therefore not possible since there exists an odd prime\footnote{This odd prime is likely $\equiv 1 \; (\bmod \; v(2^v - 1) )$ (\cite{Birkhoff_Vandiver_1904}, Theorem V).} other than $2^{2^v-1} -1$ and $2^v-1$, that divides $2^{v(2^v-1)} -1$}.  $\blacksquare$ \\ 
\end{enumerate}

\item {\bf Proposition}: For $\mu = 1$, $q^r = (2^d - 1)^r = \frac{(2^a)^v + \mu}{2^a + \mu}$ does not hold. \\
{\bf Proof}: \\

\begin{enumerate}
\item Since $p = 2^a + 1$ is prime, $a \mbox{ } (> 1)$  must be $2^w, w \geq 1$ making $p$ a sum of two squares.
\item $2^{av}+\mu = (2^v)^{2^w} + 1$ is also a sum of two squares with $gcd[(2^v)^{2^w},1] = 1$.
\item \textcolor{black}{Per Euler \cite{Euler_maa}, Proposition IV, if two squares are prime between themselves, their sum cannot be divided by any number that itself is not a sum of two squares}.
\item \textcolor{black}{Per the Fermat's theorem on sums of two squares (\cite{Fermat_thm_sum_two_sqs1},\cite{Fermat_thm_sum_two_sqs2}), an odd prime $q$ can be expressed as $q = x^2 + y^2$ with $x$ and $y$ integers, if and only if $q \equiv 1 \; (\bmod 4 )$.	
}
\item $2^d - 1 = 4k + 1$ is an impossibility since $2 \nmid 2^{d-1}-1$ for $d \geq 2$. $\blacksquare$ 
\end{enumerate}
\end{enumerate}
\end{enumerate}

\item Consider $gcd(2^a + \mu,R) = gcd(2^a + \mu,v) = v$.

\begin{enumerate}
\item \label{kaali_M1} For positive integers $\lambda$ and $\beta$, let $2^a + \mu = \lambda v$, $R = \beta v$ such that $gcd(\lambda,\beta) = 1$. Either $\beta > \lambda$,  or $\beta = \lambda = 1$. Consider cases\footnote{In a special case, when $a = \mu = 1$ and $v = R = 3$, $\lambda = \beta = 1$. $2^{av} + \mu = 2^{3} + 1 = 9$, not permitted since $a > 1$.}  when $\beta > \lambda$.
\begin{eqnarray}
R & = & \beta v = \frac{(2^a)^v + \mu}{2^a + \mu} = \frac{(\lambda v - \mu)^{v} + \mu }{\lambda v} \;\;\;\;\;\;\;\;\;\;\;\;\;\;\;\; \mbox{                         or } \nonumber \\ \nonumber \\
R & = & (\lambda v)^{v-1} - \comb[v]{1}(\lambda v)^{v-2} \mu +  \comb[v]{2}(\lambda v)^{v-3} - ... -  \comb[v]{v-2}(\lambda v)\mu + v = (G+1)v \nonumber
 \end{eqnarray}

so that $\beta = G + 1$. 

\item 
\label{kaali_M1igmr}
\textcolor{black}{For $v$ an odd prime, $v \mid \comb[v]{i}, 1 \leq i \leq v-1$. $v$ divides $(\lambda v)^{v-j}$, $j = 0, ..., v-1$ at least once and so $v \mid G$ implying $gcd(\beta,v) = 1$}.

\item $(2^a)^v + \mu = (2^a + \mu)R = \lambda \beta v^2 = p^n q^r$. $\lambda \beta v^2$ must comprise of powers of 
precisely two distinct odd primes. As $gcd(v,\beta) = gcd(\lambda,\beta) = 1$ and $\beta > \lambda \geq 1$ 
(items \ref{kaali_M1}, \ref{kaali_M1igmr}), $\lambda$ is either 1 or a power of $v$. Further, as $\lambda v = 2^a + \mu$, per Appendix \S \ref{2_k_pm_mu}, either $\lambda = v$ or $\lambda = 1$. \\

\item \label{concl1} Case I: $\lambda = v$: $v = 3, a = 3, \mu = 1$ must hold per Appendix \S\ref{2_k_pm_mu} making $m = av = 9$. $2^{av} + 1 = \boxed{2^9 + 1 = 3^3 \cdot 19} = p^nq^r$. Let $p = 19 \mbox{ and } q = 3 \mbox{ so that } n = 1 \mbox{ and } r = 3$. 
$C = 3^3 \cdot 19 > rad(N)^{(1 + \varepsilon)} = (2pq)^{(1 + \varepsilon)} = (2\cdot 19 \cdot 3)^{(1 + \varepsilon)}$ holds for $0 < \varepsilon < 0.3176$. \\  

\item \label{Kaali_M3} Case II: $\lambda = 1$: Without loss of generality, let $p = v= 2^a + \mu$ ($ \neq 3^2$) so that $n = 2$. Then, $q^r = \beta  =  \frac{2^{av} + \mu}{v^2}$.

\begin{enumerate}
\item $r$ cannot be even as it requires $2^{av} + \mu$ to be an even power, not possible for $av > 3$, per Appendix  \S\ref{2_k_pm_mu}. \\

\item {\bf Proposition}: If $q$ (prime) is of the form $q = 2^d + \gamma, \gamma = \pm 1$, $d \mbox{ (integer) } \geq 1$, and $v = 2^a + \mu$ (prime), in $q^r =  \frac{2^{av} + \mu}{v^2}$, if $r$ is odd, $r$ must be 1, yielding only one solution. \\

{\bf Proof}: \\
\begin{enumerate}
\item \label{kaali_Maa2} If $r$ is odd, $(2^d + \gamma)^r (2^a + \mu)^2 = 2^{av} + \mu$  yields 
 $2^{a+1}R_1 + 2^d R_2 = 2^{av} + \mu -\gamma$  where $R_1 = (2^{a-1} + \mu)(2^d + \gamma)^r$, odd,  and \\ 
 $R_2 = [(2^d)^{r-1} + \comb[r]{1}(2^d)^{r-2}\gamma + \comb[r]{2}(2^d)^{r-3} + ... + r]$, odd. 
\item If $\gamma = -\mu$,  as $a > 1$, $2^{a}R_1 + 2^{d-1} R_2 = 2^{av-1} + \mu$ holds only when $d = 1$. \textcolor{black}{Otherwise, the parities on both sides are different}. Then, $q = 2 + \gamma$ must be $3$ implying $\gamma = 1$. $\mu = -\gamma = -1$. \\
$q^r =  \frac{2^{av} + \mu}{v^2}$ becomes $3^r(2^a - 1)^2 = 2^{av}-1$, an impossibility, since $2^{av}-1 \equiv (-1)^{av} - 1 = -2 \; (\bmod \;3)$ noting  
 that $a$ must be an odd prime\footnote{$a$ cannot be 2 as then $p = 3$, not possible, since $q = 3$.} for $p = v = 2^a - 1$ to be prime. 
\item If $\gamma = \mu$ in item (\ref{kaali_Maa2}) above, $2^{a+1}R_1 + 2^{d} R_2 = 2^{av}$ holds when\footnote{$av > d$ and $av > a+1$ must hold. Further, if $d \neq a + 1$, different parities result. } $av > d = a+1$. Then,  $p = v = 2^a + \mu$ and $q = 2^{a+1} + \mu$. 
\item If $\mu = 1$, $a$ and $a+1$, must both be powers of 2, an impossibility since $a > 1$.
\item \label{concl2} If $\mu = -1$, $a$ and $a+1$, must  either be 2 or odd primes. Only possibility is $a = 2$ which makes $p = 3$ and $q = 7$. $v = 3$ per item (\ref{Kaali_M3}). $2^{av} + \mu = v^2 q^r \equiv p^n q^r$ becomes $\boxed{ 2^6-1 = 7^r \cdot 9}$. Thus, $r = 1$ making $A = 1, B = 7 \cdot 9, \text{ and } C = 2^6$.  $\blacksquare$
\end{enumerate}
\end{enumerate}
\end{enumerate}
\end{enumerate}


\subsection{$2^{v} + \mu = p^nq^r$}
\label{Kaali_153_a2}
In Eq. \ref{equality_3_p0} (a), consider $m = v$, an odd prime. \\
\noindent  For $\mu = 1$,
\begin{eqnarray}
2^v + 1 = (1 + 2)R = (1 + 2)(1 - 2 + 2^2 + ... + 2^{(v-1)}) = p^nq^r \nonumber
\end{eqnarray}

\noindent Per Appendix \S \ref{gcd_of_factors2}, $gcd(1+2, R) = gcd(3, v) = 1 \mbox{ or } v$.

\begin{enumerate} 
\item Let $gcd(3,R) = gcd(3,v) = 1$. \\
$3 \mid 2^v + 1$. As $2^v + 1$ is divisible by only two distinct primes $p \mbox{ and } q$,  let $p = 3$. Then, $n = 1$ noting that $gcd(3,R) = 1$ and so, $R = q^r = \frac{2^v + 1}{2 + 1}$.\\ 

\item {\bf Proposition}: In $q^r =  \frac{2^{v} + 1}{2+1}$, $r > 1$ is not possible if $q$ (prime) is of the form $2^d + \beta, \beta = \pm 1$. \\
{\bf Proof}:  
\begin{enumerate} 
\item $q^r =  \frac{2^{v} + 1}{2+1}$ becomes $3(2^d + \beta)^r = 2^v + 1$, or, $2R_1 + 2^d R_2 = 2^v + 1 - \beta^r$ where $R_1 = (2^d + \beta)^r$, which is odd, and \\
$R_2 = [(2^d)^{r-1} + \comb[r]{1}(2^d)^{r-2}\beta + \comb[r]{2}(2^d)^{r-3}\beta^2 + ... +  \comb[r]{r-2}(2^d)\beta^{r-2} + r\beta^{r-1}$. 
\item  If $\beta^r = 1$, $R_1 + 2^{d-1} R_2 = 2^{v-1}$ holds only when $2^{d-1}R_2$ is odd, that is, 
$d = 1$ and $R_2$ is odd. Then, $r$ must be odd, implying that $\beta = 1$. $q = 3$ is not possible, since $p = 3$. 
\item  If $\beta^r = -1$, which is possible only when $\beta = -1$ and $r$ odd, $2R_1 + 2^d R_2 = 2^v + 1 - \beta^r$ becomes $R_1 + 2^{d-1} R_2 = 2^{v-1} + 1$. Since $R_1$ is odd, and $R_2$ is odd (as $r$ is odd), $d \mbox{ must be } > 1$. $3(2^d + \beta)^r = 3(2^d -1 )^r = 2^v + 1$, yields $3 \cdot 2^d \cdot R_2 = 4 (2^{v-2} + 1)$ for which to hold, $d$  must  be $2$. Then, $q = 2^2 - 1 = 3$, not possible, since $p = 3$. $\blacksquare$
\end{enumerate}

\item Next, consider $gcd(3,R) = gcd(3,v) = v$. $v$ must be 3. $2^v + 1 = 9 \neq p^nq^r$. \\
\end{enumerate} 

\noindent For $\mu = -1$, 
\begin{eqnarray}
2^v - 1 = 1 + 2 + 2^2 + ... + 2^{(v-1)} = p^nq^r \nonumber
\end{eqnarray}

\begin{enumerate}
\item Consider distinct primes $p$ and $q$ of the forms $2^g + \gamma$ and $2^d + \beta$ respectively, where $\beta,\gamma = \pm 1$, and $d,g$ are positive integers such that $p, q \geq 3$. 
\item Every prime divisor of $2^v - 1 = \frac{2^v - 1}{2-1}$ is either $v$ or $\equiv 1 \; (\bmod \; v)$ (\cite{Ge_Yimin}, Corollary 6). \\ 
$2^v - 1 = 2(2^{v-1}-1) + 1$. As $v \mid 2^{v-1}-1$ per Fermat's little theorem, $v \nmid 2^v - 1$. Any prime divisor of $2^v - 1$ is then $\equiv 1 \; (\bmod \; v)$. 
\item If $p = 2^g + \gamma = \lambda v + 1$ for a positive (even) integer $\lambda$, $\gamma$ cannot be 1 as $v \nmid 2^g$. 
\item $\gamma$ must be $-1$, and likewise, $\beta = -1$ making $p = 2^g - 1$ and $q = 2^d - 1$, $g \neq d$. 
\item $g, d$ must either be 2 or odd primes.  $2^v - 1 \equiv (-1)^v - 1 = -2 \; (\bmod \;3)$ implying $3 \nmid 2^v-1$. $g, d$ must, therefore, be odd primes. 

\item \label{mers_p_div_mers_p} \textcolor{black}{ 
{\bf Proposition}: If $v$ is an odd prime, $2^d - 1$ divides $2^v - 1$ only when $d = v$. \\
{\bf Proof}: Let $v \geq d$. Then, $gcd(2^d - 1,2^v-1) = gcd(2^d - 1,2^d[2^{v-d}-1]) =  gcd(2^d - 1,2^{v-d}-1) = ... = gcd(2^d - 1,2^{v - jd}-1)$. $j$ is such that $0 \leq v - jd < d$, or, $2^{v-jd} - 1 < 2^d - 1$. $gcd(2^d - 1,2^{v - jd}-1) = 2^d - 1$ only if  $v = jd$. That is, $d \mid v$. Since $v$ and $d$ are both, odd primes, $d = v$ so that $j = 1$. \\ \\
Alternatively, $gcd(2^v-1,2^d-1) = 2^{gcd(v,d)} - 1$ from \cite{Thms_Mersenne_primes}. If $gcd(2^v-1,2^d-1) = 2^d-1$, $gcd(v,d) = d$ must hold. Since $v$ and $d$ are both primes, $d$ must be $v$.
} $\blacksquare$

\item Thus, if $d = v$, $r = 1$ and $n = 0$, not permitted since $n$ is a positive integer. Similar is the case for $g = v$. 
\end{enumerate}


\subsection{$2^{2^w} + \mu = p^nq^r$}
\label{2_2_w_plus_mu_pn_qr}
In Eq. \ref{equality_3_p0} (a), let $m = 2^w, w \geq 1$. \\ \\
\noindent {\bf Proposition}: For $\mu = 1$, $p$ (or $q$) cannot be a Mersenne or Fermat prime. \\
{\bf Proof}: 
 
\begin{enumerate}
\item Let $p=2^d + \beta$, $d \geq 1, \beta = \pm 1$. $d > 1$ if $\beta = -1$.

\item \textcolor{black}{Any prime divisor of $2^{2^w} + 1$, $w > 1$, is of the form $M2^{w+2} + 1$ for some positive integer $M$ (Lucas's \cite{Lucas_E} revision to Euler's $M2^{w+1} + 1$)}. 
\item If $p = 2^d + \beta =  M 2^{w+2} + 1$, $\beta \neq -1$ since $ M 2^{w+1} + 1$ and $2^{d-1}$ have different parities for $d > 1$. 
\item $\beta$ must be 1 so that $p = 2^d + 1$. $p$ being prime, $d$ must be $2^z, w> z \geq 0$. If $z = w$, then $n = 1$ and $r = 0$, not permitted as $r$ is a positive integer. 
\item \label{Fermat_identity} $2^{2^z} + 1 \nmid 2^{2^w} + 1$ for any $z < w$. \\
Let $F(y) : = 2^{2^y} + 1$ for any positive integer $y$. \\
Identity, $F(w) = F(0)F(1)F(2)...F(z)...F(w-1) + 2$, holds\footnote{can be shown by induction.}. \\
and therefore, $F(w) \equiv 2 \; (\bmod \; F(z) )$. \\
\end{enumerate}

\noindent For $\mu = -1$, $2^{2^w} - 1, w > 2$, can be factorized\footnote{Identity in \S \ref{2_2_w_plus_mu_pn_qr}, item (\ref{Fermat_identity})} as $(2+1)(2^2+1)(2^4+1)...(2^{2^{w-1}}+1)$. Thus, $2^{2^w} - 1 \neq p^nq^r$ for $w > 2$ since at least three distinct primes (3, 5, and 17) divide $2^{2^w}-1$. For $w = 2, 2^{2^w} - 1 = 15 = 3 \cdot 5$ so that $n = r = 1$. For $w = 1, 2^{2^w} - 1 = 3 \neq p^nq^r$.  $\blacksquare$

\noindent \textcolor{black}{ \rule{16cm}{1pt}} 
\section{$p^n + \mu q^r = 2^m$}
\label{pn_mu_qr_2m}

\noindent \textcolor{black}{Solutions to $p^n + \mu q^r = 2^m$, for $\mu = -1$, correspond to the Hugh Edgar's problem \cite{Guy_1981}. Scott and Styer \cite{Scott_Styer_2004} address a more general diophantine form $p^n + \mu q^r + \beta 2^m = 0$, $m, n, r$ positive integers, $p$ and $q$ distinct odd primes, and $\mu, \beta = \pm 1$. Using results in \cite{Scott_1993, Stroeker_and_Tijdeman_1982, Cao_1986}, they demonstrate (\cite{Scott_Styer_2004}, Theorem 8) that such a relation has at most two solutions for ($m,n,r, \mu, \beta$) if the prime set $(p, q)$ is not a permutation of $(5,3)$ [7 solutions], $(7,3)$ [4 solutions], $(11,3)$ [3 solutions] or $(13,3)$ [3 solutions]. Scott and Styer \cite{Scott_Styer_2004} assume $(p,q)$  not to be a permutation of $(11,5)$}. \\ 

\noindent Specific solutions  $(p, q, m,n,r, \mu)$ are obtained next for the diophantine form $p^n + \mu q^r = 2^m$ \textcolor{black}{with $p$ and $q$ as Mersenne and/or Fermat primes}.

\subsection{$gcd(n,r) = 2^w, w \geq 0$, must hold}
\label{a_n_r}
\noindent Let $n = au$ and $r = av$, $a$ an odd prime, and $u,v \geq 1$. \\

$p^{ua} + (\mu q^{v})^{a} = (p^u + \mu q^v)R = (p^u + \mu q^v)\left[ p^{u(a-1)} -  p^{u(a-2)} \mu q^{v} +  ... -  p^{u} \mu q^{v(a-2)} + q^{v(a-1)}     \right]$. \\

\noindent $R$ is odd since it comprises odd number of terms, all odd. $R \nmid 2^m$, thus, $gcd(n,r) = 2^w, w \geq 0$. That is, no odd prime divides $n$ and $r$.\\


\subsection{$n$ and $r$ even}
\label{n_even_r_even}
\noindent Consider $n$ and $r$ both, even. Let $p^{\frac{n}{2}} = 2G + 1$ and $q^{\frac{r}{2}} = 2H + 1$ where $G$ and $H$ are positive integers. 

\begin{enumerate}
\item For $\mu = 1$, $p^n + q^r = 4G(G+1) + 4H(H+1) + 2 = 2^m$, or, $2^{m-1} = 2G(G+1) + 2H(H+1) + 1$, impossible for $m > 1$ since the $LHS$ is even and $RHS$ odd. 

\item For $\mu = -1$, $p^n - q^r = 4G(G+1) - 4H(H+1) = 2^m$, or, $\frac{G(G+1)}{2} - \frac{H(H+1)}{2} = 2^{m-3}$, which, when rewritten, becomes $(1 + 2 + 3 + ... + G) - (1 + 2 + 3 + ... + H) = 2^{m-3}$.

\begin{enumerate} 
\item \label{Maa_1} $G > H$ must hold since $2^{m-3} > 0$. Let $G = H + T$ for $T$ a positive integer. Then, 
\begin{eqnarray}
\sum \limits^G_{i = 1} i - \sum \limits^H_{j = 1} j = \sum \limits^T_{k = 1} \left( H +  k \right) = TH + T \left( \frac{T+1}{2} \right) = 2^{m-3} \nonumber 
\end{eqnarray}

\item Consider $m > 3$. $T$ is not divisible by an odd prime since $T \mid 2^{m-3}$. \\
Let $T = 2^{a}, a \geq 0$ and $2H + T + 1 = 2^{b}$ so that $T(2H + T + 1) = 2^{a+b} = 2^{m-2}$, or, $a + b = m-2$.
\item As $2H = 2^b - 2^a - 1 > 0$, $2^{m-2-2a} > 1 + \frac{1}{2^a} > 2^{0}$, or, $b = m - 2 - a > a$ must hold.
\item If $a > 0$, $2^a + 1 \geq 3$ is odd. Since $b = m-2-a > a$, $2H = 2^b - 2^a - 1$ is an impossibility, since the $LHS$ is even, and $RHS$ odd.  
\item \label{BeeBee} $a = 0$ must hold making $T = 1$ and $H = 2^{b-1}-1$. Since $b = m-2-a = m-2$, $H = 2^{m-3}-1$ and, from item (\ref{Maa_1}), $G = H+1 = 2^{m-3}$. 
\item \label{concl3} Thus, $p^{\frac{n}{2}} = 2G + 1 = 2^{m-2} + 1$ and $q^{\frac{r}{2}} = 2H + 1 = 2^{m-2}-1$. As $2^{m-2}-1$ is not a power per Appendix \S \ref{2_k_pm_mu}, $r = 2$ must hold. For $p^{\frac{n}{2}} = 2^{m-2} + 1$,  possibilities (Appendix \S \ref{2_k_pm_mu}) are --- \\
(P1): $n = 4$ when $p = 3$ and $m = 5$ so that $q = 7$, and \\
(P2): $n = 2$.
\item \label{concl3pt5} In case (P1), $p^n - q^r = 2^m$ becomes $\boxed{3^4 - 7^2 = 2^5}$.
\item \label{concl4} In case (P2), for $p$ and $q$  to be prime, $m-2$ is either 2, or must simultaneously be a power of 2 and an odd prime respectively, the latter an impossibility for $m \geq 5$. For $m-2 = 2$, $p = 2^{m-2} + 1 = 5$ and $q = 2^{m-2} - 1 = 3$. $p^n - q^r = \boxed{5^2 - 3^2 = 2^4}$.
\item $m \leq 3$ is an impossibility since then, $H \leq 0$ per item (\ref{BeeBee}), not permitted.
\end{enumerate}
\end{enumerate}

\subsection{$p = 2^g + \gamma, q = 2^d + \beta$}
\label{p_q_in mers_ferm_form}
Possibilities for Mersenne/Fermat primes $p, q$ and positive integer exponents $m,n$ and $r$ are explored next such that for $\mu = \pm 1$
\begin{eqnarray}
p^n + \mu q^r = 2^m
\label{main_dio3}
\end{eqnarray}
\noindent holds. With $p= 2^g + \gamma$ ($g$, a positive integer, $\gamma = \pm 1$) and $q =2^d + \beta$ ($d$, a positive integer, $\beta = \pm 1$), Eq. (\ref{main_dio3}) becomes
\begin{eqnarray}
2^gR_1 + 2^d \mu R_2 + \gamma^n + \mu \beta^r = 2^m
\label{dio_3_eq1}
\end{eqnarray}
\noindent where
\begin{eqnarray}
R_1 = (2^g)^{n-1} + \comb[n]{1}(2^g)^{n-2}\gamma +  \comb[n]{2}(2^g)^{n-3} + ... + \comb[n]{n-2}(2^g)\gamma^{n-2} +  n\gamma^{n-1} \nonumber \\
R_2 = (2^d)^{r-1} + \comb[r]{1}(2^d)^{r-2}\beta +  \comb[r]{2}(2^d)^{r-3} + ... + \comb[r]{r-2}(2^d)\beta^{r-2} +  r\beta^{r-1}
\label{dio_3_eq2}
\end{eqnarray}

\noindent Cases with $n$ and $r$ both even are addressed in \S \ref{n_even_r_even}. Those wherein at least one of the $n$ or $r$ is odd, are considered next. \\

\subsection{$p^n + q^r = 2^m$ $(\mu = 1)$}
\label{main_dio3_mu_1}
\noindent For $\mu = 1$, $m \geq 3$ must hold since minimum values that $p$ and $q$ can admit are $3$ and $5$, and those for $n$ and $r$ are 1. Without loss of generality, let $n$ be odd which makes $R_1$ in Eqs. (\ref{dio_3_eq1}) and (\ref{dio_3_eq2}) odd. 

\subsubsection{$p^n + q^r = 2^m$, $n$ and $r$ odd}
\label{main_dio3_mu_1_r_odd}

\noindent Per  Eq. (\ref{dio_3_eq2}), $R_1$ and $R_2$ are both odd. Eq.  (\ref{dio_3_eq1}) becomes $2^gR_1 + 2^d R_2 + \gamma + \beta = 2^m$.

\begin{enumerate}
\item \label{Maa_2} If $\gamma = \beta$, Eq.  (\ref{dio_3_eq1}) takes the form $2^{g-1}R_1 + 2^{d-1} R_2 = 2^{m-1} - \beta$. \\
$RHS$ is odd for $m > 1$. For $LHS$ to be odd, either \\
(P1): $g = 1$ and $d > 1$ so that $p = 2+\gamma$ and $q = 2^d + \gamma$, or \\
(P2): $g > 1$ and $d = 1$ implying $p = 2^g + \gamma$ and $q = 2 + \gamma$. \\
Both cases are identical. Considering (P1), $\gamma$ must be 1 since $p$ and $q$ are odd primes making $p = 3$ and $q = 2^d + 1$. $d$ cannot be 1 ($q$ is distinct from $p$). $d$ must be $2^w, w \geq 1$. 
Eq. (\ref{main_dio3}), that is, $p^n + q^r = 2^m$ becomes  $3^n + (2^d + 1)^r = 2^m$. 

\begin{enumerate}
\item \label{Maa_3} $2^m - (2^d + 1)^r \equiv (-1)^m - [(-1)^d + 1]^r = (-1)^m - 2^r \equiv (-1)^m - (-1)^r$ $= (-1)^m + 1 \;( \bmod \; 3)$. $m$ must be odd for 3 to divide $2^m - (2^d + 1)^r$.
\item Rewriting, $3^n + (2^d + 1)^r =  (2^2 - 1)^n + (2^d + 1)^r = 2^m$, which when expanded and divided by 4, yields  \\ \\
$\frac{2^{2n} - \comb[n]{1}2^{2n-2} + ... - \comb[n]{2}2^{4} + (2^d)^r + \comb[r]{1}(2^d)^{r-1} + ... + \comb[r]{2}(2^d)^{2}}{4} + n + r2^{d-2}= 2^{m-2}$.
\item If $d > 2$, $n + r2^{d-2}$ is odd (thus $LHS$ is odd) while $RHS$ is even. 
\item \label{Maa_6} As $d$ cannot be 1 and $d$ is a power of $2$ per item (\ref{Maa_2}), $d$ must be 2. $q = 5$ and thus  $p^n + q^r = 2^m$ becomes $3^n + 5^r = 2^m$ with all exponents odd.

\item \label{Maa_E2} \label{concl5}  \textcolor{black}{Per Scott and Styer \cite{Scott_Styer_2004}, primes $(2, 3, 5)$ combine in the arrangement $p^n \pm q^r = 2^m$ in precisely seven possible ways: $\boxed{3^3 - 5^2 = 2}$ (\S \ref{main_dio3_mu_minus_1}, item \ref{Maa_17}, see footnote), $\boxed{3^2 - 5 = 2^2}$ (\S \ref{main_dio3_mu_minus_1}, item \ref{Maa_B}), $\boxed{3 + 5 = 2^3}$, $\boxed{3^3 + 5 = 2^5}$, $\boxed{3 + 5^3 = 2^7}$, $\boxed{5 - 3 = 2}$ (\S \ref{main_dio3_mu_minus_1}, item \ref{Maa_18}, see footnote) and $\boxed{5^2 - 3^2 = 2^4}$ (\S \ref{n_even_r_even}). $(m, n, r)$ has no other solution for $|3^n \pm 5^r| = 2^m$}.\\

\end{enumerate}

\item If $\gamma = -\beta$, \\  
Eq. (\ref{dio_3_eq1}) yields $2^{g}R_1 + 2^{d} R_2 = 2^{m}$. As $R_1$ and $R_2$ are both positive and odd, $m > g = d$ must 
hold\footnote{\label{foot_1n} If $g > m$ or $d > m$, $2^{g-m}R_1 + 2^{d-m} R_2 = 1$ is infeasible since $LHS > 1$. $m$ must be greater than $g$ and $d$. Further, if $g \neq d$, both sides of $2^{g}R_1 + 2^{d} R_2 = 2^{m}$ will have opposite parities upon division with $2^{min(g,d)}$.} in  
 which case $p = 2^d - \beta$ and $q = 2^d + \beta$. $p$ and $q$ being twin primes,  $d$ must 
be 2,  and for $\beta = \pm 1$, one of the primes is 3 and the other 5. Let $p = 3$ and $q = 5$. 
Eq. (\ref{main_dio3}), as, $3^n + 5^r = 2^m$ is addressed in \S \ref{main_dio3_mu_1_r_odd}, items (\ref{Maa_6}) and (\ref{Maa_E2}).
\end{enumerate}

\subsubsection{$p^n + q^r = 2^m$, $r$ even, $n$ odd}
\label{main_dio3_mu_1_r_even}

\noindent  $R_2$ in Eq. (\ref{dio_3_eq2}) is even. Eq. (\ref{dio_3_eq1}) yields $2^gR_1 + 2^d R_2 + \gamma + 1 = 2^m$.

\begin{enumerate}
\item If $\gamma = 1$, \\
$2^{g-1}R_1 + 2^{d-1} R_2 = 2^{m-1} - 1$. $RHS$ is odd, since $m \geq 3$ must hold.  $2^{d-1} R_2$ is even for $d \geq 1$. $g$ must be 1 for $2^{g-1}R_1$ to be odd making $p = 3$ and $q = 2^d + \beta$ so that Eq. (\ref{main_dio3}) takes the form $3^n + (2^d + \beta)^r = 2^m$.

\begin{enumerate}
\item \label{Maa_8} For $\beta = 1$, since $d$ must\footnote{$d$ must be greater than 1, otherwise, $q = 3$, not permitted as $p = 3$.} be a power of 2 for $2^d + 1$ to be prime, $2^m - (2^d + 1)^r \equiv (-1)^m - [(-1)^{2^w} + 1]^r \equiv (-1)^m - 1 \;( \bmod \;  3)$. 
\item \label{Maa_9} For $\beta = -1$, as $d$ must\footnote{$d$ cannot be 2 as then, $q = 3$, not permitted as $q \neq p$.} be an odd prime $v$ for $2^d - 1$ to be prime, $2^m - (2^d - 1)^r \equiv (-1)^m - [(-1)^{v} - 1]^r \equiv (-1)^m - 1 \;( \bmod \; 3)$. 
\item In both cases, $m$ must be even so that $3 \mid 2^m - (2^d + \beta)^r$. 

\item \label{Maa_10} $3^n + (2^d + \beta)^r = 2^m$, $r$ and $m$ even, when rearranged, yields \\
$3^n = [2^\frac{m}{2} - (2^d + \beta)^\frac{r}{2}][2^\frac{m}{2} + (2^d + \beta)^\frac{r}{2}]$. This expression is not possible if \\  $2^\frac{m}{2} - (2^d + \beta)^\frac{r}{2} > 1$ as $gcd[2^\frac{m}{2} - (2^d + \beta)^\frac{r}{2}, 2^\frac{m}{2} + (2^d + \beta)^\frac{r}{2}] = gcd[2^\frac{m}{2} - (2^d + \beta)^\frac{r}{2},2\cdot 2^{\frac{m}{2}}] = 1$.
\item \label{Maa_7} If $2^\frac{m}{2} - (2^d + \beta)^\frac{r}{2} = 1$, $r$ must be $2$ since $2^\frac{m}{2} - 1$ is not a perfect power, per Appendix \S \ref{2_k_pm_mu}. Further, 

\begin{enumerate}
\item If $\beta = 1$ in item (\ref{Maa_7}), $2^{d-1} = 2^{\frac{m}{2}-1}-1$, not possible for $d  > 1$ (item \ref{Maa_8}, footnote). 
\item If $\beta = -1$, for $q = 2^d - 1 = 2^\frac{m}{2} - 1$ to hold, $d= \frac{m}{2}$ must be an odd prime $v$, per item (\ref{Maa_9}), so that $m = 2v$. $3^n + (2^d + \beta)^r = 2^m$ in item (\ref{Maa_10}) becomes $3^n + (2^v - 1)^2 = 2^{2v}$, or, $3^n = 2^{v+1} - 1 = (2^{\frac{v+1}{2}} - 1)(2^{\frac{v+1}{2}} + 1)$, which again, is not  possible\footnote{$gcd(2^{\frac{v+1}{2}} - 1,2^{\frac{v+1}{2}} + 1) = 1$.} if $2^{\frac{v+1}{2}} - 1 > 1$ since $RHS$ is divisible by two or more odd primes. If $2^{\frac{v+1}{2}} - 1 = 1$, $v$ must be 1, a contradiction,
since $v$ is an odd prime. \\
\end{enumerate}
\end{enumerate}

\item \label{Maa_bobby153} If $\gamma = -1$, \\
Eq. (\ref{dio_3_eq1}) gives $2^{g}R_1 + 2^{d} R_2 = 2^{m}$. \\
Since $R_1, R_2$ are positive integers, $m > d$ and $m > g$ must be true. Further, since $R_2$ is even, $g > d$ must hold to maintain parity (footnote \ref{foot_1n}, page \pageref{foot_1n}). Then, with $p = 2^{g} - 1$ and $q = 2^d + \beta$, Eq. (\ref{main_dio3})  becomes $(2^{g} - 1)^n + ( 2^{d} + \beta)^r = 2^m$ with $m > g > d$.

\begin{enumerate}
\item \label{d_plus_k_odd} $g$ must be an odd prime. \\
Since $p = 2^{g} - 1$ is prime, $g$ is either 2 or an odd prime $v$. If $g = 2$, $p = 3$; \\
Since $d < g$, $d$ and thus $\beta$ must be 1 so that $q = 2^{1} + 1 = 3$, not permitted, as $p \neq q$. \\

\item \label{concl6} Consider $m$ even. \\
$(2^{g} - 1)^n = 2^m - ( 2^{d} + \beta)^r = [2^\frac{m}{2} - ( 2^{d} + \beta)^\frac{r}{2}][2^\frac{m}{2} + ( 2^{d} + \beta)^\frac{r}{2}]$. \\
$gcd(2^\frac{m}{2} - ( 2^{d} + \beta)^\frac{r}{2},2^\frac{m}{2} + ( 2^{d} + \beta)^\frac{r}{2}) = 1$. \\
$(2^{g} - 1)^n \neq 2^m - ( 2^{d} + \beta)^r $ if $2^\frac{m}{2} - ( 2^{d} + \beta)^\frac{r}{2} > 1$. \\
If $2^\frac{m}{2} - ( 2^{d} + \beta)^\frac{r}{2} = 1$, or, $2^\frac{m}{2} - 1 = (2^{d} + \beta)^\frac{r}{2}$, $r$ must be 2, per Appendix \S \ref{2_k_pm_mu}. \\ \\
If $\beta$  $= -1$, $d = \frac{m}{2}$. Then, $(2^{g} - 1)^n = 2^{\frac{m}{2}} + (2^d + \beta)^{\frac{r}{2}} = 2 \cdot 2^{\frac{m}{2}} - 1 = 2^{\frac{m}{2}+1} - 1$. $n$ must be 1, per Appendix \S \ref{2_k_pm_mu}. $p=2^{g} - 1$ becomes $2^{\frac{m}{2}+1} -1$ while $q$ is $2^\frac{m}{2} - 1$. \\
Since both are prime, $\frac{m}{2}$ and ${\frac{m}{2}+1}$ must either be 2, or odd prime, possible only 
when $m = 4$, making $p = 7$, $q = 3$ so that $\boxed{7 + 3^2 = 2^4}$. $g = 3, d = 2$ satisfy $m > g > d$. \\

If $\beta = 1$, $2^{\frac{m}{2}-1} - 1 = 2^{d-1}$ holds for $d = 1$, $m = 4$ making $(2^{g} - 1)^n = 2^m - (2^{d} + \beta)^r$ as $\boxed{7 = 2^4 - 3^2}$, same as the case with $\beta = -1$. \\

\item Consider $m$ odd. \\
Some examples for $x^2 = 2^m -  (2^{g} - 1)^n$, $m,n$ odd, $x$ a positive integer not necessarily prime, include (E1): $5^2 = 2^5 - 7$, (E2): $11^2 = 2^7 - 7$, (E3): $13^2 = 2^9 - 7^3$, and (E4): $181^2 = 2^{15} - 7$. In these, $gcd(n,r) = 1$. Only in (E1) is $x = q = 5$ of the form  $2^d + \beta, \beta = \pm 1$, $d \geq 1$.

\begin{enumerate}
\item $q = 2^d + \beta < p = 2^{g} - 1$ \textcolor{black}{must hold}. \\
In $(2^{g} - 1)^n + ( 2^{d} + \beta)^r = 2^m$, $m$ and $n$ odd and $r$ even, as $g > d$ (\S \ref{main_dio3_mu_1_r_even}, item \ref{Maa_bobby153}), $ 2^{g} - 1 \geq  2\cdot 2^{d} - 1 \geq 2^d + \beta$, or $2^{d} \geq 1 + \beta$ holds for $d \geq 1$ and $\beta = \pm 1$. Equality holds for $g = 2$, $d = \beta = 1$ for which $p = q$, not permitted. $p$ must be greater than $q$. \\

\item $q = 2^d + \beta$ cannot be 3. \\
Since $g$ is an odd prime (item \ref{d_plus_k_odd}), $2^m - (2^{g} - 1)^n \equiv (-1)^m - [(-1)^{g} - 1]^n = (-1) - (-2)^n = 2^n -1 $ $ \equiv (-1)^n -1 = -2 \mbox{ or $1$} \;( \bmod \; 3)$.\\

\item \label{Maa_bobby153aa_guw2020} 
\textcolor{black}{For $n = 1$, using the Theorem of Szalay \cite{szalay}, solutions to $(2^{d} + \beta)^r = x^2 = 2^m - (2^{g} - 1)^n$, $m$ any positive integer, are given by (\cite{Scott_2018}, Theorem 1.2)\\
$\bullet$ $(m, g, x) = (2t, t+1, 2^t-1)$, $t \geq 1$ an integer, not possible, if $m$ is odd. \\
$\bullet$ $(m, g, x) = (5, 3, 5)$, depicted in item (\ref{Maa_bobby153a}). \\
$\bullet$ $(m, g, x) = (7, 3, 11)$, not possible since $x = 11$ is not of the form $2^d \pm 1$. \\
$\bullet$ $(m, g, x) = (15, 3, 181)$, not possible as $x = 181$ is not of the form $2^d \pm 1$. \\
}

\item  \label{Maa_bobby153aa} \label{concl7a} \textcolor{black}{The diophantine equation $x^2 - 2^m = \pm y^n$, $x, y, m, n$, as positive integers, and $x$ and $y$ coprime, $y > 1$, $n \geq 2$, is considered by Bugeaud in \cite{Bugeaud_1997}. Basing his work on \cite{Laurent_Mignotte_Nesterenko_1995, Bugeaud_Laurent_1996}, Bugeaud shows that for $(x, y, m, n)$ to be a solution to (e1): $x^2 - 2^m = y^n$ or (e2): $x^2 - 2^m = - y^n$, $m$ and $n$ must be odd, and further, $n \leq 5 \cdot 5 \; 10^5$ for (e1) and $n \leq 7 \cdot 3 \; 10^5$ for (e2) must hold. Specifically, per \cite{Bugeaud_1997}, $x^2 = 2^m \pm y^n$ (and thus $(2^{d} + \beta)^r = 2^m - (2^{g} - 1)^n$, $r$ even, $m$ and $n$ odd) have only finitely many solutions}. \\

\item \label{Maa_bobby153a} \label{concl7} \textcolor{black}{Scott and Styer \cite{Scott_Styer_pc} contribute to the solution of $p^n + q^r = 2^m$, $p = 2^g - 1$ (prime) and $q$ any odd prime, $m, n$ odd and $r$ even, as follows: \\
In Theorem 4.1 of \cite{Scott_2018}, set $C = p^n, P = p, Q = 1$ and $m(-P) = m(-p) = g-2$. $2^{g-2}$ is the lowest power of 2 that can be factored into relatively prime principal ideals in $\mathbb{Q}(\sqrt{-p}): [2^{g-2}] = \left[\frac{1+\sqrt{-p}}{2} \,\, \frac{1 - \sqrt{-p}}{2} \right]$. $m > g$ must hold. The bound in Theorem 4.1 of \cite{Scott_2018} can be simplified to $m(-P)$ since $m$ is odd, $p \equiv 7 \;( \bmod \; 8)$ and $Q = 1$ to eliminate factors of 2 and $3^u$, and the symbol for the least common multiple. As $m-2 > g-2 = m(-P)$, any solution to $p^n + q^r = 2^m$ must belong to the exceptional cases of Theorem 4.1 \cite{Scott_2018}. Pertaining to $(2^d  + \beta)^r = 2^m -  (2^{g} - 1)^n$, the only exceptional case is $\boxed{5^2 = 2^5 - 7}$}. 

\end{enumerate}
\end{enumerate}
\end{enumerate}

\subsection{$p^n - q^r = 2^m$ $(\mu = -1)$}
\label{main_dio3_mu_minus_1}

\noindent \textcolor{black}{Given $p$ and $q$ as distinct odd primes, Scott and Styer  (\cite{Scott_Styer_2004}, Theorem 6) show that there can be at most one solution in $(n,r)$ for $p^n - q^r = 2^m$ to hold}.  With $p$ and $q$ as Mersenne or Fermat primes, following cases are considered: (i) $r$ and $n$ odd, (ii) $r$ odd, $n$ even, and (iii) $r$ even and $n$ odd. Those with $n$ and $r$ even are addressed in \S \ref{n_even_r_even}.  


\subsubsection{$p^n - q^r = 2^m$, $r$ and $n$ odd.}
\label{subsubs1_19_11_25}

In Eq. (\ref{dio_3_eq2}), $R_1$ and $R_2$ are both odd. Eq. (\ref{dio_3_eq1}) becomes $2^{g}R_1 - 2^{d} R_2 + \gamma - \beta = 2^{m}$. 

\begin{enumerate}
\item Consider $m > 1$.

\begin{enumerate}
\item \label{Maa_13} If $\gamma = -\beta$, \\
 $2^{g-1}R_1 - 2^{d-1} R_2 = 2^{m-1} + \beta$ (odd). So that $LHS$ is odd, either  \\ 
 (P1): $g = 1$ and $d > 1$, or \\
 (P2): $g > 1$ and $d = 1$.

\begin{enumerate}
\item  \label{Maa_14} In case of (P1), $p = 2 - \beta$. $\beta$ must be $-1$ as $p$ is an odd prime making $p = 3$.
$q = 2^d + \beta = 2^d - 1$. $d$ cannot be 2 since $q \neq p$. $d$ must be an odd prime. 
Eq. \ref{main_dio3} becomes  $3^n - (2^d - 1)^r = (2^2 - 1)^n - (2^d - 1)^r = 2^m$. \\ 
Expansion and division by 4 results in \\ \\
$\frac{2^{2n} - \comb[n]{1}2^{2n-2} + ... - \comb[n]{n-2}2^{4} - (2^d)^r + \comb[r]{1}(2^d)^{r-1} - ... + \comb[r]{2}(2^d)^{2}}{4} + n - r2^{d-2}= 2^{m-2}$. \\ \\
Since $d$ is an odd prime, $n - r2^{d-2}$ is odd making $LHS$ odd, but $RHS$ is even for 
$m \geq 3$, not possible. \\ \\

If $m = 2$, Eq. (\ref{main_dio3}) as $3^n - (2^d - 1)^r = 4$ with $n$ and $r$ odd is  impossible as $RHS$ is $\equiv 1 \;( \bmod \; 3)$, 
and $3^n - (2^d - 1)^r \equiv -[(-1)^d - 1]^r = - (-2)^r = 2^r \equiv (-1)^r = -1 \;( \bmod \; 3)$.\\

\item In  case of (P2), for $m > 2$, as in item (\ref{Maa_14}), $p = 2^g - 1$ and $q = 2 + 1 = 3$. Eq. (\ref{main_dio3}) becomes $(2^g-1)^n - 3^r = (2^g-1)^n - (2^2-1)^r = 2^m$. $g$ cannot be 2 ($p \neq q$) and  
thus must be  an odd prime. Expansion and division by 4 yields \\ \\
$\frac{(2^g)^n - \comb[n]{1}(2^g)^{n-1} + ... - \comb[n]{2}(2^g)^{2} - 2^{2r} + \comb[r]{1}2^{2r-2} - ... +  \comb[r]{2}2^{4}}{4} + n2^{g-2} - r = 2^{m-2}$. \\ \\
$g$ being an odd prime, $n2^{g-2} - r $ is odd (thus $LHS$ is odd), but $RHS$ is even, 
not possible. \\ \\

\label{Maa_15} \label{concl8} If $m = 2$, Eq. (\ref{main_dio3}) becomes $(2^g-1)^n = 3^r + 4$.  \\
\textcolor{black}{Per \cite{Scott_Styer_2004}, lemma 2, $N^z = 3^x + 2^y$ has no solution in positive integers $(x,y,N,z)$ for $z > 1$, except $\boxed{5^2 = 3^2 + 2^4}$ (\S \ref{n_even_r_even}). Proof employs Theorem 2 of \cite{Scott_1993} that suggests that $z$ must be 1, unless $x$ and $y$ are both even, and $z = 2$}. \\ \\

If a solution exists for $(2^g-1)^n = 3^r + 4$, $n$ and $r$ odd, $n$ must be 1. $2^g = 3^r + 5$ yields only two solutions (\S \ref{main_dio3_mu_1_r_odd}, item \ref{Maa_E2}), $(g, r) = (3, 1)$ and $(g, r) = (5, 3)$ for which $(2^g-1) = 3^r + 4$ becomes \label{concl8pt5}  $\boxed{7 = 3 + 2^2}$ and $\boxed{31 = 3^3 + 2^2}$ respectively. \\

\end{enumerate}

\item \label{Maa_bobby9} If $\gamma = \beta$, \\
Eq. (\ref{dio_3_eq1}) simplifies to $2^gR_1 - 2^d R_2 = 2^m$. Consider $m \geq 2$. As $R_1$ and $R_2$ are both odd, three possibilities arise --- \\
(P1): $m > g = d$ so that $R_1 - R_2 = 2^{m-d}$: $p = q = 2^d + \beta$, not permitted. \\

(P2): $g > m = d$ so that $2^{g-m}R_1 = R_2 + 1$: \\
Eq. (\ref{main_dio3}) becomes $(2^g + \gamma)^n - (2^d + \gamma)^r = 2^d$. \\
$\bullet$ For $\gamma = 1$, $g$ and $d$ are both powers of 2. \\
$2^d + (2^d + 1)^r \equiv (-1)^d + [(-1)^d + 1]^r = 1 + 2^r \equiv 1 + (-1)^r = 0 \;( \bmod \; 3)$ implying $2^g + \gamma = 3$, or, $g = 1$, not possible since $g > d \geq 1$. \\ \\
$\bullet$ For $\gamma = -1$, $g$ and $d$ are either 2 or odd primes. Since $g > d$, $g$ is an odd prime. In case $d = 2$, Eq. (\ref{main_dio3}) as $(2^g -1)^n  = 3^r +4$ is addressed in \S \ref{subsubs1_19_11_25} item (\ref{Maa_15}) above. In case $d$ is an odd prime, Eq. (\ref{main_dio3}) as $(2^g -1)^n - (2^d -1)^r = 2^d$ does not hold since $3$ divides the $LHS$\footnote{\label{Maa_f11} $(2^g -1)^n - (2^d -1)^r \equiv [(-1)^g - 1]^n - [(-1)^d - 1]^r = -2^n + 2^r \equiv 1^n - 1^r = 0 \;( \bmod \;  3)$.}. \\

(P3): $d > m = g$ so that $R_1 - 1 = 2^{d-m}R_2$:\\
Eq. (\ref{main_dio3}) becomes $(2^g + \gamma)^n - (2^d + \gamma)^r = 2^g$. \\
$\bullet$ For $\gamma = 1$, $g$ and $d$ are both powers of 2.  $d$ must be greater than 1 since $d > g \geq 1$. \\
If $g > 1$, $2^g + (2^d + 1)^r \equiv (-1)^g + [(-1)^d + 1]^r = 1 + 2^r \equiv 1 + (-1)^r = 0 \;( \bmod \; 3)$ implying $2^g + 1 = 3$, or, $g = 1$, a contradiction.\\
If $g = 1$, Eq. (\ref{main_dio3}) as $3^n - 2 = (2^d + 1)^r$, or, $3^n - 3 = (2^d + 1)^r - 1$ does not hold since $3 \nmid$ $RHS$\footnote{$(2^d + 1)^r - 1 \equiv [(-1)^d + 1]^r - 1 \equiv 2^r - 1 \equiv (-1)^r - 1 \equiv -2 \;( \bmod \; 3)$.}. \\ \\
$\bullet$ For $\gamma = -1$, since $d > g$, $d$ must be an odd prime while $g$ is either 2 or an odd prime. \\
If $g = 2$, Eq. (\ref{main_dio3}) becomes $3^n = 2^2 + (2^d -1)^r$. $3 \nmid RHS$ since $2^2 + (2^d -1)^r \equiv (-1)^2 + [(-1)^d - 1]^r \equiv 1 - 2^r \equiv 1 - (-1)^r = 2 \;( \bmod \; 3)$. \\
If $g$ is an odd prime, $(2^g -1)^n - (2^d -1)^r = 2^g$ does not hold (footnote \ref{Maa_f11}). \\ \\
\end{enumerate}

\item Consider $m = 1$.
\begin{enumerate}
\item For $\gamma = -\beta$, \\
Eq. (\ref{dio_3_eq1}) becomes $2^{g-1}R_1 - 2^{d-1} R_2 = 1 + \beta$. $R_1$ and $R_2$ are both odd. \\

\begin{enumerate}
\item \label{Maa_17} For $\beta = 1$,  \\
Eq. (\ref{dio_3_eq1}) simplifies to $2^{g-2}R_1 - 2^{d-2} R_2 = 1$ (odd). Either $2^{g-2}R_1$ or $2^{d-2} R_2$, but not both, must be odd. Two possibilities arise: \\ 
(P1): $g > 2, d = 2$,  or \\
(P2): $g = 2, d > 2$. \\ \\ 
In (P1), $q = 2^2 + 1$ and $p = 2^g - 1$. Eq. (\ref{main_dio3}) takes the form $(2^g-1)^n - 5^r = 2$. Since $g > 2$, it must be an odd prime as $2^g - 1$ is prime. \\

 \label{concl9} \textcolor{black}{Scott \cite{Scott_1993} states per Theorem 6 that if $a > 1, b > 1$ are relatively prime integers and if $p$ is prime, then the equation $a^x + b^y = p^z$ has at most two solutions\footnote{ \label{igmr24062024} Scott and Styer \cite{Scott_Styer_pc} contribute a more general solution for $N^z = 5^r + 2^s$ with integers $N, r, z, s > 0$ as follows: If $z > 1$, per Theorem 2 in \cite{Scott_1993}, $N^z$ is either square ($s$ must then be even) or cube ($r$ must be even while $s$ and $z$ must be odd). In the latter case, when $N^z$ is cube,  one writes $[N^{\frac{z}{3}}]^3 = 5^r + 2^s = \left[5^{\frac{r}{2}} + 2^{\frac{s-1}{2}}\sqrt{-2}\right]\left[5^{\frac{r}{2}} - 2^{\frac{s-1}{2}}\sqrt{-2}\right]$. As $\sqrt{-2}$ generates a unique factorization domain, there must exist an integer $u + v\sqrt{-2}$ such that $[u + v\sqrt{-2}]^3 = 5^{\frac{r}{2}} + 2^{\frac{s-1}{2}}\sqrt{-2}$. $u$ is either 1 or $5^x, x > 0$, and $v = 2^{\frac{s-1}{2}}$ (Lemmata 1-3 in\cite{Scott_1993}). $u = 5^x$ contradicts Theorem 1 of \cite{Scott_1993}. For $u = 1$, per Theorem 3.2 of \cite{Scott_2018}, only possibility for $N^z = 5^r + 2^s$ with $3$ dividing $z$ is $\boxed{5^2 + 2 = 3^3}$. Otherwise, $z$ must be $1$, and pertaining to $(2^g-1)^n - 5^r = 2$, $n$ must be $1$ so that the two solutions to $2^g = 5^r + 3$, $r$ odd and $g$ an odd prime, are $\boxed{2^3 = 5 + 3}$ and $\boxed{2^7 = 5^3 + 3}$. \\} in positive integers $(x, y, z)$ if $ p \neq 2$. There are exceptions for $(a, b)$ = $(3, 5)$ and $(3, 13)$. For $a = 2, b = 5$ and $p = 2^g - 1$, $x=1, y=r, z=n$ all odd,  $n$ must be $1$ so that the two solutions to $2^g = 5^r + 3$, $r$ odd and $g$ an odd prime, are (\S \ref{main_dio3_mu_1_r_odd}, item \ref{Maa_E2}) $\boxed{2^3 = 5 + 3}$ and $\boxed{2^7 = 5^3 + 3}$}. \\

In case of (P2) in item (\ref{Maa_17}), $p = 2^2 -1$ and $q = 2^d + 1, d = 2^w, w \geq 1$ since $q \neq 3$. Eq. (\ref{main_dio3}) becomes $3^n - (2^d + 1)^r = 2$. But $3 \nmid (2^d + 1)^r + 2$ since $(2^d + 1)^r + 2 \equiv [(-1)^d + 1]^r - 1 = 2^r - 1 = (-1)^r - 1 = -2 \;( \bmod \; 3)$. \\
  
\item \label{Maa_18} \label{Maa_F} \label{concl10} For $\beta = -1$,  \\ 
  Eq. (\ref{dio_3_eq1}) simplifies to $2^{g-1}R_1 - 2^{d-1} R_2 = 0$. $g = d$ and $R_1 = R_2$ must hold. 
  Then, $q = 2^d - 1$ and $p = 2^d + 1$, possible only when $d = 2$ which makes $q = 3$ 
  and $p = 5$. Eq. (\ref{main_dio3}) yields $5^n - 3^r = 2$. The only solution (\cite{Scott_Styer_2004}, Theorem 6), $\boxed{5 = 3+2}$, 
  exists for $n = r = 1$. \\
  
\end{enumerate}

\item For $\gamma = \beta$, \\
Eq. (\ref{dio_3_eq1}) yields $2^{g-1}R_1 - 2^{d-1} R_2 = 1$. As $R_1$ and $R_2$ are both odd, either \\
(P1): $2^{d-1}R_2$ alone must be odd ($d = 1, g > 1$), or \\
(P2): $2^{g-1} R_1$ alone must be odd ($g = 1, d > 1$). \\
 
\begin{enumerate}  
\item For $\beta = 1$, possibility (P1) yields $q = 3, p = 2^g + 1, g = 2^w, w \geq 1$, making Eq. (\ref{main_dio3}) as $(2^g + 1)^n - 3^r = 2$. Rewriting 3 as $2^2 - 1$ and expanding results in \\ \\
  $2^g [ (2^g)^{n-1}  + \comb[n]{1}(2^g)^{n-2} + ... + n] = 2^2 [(2^2)^{r-1}  - \comb[r]{1}(2^2)^{r-2} + ... + r]$, \\ \\
   or, $2^g H_1 = 2^2 H_2$. $H_1$ and $H_2$ are both odd since $n$ and $r$ respectively, are odd. 
  A solution to $(2^g + 1)^n - 3^r = 2$ is possible only when $g = 2$, making $2^g + 1 = 5$.
  The only solution as $5 = 3+2$ for $5^n - 3^r = 2$ is addressed in \S \ref{subsubs1_19_11_25} item (\ref{Maa_18}). \\ \\
  Possibility (P2) gives $p = 3, q = 2^d + 1, d = 2^w, w \geq 1$ so that \\
  $3^n - (2^d + 1)^r = 2$ which does not hold since $(2^d + 1)^r + 2 \equiv 2^r + 2 \equiv -2 \;( \bmod \; 3)$. \\
\item For $\beta = -1$, per  possibility (P1), $q = 2-1$, not permitted. Per  possibility (P2), $p = 2-1$,  not permitted.
\end{enumerate}
\end{enumerate}
\end{enumerate}


\subsubsection{$p^n - q^r = 2^m$, $r$ odd, $n$ even.}
\label{subsubs2_19_11_25}

\begin{enumerate}
\item Consider $m > 1$. \\
In Eq. (\ref{dio_3_eq2}), $R_2$ is odd and $R_1$ even. Eq. (\ref{dio_3_eq1}) becomes $2^gR_1 - 2^d R_2 + 1 - \beta = 2^m$.

\begin{enumerate}
\item If $\beta = -1$, the above becomes $2^{g-1}R_1 - 2^{d-1} R_2  = 2^{m-1} - 1$ (odd). $2^{g-1}R_1$ is even for $g \geq 1$. \\
 So that $2^{d-1} R_2$ is odd, $d$ must be 1. But, $q = 2 - 1$, not possible. \\
\item \label{BeeBee1} \label{concl11} If $\beta = 1$, then $2^{g}R_1 - 2^{d}R_2  = 2^m$. Two possibilities exist: \\ \\
(P1): If $g \geq m$ in $2^g R_1 - 2^m = 2^dR_2$, or in $2^{m}(2^{g-m}R_1 - 1) = 2^d R_2$, $m = d$ and $2^{g-m}R_1 - 1 = R_2$ must hold. In case $g = m$, $R_1 - 1 = R_2$ must still hold since $R_1$ is even. Eq. (\ref{main_dio3})  becomes $(2^{g} + \gamma)^n - (2^{d} + 1)^r = 2^d$. \\ \\
Case $d = m > 1$ ($d$ is a power of 2 since $2^d + 1$ is prime) is identical to that in \S \ref{subsubs1_19_11_25}, item (\ref{Maa_B}) covered later. $3 \mid 2^d + (2^{d} + 1)^r$ implying $2^{g} + \gamma$ must be $3$ making Eq. (\ref{main_dio3})  as $3^n - 2^d = (2^{d} + 1)^r$. As $n$ and $d$ are both even, $3^n - 2^d = [3^{\frac{n}{2}} - 2^{\frac{d}{2}}][3^{\frac{n}{2}} + 2^{\frac{d}{2}}] = (2^{d} + 1)^r$ is possible only when $3^{\frac{n}{2}} - 2^{\frac{d}{2}} = 1$, that is, when $n = 4$ and $d = 6$, or, when $n = 2$ (Appendix \S \ref{2_k_pm_mu}). In the former case, $d \neq 2^w$, while in the latter, $(2^{d}+1)^r = 2^{\frac{d}{2}+1} + 1$. As $r$ is odd, $r$ must be 1 (Appendix \S \ref{2_k_pm_mu}) and $d$ must be $2$ making $m = 2$ and thus\footnote{Since $g \geq m$ and $2^g + \gamma = 3$, $g$ must be 2.} Eq. (\ref{main_dio3})  as $(2^{g} + \gamma)^n = \boxed{5 + 4 = 3^2}$.\\ \\

(P2): If $g < m$ in $2^g R_1 - 2^m = 2^dR_2$, or, in $2^{g}(R_1 - 2^{m-g}) = 2^d R_2$, $d > g$ must hold since $R_1 - 2^{m-g}$ is even. \\

\begin{enumerate}
\item \label{Maa_B} \label{concl12} Consider $m$ even. \\
 $(2^{g} + \gamma)^n - 2^m =  [( 2^{g} + \gamma)^\frac{n}{2} - 2^\frac{m}{2}][( 2^{g} + \gamma)^\frac{n}{2} + 2^\frac{m}{2}]$. \\
 $( 2^{g} + \gamma)^\frac{n}{2} - 2^\frac{m}{2}$ and $( 2^{g} + \gamma)^\frac{n}{2} + 2^\frac{m}{2}$ are coprime. \\
 $(2^{d} + 1)^r \neq (2^{g} + \gamma)^n - 2^m$ if $(2^{g} + \gamma)^\frac{n}{2} - 2^\frac{m}{2} > 1$. \\
 If $(2^{g} + \gamma)^\frac{n}{2} - 2^\frac{m}{2}  = 1$, possibilities are --- (P2-Q1): $m = 6$ so that $n = 4$, or, (P2-Q2): $n = 2$ (Appendix \S \ref{2_k_pm_mu}). \\ \\
 (P2-Q1): If $n = 4$, $(2^{g} + \gamma)^2 = 3^2$, implying $p = 3$, that is, either $(g, \gamma) = (1,1)$ or $(g, \gamma) = (2,-1)$. $p^n - 2^m = \boxed{3^4 - 2^6 = 17} = (2^{d} + 1)^r$. 
 $d = 4$ and $r = 1$. \\ \\
 (P2-Q2): If $n = 2$,  $( 2^{d} + 1)^r = [( 2^{g} + \gamma)^\frac{n}{2} - 2^\frac{m}{2}][( 2^{g} + \gamma)^\frac{n}{2} + 2^\frac{m}{2}] = 2^{\frac{m}{2} + 1} + 1$. \\
 $2^{d}+1$ cannot be $3$ since $d > g \geq 1$ implying $m$ cannot be 4, per Appendix \S \ref{2_k_pm_mu}. Further,
 since $r$ is odd, $r$ must be 1.  $p$ becomes $2^{\frac{m}{2}} + 1$ while $q$ is $2^{\frac{m}{2} + 1} + 1$. 
 Since $p$ and $q$ are both prime, $m$ must be 2. Since $d = \frac{m}{2} + 1 > g$ and $m > g = \frac{m}{2}$, $g$ must be $1$ (thus $\gamma = 1$). $p^n - q^r = 2^m$ becomes \label{concl13}  $\boxed{3^2 - 5 = 2^2}$. \\
 
\item  \label{Maa_BbB} Consider $m$ odd. \\
Some examples pertaining to $x^2 = 2^m + (2^{d} + 1)^r$, $m, r$ odd, $x > 1$ (integer) include (E1): $(2^2+1)^2 = 2^3 + (2^4 + 1)^1$, (E2): $7^2 = 2^5 + (2^4 + 1)^1$, (E3): $(2^4+1)^2 = 2^5 + (2^8 + 1)^1$, (E4): $23^2 = 2^9 + (2^4 + 1)^1$, (E5): $(2^8+1)^2 = 2^9 + (2^{16} + 1)^1$ and (E6): $71^2 = 2^7 + (2^4 + 1)^3$.
 Of these, $x$ is of the form $2^g \pm 1$ (prime) in (E1), (E2), (E3) and (E5). Also, in all these examples, $d > g$ and $m > g$ hold. \\

In $(2^g + \gamma)^n = 2^m + (2^{d} + 1)^r$, $m$ and $r$ odd, $2^g + \gamma$ cannot be $3$. \\
$2^m + (2^{d} + 1)^r \equiv (-1)^m + [(-1)^{2^w} + 1]^r = -1 + 2^r \equiv -1 + (-1)^r = $
$-2 \mbox{ or } 1 \;( \bmod \; 3)$. Thus, $(g,\gamma) \neq (2,-1)$ or $(1,1)$.\\ \\

\item  \label{Maa_19} \label{concl14} \label{concl15} $(2^{g} + \gamma)^n = 2^m + (2^{d} + 1)^r$, $2^{g} \pm 1 \neq 3$,  $n$ even, $m$ odd and $r= 1$, has no solutions other than $\boxed{(2^{2^w} + 1)^2 = 2^{2^w + 1} + (2^{2^{w+1}} + 1)}$, $w \geq 1$ an integer, so that $2^{2^w} + 1$ and $2^{2^{w+1}} + 1$ are primes, and $\boxed{7^2 = 2^5 + (2^4 + 1)}$. \\

\textcolor{black}{Per Scott and Styer \cite{Scott_Styer_pc}, if $r = 1$, all solutions to $x^2 =  2^m + (2^d + 1)$ can be obtained through Theorem 1.1 (Theorem of Szalay) in \cite{Scott_2018}. In particular, solutions are $(d, m, x) = (2y, y+1, 2^y+1)$ for positive integer $y$, $(d, m, x) = (4, 5, 7)$ and $(9, 4, 23)$}. \\ \\
In the first case, as $q = 2^d + 1 = 2^{2y} + 1$ is prime, $y = 2^w$, $w \geq 1$, must hold. $w = 0$ makes $m$ even, not permitted. Since $x = p^{\frac{n}{2}}= 2^y + 1 = 2^{2^w} + 1$ is not a power since $y$ cannot be $3$ (Appendix \S \ref{2_k_pm_mu}), $n$ must be $2$. \\

For $(d, m, x) = (4, 5, 7)$,  $x^2 = 2^m + (2^d + 1)$ becomes $7^2 = 2^5 + (2^4 + 1)$. \\
In the third case, $x = 23$ is neither a Mersenne nor a Fermat prime. \\

\item \label{Maa_199999} \textcolor{black}{Siksek and Stoll \cite{Siksek_Stoll_2014} consider the form $x^2 = (-y)^3 + z^{15}$ and find the only non-trivial solution as $(x,y,z) \equiv (\pm3, -2, 1)$}. \\

\item \label{Maa_199999_igmr_29_5_24}\textcolor{black}{The equation $x^2 = 2^m + (2^{d} + 1)^r$ is addressed more generally by Bugeaud in \cite{Bugeaud_1997}. Also mentioned in \S \ref{main_dio3_mu_1}, item (\ref{Maa_bobby153aa}), for positive integer set $(x, y, m, n)$ to be a solution to $x^2 = 2^m  + y^n$, $m$ and $n$ must be odd, and further, $n \leq 5 \cdot 5 \; 10^5$ must hold. Specifically, per \cite{Bugeaud_1997}, $x^2 = 2^m \pm y^n$ (and thus $x^2 = 2^m + (2^{d} + 1)^n$, $m$ and $n$ odd) have only finitely many solutions.} 

\end{enumerate}
\end{enumerate}

\noindent \rule{14cm}{1pt}

\item Consider $m = 1$.
\begin{enumerate}
\item For $\beta = 1$, Eq. (\ref{dio_3_eq1}) becomes $2^{g-1}R_1 - 2^{d-1} R_2 = 1$. $2^{g-1}R_1 $ is even for $g \geq 1$. 
 $2^{d-1}R_2$ must be odd. As $R_2$ is odd, $d$ must be  $1$, implying $q = 3$ and $p = 2^g + \gamma$. 
 Eq. (\ref{main_dio3}) becomes $(2^g + \gamma)^n - 3^r = 2$. \\
\begin{enumerate}
\item For $\gamma = 1$, $g$ must be $2^w, w \geq 1$. $(2^g + 1)^n - 2 \equiv 2^n - 2 \equiv (-1)^n + 1 =$ 
  $ 2 \;( \bmod \;  3)$. $(2^g + 1)^n - 3^r = 2$ does not hold.
\item For $\gamma = -1$, since $p$ cannot be 3, $g$ must be an odd prime. $(2^g - 1)^n - 2 \equiv 2^n - 2 \equiv 2 \;( \bmod \;  3)$. 
 $(2^g - 1)^n - 3^r = 2$ does not hold.
\end{enumerate}

\item For $\beta = -1$, Eq. (\ref{dio_3_eq1}) takes the form $2^{g}R_1 - 2^{d} R_2 = 0$. Since $R_1$ is even, let 
 $R_1 = 2^k R_3$ such that $R_3$ is odd. Then, $d = g + k$, and $R_1 = R_3$ must hold. $p = 2^g + \gamma$ 
 and $q = 2^{g+k} - 1$ makes Eq. (\ref{main_dio3}) as $(2^g + \gamma)^n = 2 + (2^{g+k} - 1)^r$.  $g+k$ cannot be 2 
 (if so, then $g = k = 1$, and $\gamma = 1$ making $p = q = 3$, not permitted) and must be an 
 odd prime $v$. $2 + (2^{v} - 1)^r \equiv -1 + (-2)^r = -1 - (2)^r \equiv -1 - (-1)^r = 0 \;( \bmod \;  3)$
 implying $2^g + \gamma = p = 3$. Thus, Eq. (\ref{main_dio3}) becomes $(2^{v} - 1)^r + 2 = 3^n$ with $n$ even and $r$ odd. \\
 \begin{enumerate}
 	
\item {\bf Proposition}: In $(2^{v} - 1)^r + 2 = 3^n$, $v$ an odd prime, $r$ odd and $n$ even, $2 \mid n$ but $4 \nmid n$. \\
{\bf Proof}:  \\
Considering $ \;( \bmod \; 10)$, noting that \\
$\bullet$ $3^n \equiv 1 \;( \bmod \; 10)$ if $4 \mid n$, and $3^n \equiv 9 \;( \bmod \; 10)$ if $2 \mid \mid n$, \\
$\bullet$ $2^v - 1 \equiv 1 \;( \bmod \; 10)$ if $v \equiv 1 \;( \bmod \; 4)$, and $2^v - 1 \equiv 7 \;( \bmod \; 10)$ if $v \equiv 3 \;( \bmod \; 4)$, \\
$\bullet$ $(2^v - 1)^r + 2 \equiv 3 \;( \bmod \; 10)$ for any odd $r$ if $v \equiv 1 \;( \bmod \; 4)$, \\
$\bullet$ $(2^v - 1)^r + 2 \equiv 5 \;( \bmod \; 10)$ for $r \equiv 3 \;( \bmod \; 4)$ and $v \equiv 3 \;( \bmod \; 4)$, and \\
$\bullet$ $(2^v - 1)^r + 2 \equiv 9 \;( \bmod \; 10)$ for $r \equiv 1 \;( \bmod \; 4)$ and $v \equiv 3 \;( \bmod \; 4)$. \\
Only possibility for $(2^{v} - 1)^r + 2 = 3^n$ to hold is when $v, r$ and $n$ take the forms $v = 4k + 3, r = 4l + 1$, and $n = 4z + 2$ with $k, l, z$ non-negative integers. $\blacksquare$ \\

\item {\bf Proposition}: For  $(2^{v = 4k + 3} - 1)^{r = 4l + 1} + 2 = 3^{n = 4z + 2}$, $v$ must be $3$ and thus $(2^{v} - 1)^r + 2 = 3^n$ has only one solution. \\
{\bf Proof}:  \\
Since $2 \mid n$ but $4 \nmid n$, let $n = 2t$, $t$ an odd, positive integer. Rearranging as $(2^{v} - 1)^r + 1 = 3^{2t} - 1 = 9^t - 1$ and expanding yields, $2^v R_5 = 8R_6$ where \\
$R_5 = (2^v - 1)^{r-1} - (2^v - 1)^{r-2} + ... - (2^v - 1) + 1$ and \\
$R_6 = 9^{t-1} + 9^{t-2} + ... + 9 + 1$. \\
$R_5$ and $R_6$ are both, odd. Thus, $v = 3$ so that $(2^{v} - 1)^r + 2 = 3^n$ becomes $7^r + 2 = 3^{2t}$ for which, \textcolor{black}{per \cite{Scott_Styer_2004} (Theorem 6), the only solution corresponds to $r = t = 1,$ that is, \label{concl14pt5} $\boxed{7^1 + 2^1 = 3^2}$}. $\blacksquare$
\end{enumerate}
\end{enumerate} 
\end{enumerate} 

\noindent \rule{14cm}{1pt}


\subsubsection{$p^n - q^r = 2^m$, $r$ even and $n$ odd.}
\label{subsubs3_19_11_25}

\begin{enumerate} 
\item Consider $m > 1$.

\begin{enumerate} 
\item Per Eq. (\ref{dio_3_eq2}), $R_1$ is odd and $R_2$ even. Eq. (\ref{dio_3_eq1}) becomes $2^gR_1 - 2^d R_2 + \gamma - 1 = 2^m$.
\item If $\gamma = -1$, $2^{g-1}R_1 - 2^{d-1} R_2  = 2^{m-1} + 1$, which is odd. $2^{d-1} R_2$ is even for 
 $d \geq 1$. $g$ must be 1, but $p = 2^g + \gamma = 2 - 1$, not possible since $p$ is an odd prime.
\item If $\gamma = 1$, $2^{g}R_1 - 2^{d} R_2  = 2^{m}$. Eq. (\ref{main_dio3}) becomes $(2^{g} + 1)^n - (2^{d} + \beta)^r = 2^m$. \\
 
 \item \label{Maa_20} \label{concl16} \label{concl15a} \textcolor{black}{Per Lemma 4 in \cite{Scott_Styer_2004}, if $P$ and $Q$ are distinct primes (not necessarily Mersenne or Fermat), 
 $r$ is positive and even, $m$ and $n$ are positive integers,  for $n > 1$, $Q^r + 2^m = P^n$ admits 
 at most four solutions: $3^2 + 2^4 = 5^2$, $7^2 + 2^5 = 3^4$, $5^2 + 2 = 3^3$ and $11^2 + 2^2 = 5^3$. 
Proof of Lemma 4 in \cite{Scott_Styer_2004}, potentially alternative to that from Cao \cite{Cao_1986},  
uses Lemmata 1-3, Lemma 6 in \cite{Scott_1993},  and Theorem 13 in \cite{Bender_Herzberg_1979}. 
$(2^{g} + 1)^n = 2^m + (2^{d} + \beta)^r$, $r$ even, $n$ odd, and $m$ any positive integer, has only one solution for $n > 1$, that is, $\boxed{3^3 = 2 + 5^2}$}. \\ \\

\textcolor{black}{For $n = 1$, the relation $2^{g} - 2^m + 1 = x^2$ is considered by Szalay \cite{szalay}, \cite{Scott_2018} (Theorem 1.2; \S \ref{main_dio3_mu_1_r_even}, item \ref{Maa_bobby153aa_guw2020}).  For $g > m$, solutions are $(g, m, x ) = (2t, t+1,2^t-1)$ for positive integer $t > 1$, $(g, m, x) = (5, 3, 5), (7, 3, 11)$ and $(15, 3, 181)$. For $(g, m, x ) = (2t, t+1,2^t-1)$, with $x = (2^d + \beta)^{\frac{r}{2}}$, $(r, \beta, d) = (2, -1, t)$ must hold (Appendix \S \ref{2_k_pm_mu}). $t$ must be $2$ or an odd prime. But, for $g = 2t$, since $2^{2t} + 1$ is prime, $t$ cannot be an odd prime. $t$ must be $2$ for which $2^{g} + 1 = 2^m + [(2^{d} + \beta)^{\frac{r}{2}}]^2$ becomes $2^4 + 1 = \boxed{17 = 2^3 + 3^2}$ \label{concl17}. For $(g, m, x) = (5, 3, 5)$, $2^{g} + 1 = 2^m + [(2^{d} + \beta)^{\frac{r}{2}}]^2$ is $2^5 + 1 = 2^3 + 5^2$ with $2^5 + 1$ a non-prime. In the last two cases, both primes, $11$ and $181$, are not of the form $(2^{d} + \beta)^{\frac{r}{2}}$}. \\

\end{enumerate}

\item Consider $m = 1$.
\begin{enumerate}
\item \label{Maa_21}  \label{concl18} For $\gamma = 1$, Eq. (\ref{dio_3_eq1}) takes the form $2^{g-1}R_1 - 2^{d-1}R_2 = 1$ (odd). \\
 $2^{d-1}R_2$ is even for $d \geq 1$. $g$ must be 1. $p = 3$ and Eq. (\ref{main_dio3}) becomes \\
 $3^n - (2^d + \beta)^r = 2$ of which the only solution per \S\ref{subsubs3_19_11_25} item (\ref{Maa_20}) is $\boxed{3^3 = 5^2 + 2}$.\\

\item  For $\gamma = -1$, Eq. (\ref{dio_3_eq1}) takes the form $2^{g-2}R_1 - 2^{d-2}R_2 = 1$ (odd). \\
 (P1) If $2^{d-2}R_2$ is even, $g$ must be 2 implying $p = 3$. $q = 2^d + \beta$ so that Eq. (\ref{main_dio3}) is $3^n - (2^d + \beta)^r = 2$, a diophantine relation identical to that in \S\ref{subsubs3_19_11_25}  item (\ref{Maa_21}) above. \\ \\
 
 (P2) If $2^{d-2}R_2$ is odd, since $R_2$ is even, $d$ must be 1 and $R_2$ divisible by only 2 and not its power. 
 $q = 2^d + \beta = 2 + \beta = 3$ ($\beta$ must be 1), and $p = 2^g - 1$ with $g$ an odd prime ($p \neq 3$). 
 Eq. (\ref{main_dio3}), which is, $(2^g - 1)^n - 3^r = 2$, is not feasible since
   $(2^g - 1)^n - 2 \equiv (-2)^n - 2 \equiv  -(2^n + 2) \equiv 2 \;( \bmod \; 3)$.
\end{enumerate}
\end{enumerate}


\textcolor{black}{\rule{16cm}{1pt}}

\section{$2^{m}p^{n} + \mu = q^{r}$}
\label{2m_pm_plus_mu_equals_qr}

\noindent We solve the above diophantine equation for various powers of $2$, Mersenne and Fermat primes $p$ and $q$, and $\mu = \pm 1$. We first consider cases where $r$ is even.

\subsection{$2^{m}p^{n} + 1$ as an even power (of a prime)}
\label{solns_2mpn_plus_one_equals_q_2t}
Consider $\mu = 1, p \mbox{ and } q \mbox{ as distinct odd primes},  m,n,r  \mbox{ and } d$ as positive integers.

\begin{enumerate}
\item Let $2^{m}p^{n} + 1 = (2d + 1)^{2}$, or, $2^{m-2}p^{n} = d(d+1)$. $gcd(d,d+1) = 1$ and either $d$ or $d+1$ is even. $d(d+1) > 1$. If $m = 1$ or 2, $p^{n} = 2d(d+1)$ or $d(d+1)$ respectively, are both impossibilities, since $p^{n} (> 1)$ is odd. $m$ must be $ > 2$.

\item  \label{Kaali_M50} $2^{m-2}p^{n} = d(d+1)$ permits two possibilities\footnote{(i) $d = 2^{m-2}$ and $d + 1 = p^{n}  \implies p^{n} = 2^{m-2} + 1$;  (ii) $d = p^{n}$ and $d + 1 = 2^{m-2}  \implies p^{n} = 2^{m-2} - 1$} since $LHS$ is divisible by only two distinct primes: $p^{n} = 2^{m-2} + \gamma$, $\gamma = \pm 1$. Two cases arise: 

\subitem (a) $(m-2,\gamma,p,n) = (3,1,3,2)$, per Appendix \S \ref{2_k_pm_mu}, or  
\subitem (b) $n = 1$ such that $p = 2^{m-2} + \gamma$ is prime.  

\item \label{concl19} For case (a) in item (\ref{Kaali_M50}), $2^{m}p^{n} + 1 = \boxed{2^{5}3^{2} + 1 = 289 = 17^{2} } = q^r $with $q = 17$ and $r = 2$.

\item \label{Kaali_M50aa} \label{concl20} For case (b) in item (\ref{Kaali_M50}), $2^{m}p + 1 = 2^{m}(2^{m-2} + \gamma) + 1 = (2^{m-1} + \gamma)^{2} = q^{r}$. Per Appendix \S \ref{2_k_pm_mu}, $2^{m-1} + \gamma$ is a (an even) perfect power only when $m = 4$ and $\gamma = 1$ in which case, $p = 2^{m-2} + \gamma = 5$, $n = 1, q = 3$ and $r = 4$ yielding $\boxed{2^4\cdot5 +1 = 3^4}$. 

\item Otherwise, in  item (\ref{Kaali_M50aa}), $r = 2$. $p = 2^{m-2} + \gamma$ and $q = 2^{m-1} + \gamma$ must both\footnote{$gcd(2^{m-2} + \gamma,2^{m-1} + \gamma) = gcd(2^{m-2} + \gamma,\gamma) = 1$.}, be prime. If $\gamma = 1$, $m-2$ and $m-1$ must both be powers of 2. If $\gamma = -1$, one of the $m-2$ and $m-1$ must be $2$ and other, an odd prime. 

\item $m > 4$ is an impossibility since either $m-2$ or $m-1$ is odd, and the other even, making either $2^{m-2} + \gamma$ or  $2^{m-1} + \gamma$ composite.

\item \label{concl21} \label{concl22} Valid pairs $(p,q)$ for $m \leq 4$ are $(3,7)$ corresponding to $m = 4, \gamma = -1$, and $(3,5)$ for $m = 3$ and $\gamma = 1$. Correspondingly, the solutions are $\boxed{2^4 \cdot 3 + 1 = 7^2}$ and $\boxed{2^3 \cdot 3 + 1 = 5^2}$ respectively.

\end{enumerate}

\noindent \rule{16cm}{1pt}

\subsection{$2^{m}p^{n} - 1$ as an even power (of a prime)}
\label{solns_2mpn_minus_one_equals_q_2t}

Consider $\mu = -1, p \mbox{ and } q \mbox{ as distinct odd primes}, m,n,r \mbox{ and } d$ as positive integers. \textcolor{black}{In what follows, we share observations on the nature of integer exponents and primes in the diophantine relation $2^m p^n = q^{2t} + 1$}.

\begin{enumerate}

\item Let $2^{m}p^{n} - 1 = (2d + 1)^{2}$, or, $2^{m-1}p^{n} = 2d(d+1) + 1$. No solution is admitted for $m>1$ since the $RHS$ is odd. For $m = 1$, $p,q,n$ and $r = 2t$, $t$ a positive integer, must satisfy ${2p^{n} = q^{2t} + 1}$. 

\item \label{Kaali_M10} {\bf Proposition}: For $n > 2$, $2p^{n} = q^{2t} + 1$ admits no solution if an odd  prime divides $t$. \\
{\bf Proof}:

\begin{enumerate} 
\item Consider $t = av$, $a$ $(\geq 1)$ an integer, and $ v$ an odd prime. 
\begin{eqnarray}
2p^n = q^{2av} + 1 = (q^{2a}+1)R = (q^{2a}+1) \left[1 - q^{2a} + (q^{2a})^2 - ... +  (q^{2a})^{v-1} \right] \nonumber 
\end{eqnarray}
\item  \label{Kaali_M9} $q^{2a}$ and $R$ are odd\footnote{$R$ is odd since the number of terms, $(q^{2a})^x, x = 1,...,v-1$ is even.}. $R = \frac{ q^{2av}+1}{ q^{2a}+1} > q^{2a}+1 > 1$ as $q^{2a(v-1)} >  q^{2a} + 2$ for $v > 2$.

\item From Appendix \S \ref{gcd_of_factors2}, $gcd(q^{2a} + 1, R) = gcd(q^{2a} + 1, v)$ which, is either 1 or $v$.

\item Consider $gcd(q^{2a} + 1, R) = 1$. \\
$p^n = \frac{q^{2a} + 1}{2}R$ is an impossibility for any $n \geq 1$ as $\frac{q^{2a} + 1}{2}> 1, R > 1$, $gcd(\frac{q^{2a} + 1}{2},R) = 1$ making $\frac{q^{2a} + 1}{2}R$ divisible by at least two distinct odd primes. \\

\item Consider $gcd(q^{2a} + 1, R) = v$. \\
For positive integers $\lambda, \beta$, as $q^{2a} + 1$ is even, 
\begin{enumerate}
\item $q^{2a} + 1 = 2\lambda v; \;\;\;\; R = \beta v; \;\;\;\; gcd(2\lambda,\beta) = 1$, and $\beta > 2\lambda > 1$ per item (\ref{Kaali_M9}).
\item \label{Kaali_M8}  $2p^n = (q^{2a}+1)R = 2 \lambda \beta v^2$ implying\footnote{$n$ must be larger than 2 for $p^n = \lambda \beta v^2$ to hold.} either $\lambda = 1 \mbox{ or a power of } p,$ \\
$ v = p, \mbox{ and } \beta$ is a power of $p$.  As $gcd(2\lambda,\beta) = 1$ and $\beta > 2\lambda$, $\lambda$ must be 1.
\item  Also, $\beta = \frac{R}{p} = \frac{1}{p}\frac{ q^{2av}+1}{ q^{2a}+1}  = \frac{1}{p}\frac{(2p-1)^{p} +1 }{2p} = G + 1$ where \\
$G = \frac{(2p)^{p-1}}{p} - \frac{\comb[p]{1}}{p}(2p)^{p-2} +  ... -  \frac{\comb[p]{p-2}}{p}(2p)$.  
\item As $p$ is an odd prime, $p \mid \comb[p]{i}, 1 \leq i \leq p-1$. $p$ also divides $(2p)^{p-1}$ more than once implying $p \mid G$ so that $gcd(G+1,p) = gcd(\beta,p) = 1$, a contradiction to \S\ref{solns_2mpn_minus_one_equals_q_2t} item (\ref{Kaali_M8}).  $\blacksquare$
\end{enumerate}
\end{enumerate}

\item \label{Kaali_M10aaa} For $n > 2$, per  item (\ref{Kaali_M10}), $r = 2t$ must have the form $2^w, w \geq 1$ so that $2p^n = q^{2^w} + 1$. For $w = 1$ and $n = 2$, $q^2 - 2p^2 = -1$ may permit multiple solutions\footnote{\label{footnote_pell} $y^2 - 2x^2 = -1$, the negative Pell equation, has solutions of the form \cite{Beiler66}

\begin{eqnarray}
x \mbox{ or } x(g) = \frac{(\sqrt{2} + 1)^{g} + (\sqrt{2} - 1)^{g} }{2\sqrt{2}} \;\;\;\;\;\;\;\;\;\;\;\; y \mbox{ or } y(g) = \frac{(\sqrt{2} + 1)^{g} - (\sqrt{2} - 1)^{g} }{2} \nonumber 
\label{beiler_eqs}
\end{eqnarray}

where $g > 1$ is an odd integer.  Some examples that are prime pairs, are $[x(3), y(3)] = (5,7)$ and $[x(5), y(5)] = (29,41)$. For some $g$, $[x(g), y(g)]$ may be composite, e.g., $[x(7), y(7)] = (169,239)$ and $[x(9), y(9)] = (985, 1393)$. $[x(7)= 13^2, y(7) = 239]$ satisfies $ y^{2} - 2(x^{\frac{1}{2}})^{4} =  - 1$. \\ }. \textcolor{black}{Also, in $2 p^2 = q^2 + 1 \equiv C = A + B$, $rad(ABC)^{1 + \varepsilon} \equiv (2pq)^{1 + \varepsilon} > 2p^2 \equiv C$ always holds for $\varepsilon > 0$. This is because $2^\varepsilon p^{\varepsilon - 1}q^{1 + \varepsilon} > 1$, or $\varepsilon log(2pq) + log(\frac{q}{p}) > 0$ is always true since $2pq > 2$ and $\frac{q}{p} = \sqrt{2 - \frac{1}{p^2}} > 1$ for $p$ and $q$ odd primes. }

\item One now considers cases when $w > 1$.  Those with $w = 1$ are addressed in \S \ref{solns_2mpn_minus_one_equals_q_2t} item (\ref{w_equals_1}). \\
For $w \geq 2$, $q^{2^{w}}$ takes the form $x^4$, $x$ an odd prime or its (even) power so that 

\begin{eqnarray} 2p^n = x^4 + 1 \;\;\;\;\; \mbox{ or }\;\;\;\;\; p^n = \left(\frac{x^2-1}{2} \right)^2 + \left( \frac{x^2+1}{2} \right)^2  \nonumber  \end{eqnarray}

\item \label{Kaali_M52} {\bf Proposition}: In $2p^n = q^{2^w} + 1$, $p \neq 8k+t, t = 3, 5, 7$, if $n$ ($\geq1$) is odd, and $w > 1$.  \\
{\bf Proof}: 

\begin{enumerate} 
\item Let $x = 2h+1$ and $p = 2g + 1$, where $h,g$ are positive integers.  Then, \\
 $2(p^n - x^2) = (x-1)^2(x+1)^2$ becomes $\left[ gR  - 2h(h+1) \right] = 4 h^2(h+1)^2$, or, 
\begin{eqnarray}gR = 2h(h+1)\left[ 2h(h+1) + 1 \right] \nonumber \end{eqnarray} 
 where  $R = (2g)^{n-1} +  \comb[n]{1}(2g)^{n-2} + ... + \comb[n]{n-2}(2g) + n$. If $n$ is odd, $R$ is odd.
\item  As $2 \mid h(h+1)$, $4$ must divide $gR$, and since $R$ is odd, $4$ must divide $g$ implying $p \equiv 1 \; (\bmod \; 8)$. $\blacksquare$

\textcolor{black}{One also notes that since $s = \left(\frac{x^2-1}{2} \right)^2 + \left( \frac{x^2+1}{2} \right)^2$ is a sum of two squares with $\frac{x^2-1}{2}$ and $\frac{x^2 + 1}{2}$ coprime, per \cite{Euler_maa} proposition IV, $s$ cannot be divided by any number (prime) that itself is not a sum of two squares. Further, a prime is a sum of two squares, per \cite{Fermat_thm_sum_two_sqs1} if and only if $p \equiv 1 \; (\bmod \; 4)$ (that is, $p \not \equiv 3 \mbox{ or } 7\; (\bmod \; 8)$)}. 
\end{enumerate}

\item \label{Kaali_M51} {\bf Proposition}: In $2p^n = q^{2^w} + 1$, $p$ cannot be Mersenne or Fermat prime if $n \mbox{ $(\geq 1)$ }$ is odd, $w > 1$, and $q \neq 5$.\\

{\bf Proof}: 
If $q \neq 5$, $q$ is $\equiv 1, 3, 7 \mbox{ or } 9 \; (\bmod \; 10)$, its second power $\equiv 1, 9, 9 \mbox{ or } 1 \; (\bmod \; 10)$ respectively, third power $\equiv 1, 7, 3 \mbox{ or } 9 \; (\bmod \; 10)$ and fourth power $\equiv 1 \; (\bmod \; 10)$, and so on. Thus, $1 + q^{2^w} \equiv 2 \; (\bmod \; 10)$ for $w \geq 2$ implying\footnote{ \label{footnote_igmr_31_5_24}  If $p \equiv 1 \; (\bmod \; 10)$, $2p^n \equiv 2 \; (\bmod \; 10)$ for any odd $n$. \\
If $p \equiv 3 \; (\bmod \; 10)$, $2p^n \equiv 6 \; (\bmod \; 10)$ for  $n \equiv 1 \; (\bmod \; 4)$, and $2p^n \equiv 4 \; (\bmod \; 10)$ for  $n \equiv 3 \; (\bmod \; 4)$. \\
If $p \equiv 5 \; (\bmod \; 10)$ (i.e., $p = 5$), $2p^n \equiv 0 \; (\bmod \; 10)$ for any $n$. \\
If $p \equiv 7 \; (\bmod \; 10)$, $2p^n \equiv 4 \; (\bmod \; 10)$ for  $n \equiv 1 \; (\bmod \; 4)$, and $2p^n \equiv 6 \; (\bmod \; 10)$ for $n \equiv 3 \; (\bmod \; 4)$. \\
If $p \equiv 9 \; (\bmod \; 10)$, $2p^n \equiv 8 \; (\bmod \; 10)$ for any $n$ odd. \\
}  $2 p^n \equiv 2 \; (\bmod \; 10)$.  For $n$ odd, $p$ must be $\equiv 1 \; (\bmod \; 10)$. From item (\ref{Kaali_M52}), $p \equiv 1 \; (\bmod \; 8)$ must hold as well implying $p$ must be of the form\footnote{For $x, y, z$ positive integers, $p = 10x + 1 = 8y + 1. \;\; 5 \mid y. \;\; y = 5z.$} $40k+1, k \geq 1$, an integer. \\ \\
 Let $p = 2^g + \gamma, \gamma = \pm 1, g \geq 1$ ($g > 1$ if $\gamma = -1$), an integer.
 
\begin{enumerate}
\item For $\gamma = 1, p = 2^g + 1 = 40k + 1$ yields $2^{g-3} = 5k$, not possible \footnote{Per item \ref{Kaali_M52}, since $p = 2^g + 1$ is $\not \equiv 5 \mbox{ or } 3 \; (\bmod \; 8)$ for $g = 1$ or 2 respectively, and $p \neq 2^3 + 1$ for $g = 3$, $g \geq 4$ must hold.}, as $5 \nmid 2^{g-3}$.
\item For $\gamma = -1, 2^g -1 = 40k + 1$ yields $2^{g-1} = 20k + 1$, not possible for $g > 1$ due to opposite parities. $\blacksquare$
\end{enumerate}

\item {\bf Proposition}: In $2p^n = q^{2^w} + 1$, $p$ cannot be $2^g \pm 1$, $g$ a positive integer, if $n \mbox{ $(\geq 1)$ }$ is odd, $w > 1$, and $q = 5$. \\
{\bf Proof}:
\begin{enumerate}
\item  If $q = 5$, $1 + 5^{2^w} \equiv 6 \; (\bmod \; 10)$ so that $2p^n \equiv 6 \; (\bmod \;10)$, 
  possible only when  either $p \equiv 3\; (\bmod \; 10)$ or $p \equiv 7 \; (\bmod \;10)$ since $n$ $(\geq 1)$ is odd 
  (item \ref{Kaali_M51}; footnote \ref{footnote_igmr_31_5_24}). From item (\ref{Kaali_M52}), $p$ must be $\equiv 1 \; (\bmod \; 8)$.  \\
 In addition, let $p = 2^g + \gamma$, $\gamma = \pm 1, g \geq 1$, an integer.

\item  If $p \equiv 3 \; (\bmod \;10)$ and $p \equiv 1 \; (\bmod \;8)$, $p$ must be of the form\footnote{For positive integers $x$ and $y$, $p = 10x + 3 = 8y + 1$. 5 must divide $4y-1$, possible for $y = 5z + 4, z \geq 0$ and integer. Then, $p = 8(5z + 4) + 1 = 40z + 33$.} $40 k + 33$.

\begin{enumerate}
\item If $p = 40k + 33 = 2^g - 1, 2^{g-1} = 20 k + 17$, not possible for $g > 1$ due to opposite parities. $g \neq 1$ since $p \geq 3$. 
\item If $p = 40k + 33 = 2^g + 1, 4(2^{g-5} - 1) = 5k$. $5 \mid (2^{g-5} - 1)$ for $g = 4l + 5, l \geq 0$ and integer. $g$ is odd, not possible, since $p = 2^g + 1$ (prime)  requires $g$ to be a power of 2.
\end{enumerate}
	
\item $p \equiv 7 \; (\bmod \; 10)$ and $p \equiv 1 \; (\bmod \; 8)$: $p$ must be of the form\footnote{$p = 10x + 7 = 8y + 1$. 5 must divide $4y-3$, possible for $y = 5z + 2$, making $p = 8(5z + 2) + 1 = 40z + 17$.} $40 k + 17$.

\begin{enumerate}
\item If $p = 40k + 17 = 2^g- 1, 2^{g-1} = 20 k + 9$, impossible for $g > 1$ due to opposite parities. $g$ cannot be 1.
\item If $p = 40k + 17 = 2^g + 1, 2^{g-3} -2 = 5 k$. $k$ must be even making $p$ of the form $80k+17$. Further, $5 \mid 2^{g-4} - 1$ for $g = 4l + 4, l \geq 0$ and integer. \\ 
$p = 2^g + 1$ being prime, $g = 4(l + 1) \geq 4$ must be a  power of 2. \\
For $g = 2^z, z \geq 2$, $2p^n = 1 + q^{2^w}$ for $q = 5$ becomes $2(2^{2^z} + 1)^n = 1 + (2^2 + 1)^{2^w}$. $2(2^{2^z} + 1)^n - 1$ is divisible by 3 since $2(2^{2^z} + 1)^n - 1 \equiv 2 \cdot 2^n -1 $ \\
$\equiv (-1) \cdot (-1)^n -1 = 0 \; (\bmod \; 3)$ and thus, $2(2^{2^z} + 1)^n - 1$ is not a power of 5. $\blacksquare$
\end{enumerate}
\end{enumerate}

\item \label{check1} {\bf Proposition}: For $n$ even and $w > 1$, $2p^n = q^{2^w} + 1 = x^4 + 1$ has no solution for $p = 2^g + \gamma$, $\gamma = \pm 1$. \\

{\bf Proof}: \\
 A related equation is the Ljunggren equation \cite{Ljunggern, Ljunggern1, Ljunggern2, Ljunggern3}, $2Y^4 = X^2 + 1$, for $Y = p, n = 4$ and $X = x^2$. \textcolor{black}{It is shown therein that the only solution $(Y,X) > (1,1)$ is $(13,239)$}. \\
 
\begin{enumerate}
\item \label{BeeBee4} Thus, if $4 \mid n$, $2p^n = (x^2)^2 + 1$ has no solution.

\item \textcolor{black}{That $x^4 + 1$ is the sum of two squares with $1$ and $x^2$ ($x > 1$) coprime, any odd prime divisor must be a sum of two squares (\cite{Euler_maa}, Proposition IV) and thus of the form $4K + 1$ (\cite{Fermat_thm_sum_two_sqs2}), $K \geq 1$, an integer}. Then, if $p = 2^g + \gamma = 4K + 1$, $\gamma$ must be $1$ since otherwise (if $\gamma = -1$), $2^{g-1} = 2K + 1$ is not possible. For $2^g + 1$ to be prime, $g$ must be a power of 2.
\item Since $x = q$ (a prime) or its even power (\S \ref{solns_2mpn_minus_one_equals_q_2t}, item \ref{Kaali_M10aaa}), per \S \ref{solns_2mpn_minus_one_equals_q_2t}, item (\ref{Kaali_M51}), if $q \neq 5$, \\
$x^4 + 1 \equiv 2 \; (\bmod \; 10)$ so that $2p^n \equiv 2  \; (\bmod \;10)$. If $2 \mid n$ but $4 \nmid n$ (item \ref{BeeBee4}), for $k \geq 1$ an integer, $p$ cannot be of the forms $10k + 3$, $10k+5$ or $10k + 7$. $p$ must then either be of the form $10k + 1$ or $10k + 9$. 
\begin{enumerate}
\item $p = 2^g + 1 = 10k + 1$ is not possible as $5 \nmid 2^{g-1}$.
\item If $p = 2^g + 1 = 10k + 9$, $5 \mid 2^{g-1} - 4 = 4(2^{g-3} - 1)$, that is, $g$ must be of the form $4l + 3$, $l \geq 0$ an integer. But $g$ cannot be odd since $2^g + 1$ is prime. 
\end{enumerate}
\item \label{igmr_10_6_24} If $q = 5$, $x^4 + 1 \equiv 6 \; (\bmod \; 10)$ so that $2p^n \equiv 6 \; (\bmod \; 10)$, not possible if $2 \mid n$. $\blacksquare$
\end{enumerate}

\item \label{w_equals_1} One now considers the form $2p^n = q^2 + 1$, corresponding to $w = 1$ in \S \ref{solns_2mpn_minus_one_equals_q_2t} item (\ref{Kaali_M10aaa}). \\
Since $q^2 + 1$ is a sum of two squares prime between themselves,  per Euler \cite{Euler_maa}, Proposition IV, $p$ must be a sum of two squares, and since $p$ is an odd prime, it must be $\equiv 1 \; (\bmod \;4)$ \cite{Fermat_thm_sum_two_sqs2}. Per item  (\ref{Kaali_M51}), if $q = 5$, $q^2 + 1 = 26$, implying $n = 1$ and $p = 13$, not a Mersenne or Fermat prime. If $q \neq 5$, $q^2 + 1$ is either congruent to $2$ or $0$ $\; (\bmod \; 10)$. We consider two cases: $n$ odd, and $n$ even. 

\begin{enumerate}
\item $n \;(\geq 1)$ odd:\\
Per item (\ref{Kaali_M51}), if $n = 4k + 1$, $k$ a non-negative integer, for $p \equiv 1, 3, 5, 7 \mbox{ or } 9\; (\bmod \; 10) $, $p^{n} \equiv 1, 3, 5, 7 \mbox{ or } 9\; (\bmod \; 10) $ respectively so that $2p^{n} \equiv 2, 6, 0, 4 \mbox{ or } 8\; (\bmod \; 10) $ respectively. \\
If $n = 4k + 3$, for $p \equiv 1, 3, 5, 7 \mbox{ or } 9\; (\bmod \; 10) $, $p^n \equiv 1, 7, 5, 3 \mbox{ or } 9\; (\bmod \; 10) $ respectively so that $2p^n \equiv 2, 4, 0, 6 \mbox{ or } 8\; (\bmod \; 10) $ respectively.
\begin{enumerate}
\item \label{as_ep1} If  $q^2 + 1 \equiv 2 \; (\bmod \; 10)$, $2p^n$ must be congruent to $2 \; (\bmod \; 10)$ implying $p^n$ and thus $p$  $\equiv 1 \; (\bmod \; 10)$ must hold. If  $p \equiv 1 \; (\bmod \; 10)$ and  $p \equiv 1 \; (\bmod \; 4)$, $p$ must be of the form $p = 20K + 1$, $K$ a positive integer. Further, if $p$ is of the form $2^d + \beta$, $d$ an integer and $\beta = \pm 1$, \\
$\bullet$ If $\beta = -1$, $p = 20K + 1 = 2^d - 1$, or, $10K + 1 = 2^{d-1}$ does not hold for $d > 1$ due to opposite parities. \\
$\bullet$ If $\beta =  1$, $p = 20K + 1 = 2^d + 1$, or, $5K = 2^{d-2}$ does not hold for $d \geq 2$ since $5 \nmid  2^{d-2}$.
\item \label{new_soln_eva1} If  $2p^n = q^2 + 1 \equiv 0 \; (\bmod \; 10)$, $p$ must be $5$. \\
If $q = 2^d + \beta$, $d$ an integer and $\beta = \pm 1$, $2p^n = q^2 + 1$ yields, \\
$5^n = 2^{2d-1} + 2^d \beta + 1$, or, $4R = 2^d(2^{d-1} + \beta)$ where $R = (5^{n-1} + 5^{n-2} + ... + 1)$. 
$R \geq 1$ is odd. If $d > 1$, $d$ must be $2$, and since $q$ cannot be $5$, $\beta = -1$ must hold making $q = 3$. Then, for $2 \cdot 5^n = 3^2 + 1$ to hold, $n$ must be $1$. If $d = 1$, $\beta$ must be $1$ making $q = 3$, resulting in the same solution, $\boxed{2 \cdot 5 = 3^2 + 1}$. 
\end{enumerate}
\item $n$ even: \\
Per item \ref{check1}, if $4 \mid n$,  $2p^n = q^2 + 1$ has only one solution: $2 \cdot 13^4 = 239^2 + 1$. \\
One therefore considers cases when $2 \mid \mid n$. Per item (\ref{Kaali_M51}), if $p \equiv 1 \mbox{ or } 9 \; (\bmod \; 10)$, for 
$n = 4k + 2$, $k$ a non-negative integer, $p^{n} \equiv 1 \; (\bmod \; 10) $. \\
If $p \equiv 3 \mbox{ or } 7 \; (\bmod \; 10)$, $p^{4k+2} \equiv 9 \; (\bmod \; 10) $. 
\begin{enumerate}
\item If $q^2 + 1 \equiv 2 \; (\bmod \; 10)$, \\
for $2p^{4k+2} \equiv 2 \; (\bmod \; 10)$ to hold, $p$ must be $\equiv 1 \mbox{ or } 9 \; (\bmod \; 10)$.
\begin{enumerate}
\item $p \equiv 1 \; (\bmod \; 4) $ and $p \equiv 1 \; (\bmod \; 10) $: \\
Per item (\ref{as_ep1}), $p$ cannot be a Mersenne or Fermat prime.
\item $p \equiv 1 \; (\bmod \; 4) $ and $p \equiv 9 \; (\bmod \; 10) $: \\
$p$ must be of the form $20K + 9$, $K$ a positive integer.  Further, if $p = 2^d + \beta$, $d$ an integer and $\beta = \pm 1$,\\
$\bullet$  if $\beta = 1$, $2^d + 1 = 20K + 9$, or, $2(2^{d-3} - 1) = 5K$ must hold. $d$ must be $4l + 3$, $l \geq 0$, an integer. $d$ is odd, not possible, as then $2^d + 1$ is not prime.\\
$\bullet$  if $\beta = -1$, $2^d - 1 = 20K + 9$, or, $2^{d-1} = 5(2K + 1)$ must hold, not possible, since $5 \nmid 2^{d-1}$ for $d \geq 1$. 
 \end{enumerate}
 \item If $q^2 + 1 \equiv 0 \; (\bmod \; 10)$, \\
$p$ must be $5$. Per item (\ref{new_soln_eva1}), $n$ cannot be even. 

\item \textcolor{black}{$q$ cannot be $5$ since then $q^2 + 1 \equiv 6 \, (\bmod 10)$, not possible for $2p^n $, $p$ prime and $n$ even, per item  \ref{igmr_10_6_24}}. 
\end{enumerate}
\end{enumerate}
\end{enumerate}

\noindent \rule{16cm}{1pt}

\subsection{$2^{m}p^{n} + \mu, \mu = \pm 1,$ as an odd power (of a prime)}
\label{2mpn_plus_mu_equals_qr_r_odd}

Consider $2^{m}p^{n} + \mu = q^{r} = (2D+\mu)^{r}$, $r > 1$ and odd; $D \geq 1$ ($D > 1$ for $\mu = -1$). $p$ and $q$ are distinct (not necessarily Mersenne or Fermat primes) odd primes, and $m$ and $n$ are positive integers. 

\begin{enumerate}
\item \label{2mpn_plus_mu_equals_qr_r_odd_gcd_mn} {\bf Proposition}: In $2^{m}p^{n} + \mu = q^r$, $gcd(m,n) = 1$. \\
{\bf Proof}:
\begin{enumerate}
\item \label{Kaali_M11} Let $m = au$, $n = av$, $a$ an odd prime, $u,v \geq 1$. \\
$2^mp^n + \mu = (2^up^v)^{a} + \mu = (2^up^v + \mu)R$ where \\
$R = (2^up^v)^{a-1} - \mu (2^up^v)^{a-2} +  (2^up^v)^{a-3} - ... - \mu (2^up^v) + 1$. \\
$R > 2^up^v + \mu$ as 
$(2^u p^v)^{a-2} > 2 > 1 + \frac{2 \mu}{2^u p^v} + \frac{1 - \mu}{(2^u p^v)^2}$. 
\item Per Appendix \S \ref{gcd_of_factors2}, $gcd(2^up^v + \mu,R) = gcd(2^up^v + \mu,a) = 1$ or $a$.
\item If  $gcd(2^up^v + \mu,R) = 1$, $q^r \neq 2^mp^n + \mu$ as two or more primes divide the $RHS$.
\item If  $gcd(2^up^v + \mu,R) = a$, for positive integers $\lambda \mbox{ and } \beta$, let $2^up^v + \mu = \lambda a$, $R = \beta a$, $gcd(\lambda, \beta) = 1$ and $\beta  > \lambda $ (item \ref{Kaali_M11}).  Then, $q^r = 2^mp^n + \mu = (2^up^v + \mu)R = \lambda \beta a^2$. \\
$a$ must be $q$, and $q$ must divide $\beta$. \\ 
$\beta = \frac{R}{a} = \frac{1}{a}\frac{(2^up^v)^{a} + \mu}{2^up^v + \mu} = \frac{(\lambda a - \mu)^a + \mu}{\lambda a^2} = \frac{(\lambda a)^{a-1} -  \comb[a]{1} \mu(\lambda a)^{a-2} + ... - \comb[a]{a-2} \mu(\lambda a)}{a} + 1 = G + 1$. \\
$a \mid \comb[a]{i}, 1 \leq i \leq a-1$ and $a \mid (\lambda a)^j, j = 1, ... , a-1$. $a \mid G$ and so, \\ $gcd(\beta, a = q) = 1$, a contradiction. 
\end{enumerate}

From the foregoing, $gcd(m,n)$ must be $2^w, w \geq 0$. If $m$ and $n$ are even, 
\begin{enumerate}
\item If $\mu = -1$, $2^mp^n - 1 = (2^{\frac{m}{2}}p^{\frac{n}{2}} -1)(2^{\frac{m}{2}}p^{\frac{n}{2}} + 1) \neq q^r$ since \\
$gcd(2^{\frac{m}{2}}p^{\frac{n}{2}} -1,2^{\frac{m}{2}}p^{\frac{n}{2}} +1) = 1$ and $2^{\frac{m}{2}}p^{\frac{n}{2}} -1 > 1$.
\item If $\mu = 1$, $\left[ 2^{\frac{m}{2}}p^{\frac{n}{2}} \right]^2 + 1 \neq  q^r$  \textcolor{black}{per the Catalan-Mih\u{a}ilescu theorem (\cite{Catalan_Mihailescu})}. $\blacksquare$
\end{enumerate}

\item \label{Kaali_M14} Expanding $2^{m}p^{n} + \mu = q^{r} = (\mu + 2D)^{r}$, as $\mu^{j} = \mu$ (or 1) if $j$ is odd (or even), and $\comb[r]{1} = r$,
\vspace{-1mm}
\begin{eqnarray} 
2^{m}p^{n} = (2D)H = 2D\left[ r + \comb[r]{2}\mu (2D) + \comb[r]{3}(2D)^{2} + \comb[r]{4}\mu(2D)^{3} + ... + (2D)^{r-1} \right] \nonumber
\end{eqnarray}

 $H$ is odd since $r$ is odd. Also, $H > 1$ as $r > 1$. $2^{m} \nmid H$ and so $2^{m}$ must divide $2D$. Let $2D = 2^{m}\lambda$, $\lambda$ a(n odd) positive integer. Then
\vspace{-1mm}
\begin{eqnarray} 
p^{n} = \lambda \left[  r + \comb[r]{2}\mu (2^{m}\lambda) + \comb[r]{3}(2^{m}\lambda)^{2} + ... + (2^{m}\lambda)^{r-1} \right] \nonumber
\end{eqnarray}

$\lambda$ must assume the form $p^{g}$, $0 \leq g < n$. With $q = 2D + \mu = 2^mp^g + \mu$ (prime), $2^mp^n + \mu = q^r $ becomes (see item \ref{Kaali_M12})
\begin{eqnarray}
p^{n-g} = \frac{(2^mp^g + \mu)^r - \mu}{(2^mp^g + \mu) - \mu} = \frac{q^r - \mu}{q - \mu} \nonumber 
\end{eqnarray}

\item \label{Kaali_M18} {\bf Proposition}: \textcolor{black}{In $p^{n-g} = \frac{q^r - \mu}{q - \mu},$ $r$ must be an odd prime}. \\
{\bf Proof}: \\
Consider $r = av$, $a > 1$ and odd, and $v$ an odd prime. 
\begin{enumerate}
\item \label{Kaali_M13} $q^{av} - \mu = (q^a - \mu)(1 + q^a \mu + q^{2a} + q^{3a} \mu + ... + q^{a(v-1)}) = (q^a- \mu)R$.  $R$ is odd\footnote{$q^a \mu,  q^{2a}, q^{3a} \mu, ..., q^{a(v-1)}$ are even number of terms, each odd.}. $\;$ $R > q^a-\mu > 1$ since $q^{av} - \mu > (q^a-\mu)^2$, or, $q^{a(v-2)}  > 2 > 1 - \frac{2}{q^a}\mu + \frac{1 + \mu}{q^{2a}}$.
\item Per Appendix \S \ref{gcd_of_factors2}, $gcd(q^a-\mu,R) = gcd(q^a-\mu,v) = 1$ or $v$.
\item If $gcd(q^a-\mu,R) = 1$, $(q^a-\mu)R$ is divisible by two or more distinct primes. \\
Further, $\frac{q^a-\mu}{q-\mu} > 1$. $p^{n-g} = \frac{q^{av} - \mu}{q - \mu}$ does not hold.
\item Consider $gcd(q^a-\mu,R) = v$. Then, $q^a - \mu = \lambda v, R = \frac{q^{av}-\mu}{q^a-\mu} = \beta v$, $\beta > \lambda$ from item (\ref{Kaali_M13}) and  $gcd(\lambda,\beta) = 1$  for positive integers $\lambda$ (even) and $\beta$ (odd).
\begin{eqnarray} \beta = \frac{1}{v}\frac{(\lambda v + \mu)^v-\mu}{\lambda v} =  \frac{1}{v}[ (\lambda v)^{v-1} + \comb[v]{1} (\lambda v)^{v-2}\mu + ... + \comb[v]{v-2} (\lambda v) \mu]  + 1  \nonumber 
\end{eqnarray}

$v^2$ divides $\left[ (\lambda v)^{v-1} + \comb[v]{1} (\lambda v)^{v-2}\mu + ... + \comb[v]{v-2} (\lambda v) \mu \right]$ as $v \mid \comb[v]{i}, 1 \leq i \leq v-1$. Thus, $gcd(\beta,v) = 1$. $\frac{q^{av} - \mu}{q-\mu} = \frac{\lambda \beta v^2}{q-\mu} = \frac{\beta v(q^a-\mu)}{q-\mu} = \beta v (1 + q\mu + q^2 + ... + q^{a-1})$ is  divisible by at least two distinct odd primes. $p^{n-g} = \frac{q^{av} - \mu}{q - \mu}$ does not hold. $\blacksquare$
\end{enumerate}

\item \label{Kaali_M19} Consider $q = 2^d + \beta$, $\beta = \pm 1$ and $d \geq 1$ (integer). Then, from \S \ref{2mpn_plus_mu_equals_qr_r_odd} item (\ref{Kaali_M14}), $q = 2^mp^g + \mu = 2^d + \beta$.\\
If $\beta = -\mu$, 
\begin{enumerate}
\item If $d > 1$,\\
$2^{m-1}p^g= 2^{d-1} -\mu$. Since $RHS$ is odd ($d > 1$), $m$ must be 1 implying $p^g = 2^{d-1} -\mu$. Per Appendix \S \ref{2_k_pm_mu}, possibilities are (P1): $(p,g,d,\mu) = (3, 2, 4, -1)$, or (P2): $g = 1$. 

\begin{enumerate}
\item \label{Kaali_M15} In (P1), $q = 2\cdot 9 - 1 = 17$ so that $2^{m}p^{n} + \mu = q^{r}$ becomes $2\cdot 3^n = 17^r + 1$.  
\item In (P2), $q = 2^{d} - \mu$ and $p = 2^{d-1} - \mu$ must both be prime. 
\begin{enumerate} 
\item \label{Kaali_M16} For $\mu = 1$, $d-1$ or $d$ must be 2 or an odd prime, only possibility being 
$d = 3$, so that $(p,q) = (2^{d-1} - 1, 2^{d} - 1) = (3, 7)$. $2^{m}p^{n} + \mu = q^{r}$ becomes $2\cdot 3^n = 7^r - 1$.
\item \label{Kaali_M17} For $\mu = -1$, $d-1$ and $d$ must both be powers of 2, only possibility being 
$d-1 = 2^{0}$ and $d = 2^{1}$ so that $(p,q) = (2^{d-1} +1, 2^{d} +1) = (3,5)$.   $2^{m}p^{n} + \mu = q^{r}$ 
becomes $2\cdot 3^n = 5^r + 1$.
\end{enumerate}

\item $2\cdot 3^n = a^r + \mu$ generically represents the diophantine relations in items (\ref{Kaali_M15}), (\ref{Kaali_M16}) and (\ref{Kaali_M17}) where $a = 17, 7 \mbox{ and } 5$ and \textcolor{black}{$\mu = 1, -1 \mbox{ and } 1$} respectively.\\ 
\textcolor{black}{$2 \mid (17^r + 1)$, $2 \mid (7^r - 1)$ and $2 \mid (5^r + 1)$ for any positive integer $r$. Likewise, $3 \mid (17^r + 1)$, $3 \mid (7^r - 1)$ and $3 \mid (5^r + 1)$ for any odd $r$ (prime, item \ref{Kaali_M18}). Per Zsigmondy's theorem (\cite{Zsigmondy_1892, Birkhoff_Vandiver_1904}), there exists a prime other than $2$ and $3$ that divides $ (17^r + 1)$, $ (7^r - 1)$ and $ (5^r + 1)$ for $r > 3$. $2 \cdot 3^n = a^r + \mu$ for $a = 17, 7 \mbox{ and } 5$ and \textcolor{black}{$\mu = 1, -1 \mbox{ and } 1$} respectively, does not hold. }
\end{enumerate}

\item If $d = 1$,\\
$\beta = 1$ must hold so that $q =3$. Then, $2^m p^g = 3 - \mu$. If $\mu = 1$, $m = 1$ and $g = 0$. If $\mu = -1$, $m = 2$ and $g = 0$. Cases $p^n = \frac{3^r - 1}{2}$ and $p^n = \frac{3^r + 1}{4}$ are addressed next, in \S \ref{app_conjecture_A}.

\end{enumerate}

\item \label{Kaali_M12} \textcolor{black}{In item (\ref{Kaali_M19}), if $\beta = \mu$, $d = m$ and $g = 0$. Then, $p^{n} = \frac{q^r - \mu}{q - \mu} =  \frac{(2^m + \mu)^r - \mu}{2^m}$. The equation is a specific case of Nagell–Ljunggren (\cite{Waldschmidt_2009, Bennett_Levin_2013}) relation\footnote{ \label{foot_ljunggern} a total of four solutions are known (though not yet proved that these are the only ones \cite{Bennett_Levin_2013})  which are: $\frac{7^4-1}{7-1} = 20^2$, $\frac{18^3 - 1}{18-1} = 7^3 = \frac{(-19)^3 - 1}{-19 - 1}$ and $11^2 = \frac{3^5 - 1}{3 - 1}$.} $y^z = \frac{x^n - 1}{x - 1}$, $x,y > 1$ and $n > 2$. 
$y^z = \frac{x^n - 1}{x - 1}$, when rewritten as $x^n - (x-1)y^z = 1$, is also a specific case of the generalized Pillai equation (\cite{gp1,gp2,gp3,gp4})}. 


\item \label{Kaali_M12pqr} \textcolor{black}{Following Nagell's approach \cite{Nagell1, Nagell2}  to $y^z = \frac{x^n - 1}{x - 1}, x > 1, y > 1, n > 2, z \geq 2$, Ljunggren \cite{Ljunggern4}, using\footnote{translation of a part solution is provided in $https://mathoverflow.net/questions/206645/on-a-result-attributed-to-w-ljunggren-and-t-nagell;$, accessed 11th June, 2024} Mahler's theorem \cite{Mahler} shows (\cite{Bugeaud_Mihailescu}) that apart from certain solutions (footnote \ref{foot_ljunggern}), $y^z = \frac{x^n - 1}{x - 1}$ has no other solution if (i) $z = 2$, (ii) $3 \mid n$, (iii) $4 \mid n$ and (iv) $z = 3$ and $n \not \equiv 5 \; (\bmod \; 6)$}. 

\item With $p = 2^g+ \gamma$ and $q = 2^m + \mu$ as Mersenne/Fermat primes, it behooves to consider $n$ odd in $p^{n} = \frac{q^r - \mu}{q - \mu}$ or $ (2^g + \gamma)^n =  \frac{(2^m + \mu)^r - \mu}{2^m}$ to further analyze the diophantine for $(g, m,  r, n, \gamma, \mu)$. 

\begin{enumerate}
\item Consider $\mu = 1$ so that $q = 2^m + 1$. $p^{n} =\frac{(2^m + \mu)^r - \mu}{2^m}$ becomes $2^mp^n + 1 = (2^m + 1)^r$.\\
For $2pq > \sqrt{2^mp^n + 1}$ to hold, $4p^2q^2 > 2^mp^n + 1 = (q-1)p^n + 1$, or, \\
$4q >  \frac{(q-1)}{q}p^{n-2} + \frac{1}{qp^2}$ must be true. The inequality holds for $p, q \geq 3$ and $n \leq 2$ since $\frac{(q-1)}{q}p^{n-2} + \frac{1}{qp^2} < 2$.
\item With $\mu = -1$,  $q = 2^m - 1$. $p^{n} =\frac{(2^m + \mu)^r - \mu}{2^m}$ becomes $2^mp^n = (2^m - 1)^r + 1$. \\
$2pq > \sqrt{2^mp^n}$ holds when $4p^2q^2 > 2^mp^n = (q + 1)p^n$, or, $4q >  \frac{(q + 1)}{q}p^{n-2}$, true for $p, q \geq 3$ and $n \leq 2$ as $ \frac{(q + 1)}{q}p^{n-2} < 2$.
\end{enumerate}

\end{enumerate}

\noindent \rule{16cm}{1pt}

\subsubsection{For $p^{n} = \frac{q^r - \mu}{q - \mu}$ where $p = 2^g + \gamma$, $q = 2^m + \mu$ ($p \neq q$) and $r$ are all odd primes, $n$ is odd and $\mu, \gamma = \pm 1$, two solutions are $(g, r, m, \mu, n) = (3, 3, 2, -1, 1) \mbox{ and } (5, 3, 2, 1, 1)$.}
\label{app_conjecture_A}

\begin{enumerate}

\item \label{Maa_bobby2} \textcolor{black}{Every prime divisor of $\frac{(2^m + \mu)^r - \mu}{(2^m + \mu) - \mu}$, $r$ an odd prime, is either $r$ or $D_1 r + 1$, $D_1 \geq 1$, an integer  (\cite{Ge_Yimin}, Corollary 6)}.
\begin{enumerate} 
\item \label{Kaali_M53} $(2^m + \mu)^r - \mu = (2^m + \mu)[(2^m + \mu)^{r-1} - 1] + 2^m$. \\
$r$ divides\footnote{$r$ and $2^m + \mu$ are odd primes. Either $r = 2^m + \mu$, or, $r \mid (2^m + \mu)^{r-1} - 1$ per the Fermat's little theorem.} $(2^m + \mu)[(2^m + \mu)^{r-1} - 1]$, but $r \nmid 2^m$. Thus, $p = 2^g + \gamma = D_1 r + 1$.
\item \label{Kaali_M53_eves} $\gamma$ cannot be $1$ since $r \nmid 2^g$. $p$ must be $2^g - 1$. Then, $D_1 r = 2(2^{g-1} - 1)$. \\
Also, $gcd(p,r) = gcd(2^g-1 = D_1r +1,r) = 1$, implying that $p$ and $r$ are coprime. 
\item \label{Kaali_M53_aaa} Since $p$ is prime, $g$ must be 2 or an odd prime. If $g = 2$, $D_1 r = 2$ is not possible. \\ 
$g$ must be an odd prime. 
\item \label{Kaali_M53abcdef} Per item (\ref{Kaali_M53_eves}), $r \mid 2^{g-1} - 1 = (2^{\frac{g-1}{2}} - 1)(2^{\frac{g-1}{2}} + 1)$, product of two coprimes, implying either $r \mid (2^{\frac{g-1}{2}} - 1)$ or $r \mid (2^{\frac{g-1}{2}} + 1)$. Moreover, $r \leq 2^{\frac{g-1}{2}} + 1$. 
\item Since $g \mid 2^{g-1} - 1$ (Fermat's Little theorem), if $r \neq g$, $g \mid D_1$. Thus, $p = 2^g - 1 = D_1 r + 1 \equiv 1 \; (\bmod \; r)$ and $\equiv 1 \; (\bmod \; g)$.
\item \label{Kaali_M53abcd} $2^g - 1 > \left(2^{\frac{g-1}{2}} + 1\right)^2$ for $g > 3$. Thus, $p = 2^g - 1 > \left(2^{\frac{g-1}{2}} + 1\right)^2 \geq r^2$ per item (\ref{Kaali_M53abcdef}) for $g > 3$.

\end{enumerate}

\item For a positive odd integer $n$, $(2^g - 1)^n = \frac{(2^m + \mu)^r - \mu}{2^m}$ can be rearranged as
\begin{enumerate}
\item \label{Maa_bobbyB} $2^g[(2^g)^{n-1} - \comb[n]{1} (2^g)^{n-2} +  \comb[n]{2} (2^g)^{n-3} -  ... + n] = $ \\
$2^m[(2^m)^{r-2} + \comb[r]{1} (2^m)^{r-3}\mu + \comb[r]{2} (2^m)^{r-4} + ... + \comb[r]{r-3}2^m + \comb[r]{r-2}\mu] + r + 1$.
\item If $m \geq g$, $2^g$ must divide $ r + 1$, or, $r = 2^g D_2 -1 $ for $D_2 \geq 1$, an integer. \\
Then, $r = 2^g D_2 - 1 \geq p = 2^g-1$.  But, $p = D_1  r + 1 > r$, (item \ref{Kaali_M53}), and moreover $p > r^2$ for $g > 3$ (item \ref{Kaali_M53abcd}), a contradiction. 
\item \label{Maa_bobby1} $g > m$ must hold in which case $2^m$ must divide $r + 1$, or, $r = 2^m D_3 -1$ for $D_3 \geq 1$, an integer.
\item \label{Maa_bobby111} Per item (\ref{Kaali_M53abcdef}), $2^{\frac{g-1}{2}} + 1 \geq r = 2^m D_3 -1$ implies  $2^{\frac{g-3}{2}} \geq 2^{m-1} D_3 - 1$. \\
$2^{m-1} D_3 - 1 \geq 2^{m-2}$ for $m \geq 2$ \textcolor{black}{(item \ref{chiggs3} below)} since $D_3 \geq 1$. Thus, $\frac{g-3}{2} \geq m - 2$, or, $g \geq 2m - 1$.
\item As $2^m + \mu$ is prime, either $m$ is an odd prime, or, it is a perfect power of $2$. In case the latter, $gcd(m,r) = 1$. In case $m$ is an odd prime, either $r = m$ or $gcd(m,r) = 1$. $r = 2^m D_3 - 1 = m$ does not hold for $m > 1$ since $D_3 = \frac{m+1}{2^m} < 1$, not permitted. $m$ cannot be $1$ since then $r = 1$, not possible. $m$ and $r$ are therefore coprimes. 

\end{enumerate}

\item  \label{chiggs2} Equation in item (\ref{Maa_bobbyB}) can be rewritten as, 
\begin{enumerate}
\item \label{chiggs2aaa} $2^{g-m}[(2^g)^{n-1} - \comb[n]{1} (2^g)^{n-2} +  \comb[n]{2} (2^g)^{n-3} - ... + n] = $ \\
$2^m[(2^m)^{r-3} + \comb[r]{1} (2^m)^{r-4}\mu + \comb[r]{2} (2^m)^{r-5} + ... + \comb[r]{r-3}] + \comb[r]{r-2}\mu + D_3$.\\
As $g > m$ per item (\ref{Maa_bobby1}), for $m \geq 1$, $\comb[r]{r-2}\mu + D_3$ must be even. Also, since $r = 2^m D_3 - 1$, for $m > 1$, $r \equiv 3 \; (\bmod \; 4)$. 

\item \label{chiggs3} Thus, if $r = 4K + 3$, $K$ a non negative integer, $\comb[r]{r-2} = \comb[r]{2} = (4K+3)(2K + 1)$ is odd implying $D_3$ must be odd. In case $m = 1$, $\mu$ must be 1 so that $q = 3$. $\comb[r]{r-2}\mu + D_3 = (2D_3 - 1)(D_3 - 1) + D_3 = 2D^{2}_3 - 2D_3 + 1$  cannot be even for $D_3 \geq 1$. Thus,  $m \geq 2$ must hold. 

\item \label{igmr18062024b} Pertaining to item (\ref{chiggs2aaa}), $g - m \neq m$ since $g$ is an odd prime (item \ref{Kaali_M53_aaa}). Thus, either $g - m < m$ or $g - m > m$. In particular, for the latter case, as $m \nmid g$, $(I+1)m > g-m > I m$ must hold for some positive integer $I$. 

\end{enumerate}

\item $g - m < m$ 
\begin{enumerate}
\item \label{Shambhu1} Per item (\ref{Maa_bobby111}), $g \geq 2m - 1$. Since $g - m < m$, $g = 2m - 1$ must hold. Equation in item (\ref{chiggs2aaa}) becomes\\
$2^{m-1}[(2^g)^{n-1} - \comb[n]{1} (2^g)^{n-2} +  \comb[n]{2} (2^g)^{n-3} - ... - \comb[n]{n-2} (2^g) + n] = $ \\
$2^m[(2^m)^{r-3} + \comb[r]{1} (2^m)^{r-4}\mu + \comb[r]{2} (2^m)^{r-5} + ... + \comb[r]{r-3}] + \comb[r]{r-2}\mu + D_3$, or, \\
$2^{m-1}R_1 = 2^m R_2 + \comb[r]{r-2}\mu + D_3$ with \\ \\
$R_1 = (2^g)^{n-1} - \comb[n]{1} (2^g)^{n-2} +  \comb[n]{2} (2^g)^{n-3} - ... - \comb[n]{n-2} (2^g) + n$ and \\
$R_2 = (2^m)^{r-3} + \comb[r]{1} (2^m)^{r-4}\mu + \comb[r]{2} (2^m)^{r-5} + ... + \comb[r]{r-4} (2^m) \mu + \comb[r]{r-3}$.
\item For $n$ odd, $R_1$ is odd. As $r = 2^mD_3 - 1 = 4K + 3$ (item \ref{chiggs3}) for some positive integer $K$, $\comb[r]{r-3} = \comb[r]{3} = \frac{(4K+3)(2K+1)(4K+1)}{3}$ is odd making $R_2$ odd. \\
$2^{m-1}$ must divide $\comb[r]{r-2}\mu + D_3$ and further, $\frac{\comb[r]{r-2}\mu + D_3}{2^{m-1}}$ must be odd. That is, $ 2^{m}D_3^2 \mu - 2D_3 \mu - D_3 \mu + \frac{D_3 + \mu}{2^{m-1}}$ must be odd. As $D_3$ is odd (item \ref{chiggs3}), $\frac{D_3 + \mu}{2^{m-1}}$ must be even, or, $D_3 = 2^{m}D_4 - \mu$, $D_4$, an integer (including $0$; \textcolor{black}{Note that if $D_4 = 0$, $\mu$ must be $-1$ since $D_3 \geq 1$ per item (\ref{Maa_bobby1})}). Then, $r = 2^mD_3 - 1 = 2^{2m}D_4 - 2^m \mu - 1$. 

\item Per item (\ref{Kaali_M53}), $2^g - 1 = D_1 r + 1$ where $D_1$ is an even positive integer. Using $g = 2m-1$ (item \ref{Shambhu1}), $2^{2m-1} - 1 = D_1 (2^{2m}D_4 - 2^m \mu - 1) + 1$, or, $D_1 = 2 \frac{(2^{2m-2} - 1)}{2^{2m}D_4 - 2^m \mu - 1}$. As $D_1 \geq 2$, $D_4 \leq \frac{2^{2m-2} + 2^m \mu}{2^{2m}} = 2^{-2} + 2^{-m} \mu$. 
\item $2^{-2} + 2^{-m} \mu < 1$ for $\mu = \pm 1$ and $m \geq 2$.  As $D_3 \geq 1, D_4 \geq \frac{1 + \mu}{2^m} \geq 0$. $D_4$ then must be zero making $D_3 = - \mu$, and since $D_3$ is positive (item \ref{Maa_bobby1}) and odd (item \ref{chiggs3}), $\mu$ must be $-1$. $r$ then becomes $2^m - 1$.

\item Per item (\ref{Kaali_M53abcdef}), $r \mid 2^{g - 1} - 1$, or $2^m - 1$ divides $2^{2m-2} - 1 = (2^{m-1} - 1)(2^{m-1} + 1)$. \\
$2^m - 1$ cannot divide $2^{m-1} - 1$. For $2^m - 1$ to divide $2^{m-1} + 1$, $2^{m-1} + 1 \geq 2^{m} - 1$, or, $2 \geq 2^{m-1}$ must hold, possible only for $m = 2$. 
\item \label{concl23} (S1) $\left(2^{2m-1} - 1 \right)^n = \frac{(2^m - 1)^{2^m - 1} + 1}{2^m}$ becomes $\boxed{\left(2^{3} - 1 \right)^n = \frac{3^3 + 1}{2^2}}$ for which $n = 1$. \\
\end{enumerate}

\item $2m > g - m > m$
\begin{enumerate}
\item  \label{just_label1} In \\
$2^{g-m}[(2^g)^{n-1} - \comb[n]{1} (2^g)^{n-2} +  \comb[n]{2} (2^g)^{n-3} - ... + n] = $ \\
$2^m[(2^m)^{r-3} + \comb[r]{1} (2^m)^{r-4}\mu + \comb[r]{2} (2^m)^{r-5} + ... + \comb[r]{r-4}2^m \mu + \comb[r]{r-3}] + \comb[r]{r-2}\mu + D_3$ \\ \\
(item \ref{chiggs2aaa}), as $g - m > m$, $\comb[r]{r-2}\mu + D_3$ must be divisible by $2^m$. \\
Using $r = 2^mD_3 - 1$ (item \ref{Maa_bobby1}), $\comb[r]{r-3} = \frac{(2^mD_3 - 1)(2^{m-1}D_3 - 1)(2^mD_3 - 3)}{3}$ is odd for $m \geq 2$ (item \ref{chiggs3}). $\frac{\comb[r]{r-2}\mu + D_3}{2^m}$ must be odd. 
\item \label{just_label1a} As $D_3$ is an odd positive integer (item  \ref{chiggs3}),  $2^m$ divides $\comb[r]{r-2}\mu + D_3 = 2^{2m-1}D_3^2 \mu - 2^m D_3 \mu - 2^{m-1}D_3 \mu + D_3 + \mu$ if $D_3 + \mu = 2^{m-1}D_4$, $D_4$ an odd positive integer\footnote{$D_4$ cannot be zero for if it is, then $D_3 = -\mu = 1$. $r = 2^mD_3 - 1 = 2^m - 1$. $2^m \nmid \comb[r]{2}\mu + D_3 = -2^{2m-1} + 2^m + 2^{m-1}$.}, and $D_4 - D_3 \mu$ is even. Then, $r = 2^mD_3 - 1 = 2^m(2^{m-1}D_4 - \mu) - 1 = 2^{2m-1}D_4 - 2^m \mu - 1$. \\

\item {\bf Proposition}: If $2^t$ precisely divides $\comb[r]{r-3} + \frac{\comb[r]{r-2}\mu + D_3}{2^m}$ leaving an odd quotient for a positive integer $t$ and if $t < m$, $(2^g - 1)^n = \frac{ (2^m + \mu)^r - \mu }{2^m} $ yields only one solution. \\

{\bf Proof}:
\begin{enumerate}
\item \label{igmr18062024} Let $\comb[r]{r-3} + \frac{\comb[r]{r-2}\mu + D_3}{2^m} = 2^t X$, $X$, an odd integer. For the expression in item (\ref{just_label1}) to hold, $g - m = m + t$, or $g = 2m + t$ must be true. $t$ must be odd since $g$ is an odd prime (item \ref{Kaali_M53_aaa}). 
\item Per item (\ref{Kaali_M53abcdef}) $r \mid 2^{g-1} - 1$. That is, $2^{2m-1}D_4 - 2^m \mu - 1$ divides \\
$2^{2m + t-1} - 1 = \left[2^{\frac{2m + t - 1}{2}} - 1\right] \left[ 2^{\frac{2m + t - 1}{2}} + 1 \right]$. Since \\
$gcd(2^{\frac{2m + t - 1}{2}} - 1, 2^{\frac{2m + t - 1}{2}} + 1) = 1$, $r$ divides either $2^{\frac{2m + t - 1}{2}} - 1$ or $2^{\frac{2m + t - 1}{2}} + 1$. 
\item \label{just_label1b} If $r$ divides $2^{\frac{2m + t - 1}{2}} - 1$, $ 2^{2m-1}D_4 - 2^m \mu - 1 \leq 2^{\frac{2m + t - 1}{2}} - 1 < 2^{\frac{3m - 1}{2}} - 1$ must hold. That is, $D_4 < \frac{1}{2^{\frac{m-1}{2}}}  + \frac{1}{2^{m-1}} \mu \leq  \frac{1}{2^{\frac{m-1}{2}}}  + \frac{1}{2^{m-1}}$. \textcolor{black}{If $m = 2, D_4 \leq \frac{1}{\sqrt{2}} + \frac{1}{2} \approx 1.2071$, and if $m = 3, D_4 \leq  \frac{1}{2^{\frac{3-1}{2}}}  + \frac{1}{2^{3-1}} = 0.75$.}

As $D_4 \geq 1$ and odd (item \ref{just_label1a}), and since $m \geq 2$ (item \ref{chiggs3}), only possibility is $m = 2$. Then, $D_4 = 1$ and $t = 1$ (since $t < m$) making $g = 2m + t = 5$. $r = 2^{2m-1}D_4 - 2^m \mu - 1 = 2^3 - 2^2 \mu - 1$. Possibilities are
\begin{enumerate}
\item \label{concl24} $\mu = 1$: (S2)  $r = 3$. $(2^g - 1)^n = \frac{ (2^m + \mu)^r - \mu }{2^m} $ becomes $\boxed{(2^5 - 1)^n = \frac{ 5^3 - 1}{2^2} }$ for which $n = 1$.
\item $\mu = -1$:  $r = 11$. $\frac{ 3^{11}  + 1}{2^2} = 67 \times 661$ is not of the form $(2^g - 1)^n$.
\end{enumerate}
\item If $r$ divides $2^{\frac{2m + t - 1}{2}} + 1$, $ 2^{2m-1}D_4 - 2^m \mu - 1 \leq 2^{\frac{2m + t - 1}{2}} + 1 < 2^{\frac{3m - 1}{2}} + 1$ must hold. That is, $D_4 < \frac{1}{2^{\frac{m-1}{2}}}  + \frac{1}{2^{m-1}} \mu + \frac{2}{2^{2m-1}} \leq  \frac{1}{2^{\frac{m-1}{2}}}  + \frac{1}{2^{m-1}} + \frac{2}{2^{2m-1}} $. \textcolor{black}{For $m = 2, \frac{1}{2^{\frac{m-1}{2}}}  + \frac{1}{2^{m-1}} + \frac{1}{2^{2m-2}} \approx  1.4571$, and for $m=3, \frac{1}{2^{\frac{m-1}{2}}}  + \frac{1}{2^{m-1}} + \frac{1}{2^{2m-2}} = 0.8125$.} \\
As $D_4 \geq 1$ and odd, and $m \geq 2$, the only possibility is $m = 2$, addressed in item (\ref{just_label1b}). $\blacksquare$ 
\item \textcolor{black}{To consider the case $t = m$, one rewrites the expression in item (\ref{just_label1}) as \\
$2^{g-2m-t}[(2^g)^{n-1} - \comb[n]{1} (2^g)^{n-2} +  \comb[n]{2} (2^g)^{n-3} - ... + n] = $ \\
$2^{m-t}[(2^m)^{r-4} + \comb[r]{1} (2^m)^{r-5}\mu + \comb[r]{2} (2^m)^{r-6} + ... + 2^m\comb[r]{r-5} + \comb[r]{r-4} \mu] + X$ \\
From item (\ref{igmr18062024}), $X$ is odd. $\comb[r]{r-4} = \frac{(2^mD_3-1)(2^{m-1}D_3-1)(2^mD_3-3)(2^{m-2}D_3-1)}{3}$ is odd if $m > 2$. If $t = m$, for parity of the above expression to hold, $g > 3m$ must be true, which is not permitted. If $t = m = 2$, $\comb[r]{r-4}$ is even since $D_3$ is odd (item \ref{chiggs3}). To maintain parity, again, $g$ must be $3m = 6$ which is not permitted (item \ref{igmr18062024b}). }

\end{enumerate}
\end{enumerate}

\item $3m > g-m > 2m$
\begin{enumerate}
\item \label{just_label2} Noting that $2^m \mid \comb[r]{r-2}\mu + D_3$ (item \ref{just_label1}),  in \\ \\
$2^{g-m}[(2^g)^{n-1} - \comb[n]{1} (2^g)^{n-2} +  \comb[n]{2} (2^g)^{n-3} - ... + n] = $ \\
$2^m[(2^m)^{r-3} + \comb[r]{1} (2^m)^{r-4}\mu + \comb[r]{2} (2^m)^{r-5} + ... + \comb[r]{r-4}2^m \mu + \comb[r]{r-3}+ \frac{\comb[r]{r-2}\mu + D_3}{2^m} ]$ \\ \\
as $g - m > 2m$, $2^m$ must divide $\comb[r]{r-3}+ \frac{\comb[r]{r-2}\mu + D_3}{2^m} $. 
\item \label{just_label2a} Considering $(\bmod \; 2^m)$, using $r = 2^mD_3 - 1$ (item \ref{Maa_bobby1}) and $D_3 = 2^{m-1}D_4 - \mu$ (item \ref{just_label1a}) where $D_3$ and $D_4$ are both, odd and positive integers, as $m \geq 2$,\\ \\
$\comb[r]{r-3}+ \frac{\comb[r]{r-2}\mu + D_3}{2^m}  = \frac{(2^mD_3 - 1)(2^{m-1}D_3 - 1)(2^mD_3 - 3)}{3} + \frac{(2^mD_3 - 1)(2^{m-1}D_3 - 1) \mu + D_3}{2^m}$\\
$\equiv \frac{(-1) (2^{m-1}D_3 - 1)(-3)}{3} + 2^{m-1}D_3^2 \mu - D_3 \mu - \frac{D_3 \mu}{2} + \frac{D_3 + \mu}{2^m}$ \\
$\equiv 2^{m-1}\left(2^{m-1}D_4 - \mu \right) - 1 + 2^{m-1} \left(2^{2m-2}D_4^2 - 2^m D_4\mu + 1\right) \mu - \frac{3}{2}D_3 \mu + \frac{D_3 + \mu}{2^m}$\\
$\equiv -2^{m-1} \mu -1 + 2^{m-1} \mu -2 (2^{m-1}D_4 - \mu) \mu + \frac{D_3 \mu}{2} + \frac{2^{m-1}D_4}{2^m}$ \\
$\equiv  -1  -2 (2^{m-1}D_4 - \mu)\mu + \frac{(2^{m-1}D_4 - \mu) \mu}{2} + \frac{D_4}{2} $ \\
$\equiv -1  + 2 + \frac{(2^{m-1}D_4 - \mu) \mu}{2} + \frac{D_4}{2}  \;\;\;\; = \frac{2^{m-1}D_4 \mu + D_4 + 1}{2}$. \\
Thus, for $2^m$ to divide $\comb[r]{r-3}+ \frac{\comb[r]{r-2}\mu + D_3}{2^m}$, $D_4 + 1$ must be $2^{m-1}D_5$, $D_5$ an odd positive integer such that $D_4 \mu + D_5$ is even and divisible by 4. Then, \\
$r = 2^mD_3 - 1 = 2^{2m-1}D_4 - 2^m \mu - 1 = 2^{2m-1}(2^{m-1}D_5 - 1) - 2^m \mu - 1 = 2^{3m-2}D_5 - 2^{2m-1} - 2^m \mu - 1$.\\

\item {\bf Proposition}: Let $\comb[r]{r-4} \mu + \frac{1}{2^m}\left\{ \comb[r]{r-3}+ \frac{\comb[r]{r-2}\mu + D_3}{2^m} \right\} $ be precisely divisible by $2^t$, $t \geq 1$ an integer and $t < m$, so that 
the quotient is odd. Then, $(2^g - 1)^n = \frac{ (2^m + \mu)^r - \mu }{2^m} $ provides no solution. \\

{\bf Proof}:
\begin{enumerate}
\item $g$ must be $3m + t$ with $t$ such that $3m + t$ is an odd prime (item \ref{Kaali_M53_aaa}). 
\item \label{just_label2b} Per item (\ref{Kaali_M53abcdef}), $r \mid 2^{\frac{g-1}{2}} - \delta$, $\delta = 1$ or $-1$.\\ 
Thus, if $r = 2^{3m-2}D_5 - 2^{2m-1} - 2^m \mu - 1$ (item \ref{just_label2a}) divides $2^{\frac{3m + t - 1}{2}} - \delta$, \\
$2^{3m-2}D_5 - 2^{2m-1} - 2^m \mu - 1 \leq 2^{\frac{3m + t - 1}{2}} - \delta$ must hold. Noting that $m \geq 2$ and $t < m$, \\
$2^{3m-2}D_5 - 2^{2m-1} - 2^m \mu - 1 \leq 2^{\frac{3m + t - 1}{2}} - \delta < 2^{\frac{4m - 1}{2}} - \delta$, or, \\
$D_5 < \frac{1}{2^{m - \frac{3}{2}}} + \frac{1}{2^{2m - 2}} \mu + \frac{1}{2^{m-1}} + \frac{1 - \delta}{2^{3m-2}}
\leq \frac{1}{2^{m - \frac{3}{2}}} + \frac{1}{2^{2m - 2}} + \frac{1}{2^{m-1}} + \frac{1}{2^{3m-3}}
\leq \frac{1}{\sqrt{2}} + \frac{1}{2} + \frac{1}{4} + \frac{1}{8} = 1.5821$. \\ \\
For $m > 2$, $D_5 \approx 0.6817 < 1$, not possible since $D_5$ is an odd positive integer (item \ref{just_label2a}). $m$ must be $2$ and $D_5 = 1$. $t = 1$ making $g = 3\cdot 2 + 1 = 7$. $r = 2^{4} - 2^{3} - 2^2 \mu - 1$ which is $3$ for $\mu = 1$ and $11$ for $\mu = -1$. $\frac{(2^m + \mu)^r - \mu}{2^m}$ becomes $\frac{5^3 - 1}{4}$ and $\frac{3^{11} + 1}{4}$ respectively, and neither is of the form $(2^7 - 1)^n$. $\blacksquare$
\end{enumerate}
\end{enumerate}

\item $4m > g-m > 3m$
\begin{enumerate}
\item \label{just_label3} In \\ \\
$2^{g-m}[(2^g)^{n-1} - \comb[n]{1} (2^g)^{n-2} +  \comb[n]{2} (2^g)^{n-3} - ... + n] = $ \\
$2^{2m}[(2^m)^{r-4} + \comb[r]{1} (2^m)^{r-5}\mu +  ... + \comb[r]{r-5}2^m + \comb[r]{r-4} \mu  + \frac{1}{2^m}\left\{ \comb[r]{r-3}+ \frac{\comb[r]{r-2}\mu + D_3}{2^m} \right\} ]$\\ \\
(item \ref{just_label2}), $2^m$ must divide $\comb[r]{r-4} \mu  + \frac{1}{2^m}\left\{ \comb[r]{r-3}+ \frac{\comb[r]{r-2}\mu + D_3}{2^m} \right\} $. 
\item From items (\ref{just_label1a}) and (\ref{just_label2a}), using $r = 2^mD_3 - 1$, $D_3 = 2^{m-1}D_4 - \mu$ and $D_4 = 2^{m-1}D_5 - 1$ with $D_4$ and $D_5$ both positive and odd integers, \\ \\
$\comb[r]{r-2}\mu + D_3 =  (2^mD_3 - 1)(2^{m-1}D_3 - 1)\mu + D_3 $ \\
$= 2^m\left[ 2^{m-1}D_3^2 \mu + \frac{D_4 - 3D_3 \mu}{2}\right] = 2^m\left[ 2^{m-1}D_3^2 \mu + \frac{D_4 - 3 \mu 2^{m-1}D_4  + 3}{2} \right]$, and \\ \\
$\comb[r]{r-3} = \frac{(2^mD_3 - 1)(2^{m-1}D_3 - 1)(2^mD_3 - 3)}{3}$\\
$= \frac{1}{3}2^m \left( 2^{2m-1}D_3^3 + D_3 \right) - 2^m \left( 2^mD_3^2 - D_3 \right) + 2^{m-1}D_3 - 1$ so that \\ \\

$\frac{1}{2^m}\left\{ \comb[r]{r-3}+ \frac{\comb[r]{r-2}\mu + D_3}{2^m} \right\} = $\\
$\frac{1}{3}\left( 2^{2m-1}D_3^3 + D_3 \right) -2^mD_3^2 + \frac{D_3}{2}(D_3 \mu + 3) + \frac{D_5 - 3 \mu D_4}{4}$. 
\item Considering $(\bmod \; 2^m)$, noting that $m \geq 2$,
\begin{enumerate}
\item $\frac{D_5 - 3 \mu D_4}{4} \equiv \frac{D_5 + \mu 2^{m-1}D_5 - \mu }{4} - 2^{m-1}D_5\mu + \mu$ $(\bmod \; 2^m)$,
\item $\frac{D_3}{2}(D_3 \mu + 3) = 2^m\left[ 2^{3m-5}D_5^2 \mu - 2^{2m-3}D_5 \mu + 2^{m-2}\left( \frac{D_5 + \mu}{2} \right) \right] - 2^{m-2} - \mu$ \\
$\equiv - 2^{m-2} - \mu $ $(\bmod \; 2^m)$,
\item $2^mD_3^2 \equiv 0 $ $(\bmod \; 2^m)$,
\item $\frac{1}{3}\left( 2^{2m-1}D_3^3 + D_3 \right) = $ \\
$2^{m} \left[ -2^{6m-6}D_5^2 - 2^{2m-2} + 2^{5m-5}D_5 - 2^{3m-3} \mu + 2^{4m-3}D_5 \mu + 2^{3m-3}D_5 \right] $\\
$- 2^m \left[ 2^{5m-5}D_5^2 \mu \right] + \frac{1}{3}\left[ 2^m \left( 2^{7m-7}D_5^3 - 2^{4m-4} - 2^{m-1} \mu + 2^{m-2}D_5 \right) - \left( 2^{m-1} + \mu \right) \right]$ \\
$ \equiv \frac{1}{3}\left[ 2^m \left( 2^{7m-7}D_5^3 - 2^{4m-4} - 2^{m-1} \mu + 2^{m-2}D_5 \right) - \left( 2^{m-1} + \mu \right) \right]$ $(\bmod \; 2^m)$. \\ \\
Noting that\footnote{(a). If $m > 2$ is even, $\mu = 1$. $2^{m-1} + \mu \equiv -1 + 1 = 0 \; (\bmod \; 3)$.  (b). If $m > 2$ is odd, $\mu = -1$. $2^{m-1} + \mu \equiv 1 - 1 = 0 \; (\bmod \; 3)$. (c) If $m = 2$ and $\mu = 1$, $2^{m-1} + \mu = 3$. (d) $m = 2$ and $\mu = -1$ is not a possibility since $q = 2^m - \mu = 1$ is not permitted.} $3 \mid \left( 2^{m-1} + \mu \right)$ for $m \geq 2$ and $\mu = \pm 1$, $3$ must divide \\
$\left( 2^{7m-7}D_5^3 - 2^{4m-4} - 2^{m-1} \mu + 2^{m-2}D_5 \right)$. Thus, \\
$\frac{1}{3}\left( 2^{2m-1}D_3^3 + D_3 \right) \equiv - \frac{1}{3}(2^{m-1} + \mu )$ $(\bmod \; 2^m)$. \\
\item $\comb[r]{r-4}\mu = \frac{(2^mD_3 - 1)(2^{m-1}D_3 - 1)(2^mD_3 - 3)(2^{m-2}D_3 - 1)\mu}{3}$\\
$\equiv \frac{(-1)(2^{m-1}D_3 - 1)(-3)(2^{m-2}D_3 - 1)\mu}{3} = \left[ 2^{2m-3}D_3^2 - 3 \cdot 2^{m-2}D_3 + 1 \right] \mu$\\
$ = \left[ 2^{2m-3}(2^{2m-2}D_4^2 - 2^{m}D_4 \mu + 1) - 3 \cdot 2^{m-2}(2^{m-1}D_4 - \mu) + 1 \right] \mu $\\
$\equiv \left[ 2^{2m-3} -  2^{m}(2^{m-1}D_4 - \mu) + 2^{m-2}(2^{m-1}D_4 - \mu) + 1\right]\mu$\\
$\equiv \left[ 2^{2m-3} + 2^{2m-3}D_4  - 2^{m-2}\mu + 1\right]\mu = \left[2^{2m-2}(\frac{1 + D_4}{2})  - 2^{m-2}\mu + 1 \right]\mu$\\
$\equiv - 2^{m-2} + \mu \; (\; \bmod \; 2^m)$.
\end{enumerate}
\item \label{just_label3aa} Thus, $\comb[r]{r-4} \mu  + \frac{1}{2^m}\left\{ \comb[r]{r-3}+ \frac{\comb[r]{r-2}\mu + D_3}{2^m} \right\} $ $(\bmod \; 2^m)$ is\\
$\equiv \frac{D_5 + \mu 2^{m-1}D_5 - \mu }{4} - 2^{m-1}D_5\mu + \mu  - 2^{m-2} - \mu - \frac{1}{3}(2^{m-1} + \mu )  - 2^{m-2} + \mu$\\
$= \frac{D_5 + \mu 2^{m-1}D_5 + 3\mu }{4} - 2^{m}(\frac{D_5\mu  + 1}{2}) - \frac{1}{3}(2^{m-1} + \mu )$\\
$\equiv \frac{D_5 + \mu 2^{m-1}D_5 + 3\mu }{4} - \frac{1}{3}(2^{m-1} + \mu ) = \frac{1}{12}\left[(3D_5 + 5 \mu) + 2^{m-1}(3D_5 \mu - 4) \right]$. \\

\item \label{igmr19062024} If $2^m \mid \comb[r]{r-4} \mu  + \frac{1}{2^m}\left\{ \comb[r]{r-3}+ \frac{\comb[r]{r-2}\mu + D_3}{2^m} \right\} $, $3D_5 + 5 \mu = 2^{m-1}X$ must hold where $X$ is an odd integer (since $3D_5 \mu - 4$ is odd and $\frac{X + 3D_5 \mu - 4}{12}$ must be even). Setting $X = 2Y + 1$, $Y$ an integer, yields $D_5 = \frac{2^mY}{3} + \frac{2^{m-1} + \mu}{3} - 2 \mu$, noting that $3 \mid 2^{m-1} + \mu$. With $Y = 3D_6$, $D_6$ an integer, $D_5 = 2^mD_6 + \frac{2^{m-1} + \mu}{3} - 2 \mu$.  \\
Since $D_5 \geq 1$ (item \ref{just_label2a}), noting that $\mu = \pm 1$ and $m \geq 2$, \\
$D_6 \geq \frac{1 + 2 \mu}{2^m} - \frac{1}{3} \frac{2^{m-1} + \mu}{2^m} =  -\frac{1}{6} + \frac{1}{2^m}\left( 1 + \frac{5 \mu}{3} \right) \geq  -\frac{1}{6} - \frac{1}{2^m}\left(\frac{2}{3} \right) > -1$. \\

\item $r = 2^mD_3 - 1 = 2^m\left[ 2^{m-1}D_4 - \mu \right] - 1 = 2^m\left[ 2^{m-1}\left\{ 2^{m-1}D_5 - 1 \right\} - \mu \right] - 1$\\
$= 2^m\left[ 2^{m-1}\left\{ 2^{m-1}(2^mD_6 + \frac{2^{m-1} + \mu}{3} - 2 \mu) - 1 \right\} - \mu \right] - 1$\\
$ = 2^{4m-2}D_6 + 2^{3m-2}\left(\frac{2^{m-1} + \mu}{3} \right) - 2^{3m-1} \mu - 2^{2m-1} - 2^m \mu  - 1$. 

\item {\bf Proposition}: Let $ \comb[r]{r-5} + \frac{1}{2^m}\left[ \comb[r]{r-4} \mu  + \frac{1}{2^m}\left\{ \comb[r]{r-3}+ \frac{\comb[r]{r-2}\mu + D_3}{2^m} \right\}  \right]$ be precisely divisible by $2^t$, $1 \leq t < m$ an integer, so that the quotient is odd.  $(2^g - 1)^n = \frac{ (2^m + \mu)^r - \mu }{2^m} $ provides no solution. \\
{\bf Proof}:
\begin{enumerate}
\item $g = 4m + t$ must hold with $t$ such that $4m + t$ is an odd prime. 
\item \label{just_label3d} As in item (\ref{just_label2b}), for $\delta = 1$ or $-1$, $r \leq 2^{\frac{g-1}{2}} - \delta = 2^{\frac{4m + t -1}{2}} - \delta < 2^{\frac{5m -1}{2}} - \delta$ yields\\
$D_6 < \frac{1}{2^{\frac{3m-3}{2}}} + \frac{1}{2^{2m-1}} - \frac{1}{6} + \mu\left(\frac{1}{2^{3m-2}} + \frac{1}{2^{m-1}} - \frac{1}{3} \frac{1}{2^m} \right) + \frac{1 - \delta}{2^{4m-2}}\leq 0.8223$ (limit obtained for $m=2, \delta = 1, \mu = 1$).
\item \label{just_label3ab} Thus, $D_6 < 1$. Per item (\ref{igmr19062024}), $D_6 > -1$. $D_6 = 0$ must hold making $D_5 = \frac{2^{m-1} + \mu}{3} - 2 \mu$ (item \ref{igmr19062024}) and thus $r = 2^{3m-2}\left(\frac{2^{m-1} + \mu}{3} \right) - 2^{3m-1} \mu - 2^{2m-1} - 2^m \mu  - 1$.  

\item $r = 2^{3m-2}\left(\frac{2^{m-1} + \mu}{3} \right) - 2^{3m-1} \mu - 2^{2m-1} - 2^m \mu  - 1 \leq 2^{\frac{5m-1}{2}} - \delta $, $\delta = 1$ or $-1$, further yields \\
$\frac{2^{m-1} + \mu}{3} \leq \frac{1}{2^\frac{m-3}{2}} + \frac{1}{2^{m-1}} + \mu \left( 2 + \frac{1}{2^{2m-2}} \right) + \frac{1 - \delta}{2^{3m - 2}}$. \\
For $\mu = -1$, the right hand side of the inequality is negative for $m \geq 2$ and $\delta = 1$ or $-1$. \\
For $\mu = 1$, only $m=2$ and $m = 3$ satisfy the above inequality. But $m$ cannot be $3$ since $q = 2^m + 1$ is prime. \\
For $m = 2$ ($t = 1$), $r = 2^{4}\left(\frac{2^{1} + 1}{3} \right) - 2^{5} - 2^{3} - 2^2  - 1 = -29 < 0$, a contradiction, since $r$ must be positive and an odd prime (\S\ref{2mpn_plus_mu_equals_qr_r_odd}, item \ref{Kaali_M18}).

\item $ \comb[r]{r-5} + \frac{1}{2^m}\left[ \comb[r]{r-4} \mu  + \frac{1}{2^m}\left\{ \comb[r]{r-3}+ \frac{\comb[r]{r-2}\mu + D_3}{2^m} \right\}  \right]$ is therefore not divisible by $2^t$ for any $t \geq 1$. In particular, \\
$2^m$ does not divide $ \comb[r]{r-5} + \frac{1}{2^m}\left[ \comb[r]{r-4} \mu  + \frac{1}{2^m}\left\{ \comb[r]{r-3}+ \frac{\comb[r]{r-2}\mu + D_3}{2^m} \right\}  \right]$ and thus a solution to $(2^g - 1)^n = \frac{(2^m + \mu)^r - \mu}{2^m}$ is unlikely for $g > 4m$. $\blacksquare$
\end{enumerate}
\end{enumerate}
\end{enumerate}
\noindent \rule{16cm}{1pt}
\section{Concluding Remarks} 

For $p$ and $q$ as distinct Mersenne and/or Fermat primes, $\mu = \pm 1$ and $m, n, r$ as positive integers, per sections  \S\ref{2m_plus_mu_equals_pnqr}-\ref{2m_pm_plus_mu_equals_qr} and \cite{Bugeaud_1997}, number of integer solutions to (i) $2^{m} + \mu = p^{n}q^{r}$, (ii) $ 2^{m}p^{n} + \mu = q^{r}$ and (iii) $2^{m} = p^{n} + \mu q^{r}$ is shown to be finite using elementary methods and other approaches in the literature. 


When (iii) takes the form $2^{y+1} = (2^y + 1)^2 - (2^{2y} + 1)$, $y$ such that $2^y + 1$ and $2^{2y} + 1$ are primes, $rad(ABC)^{1 + \varepsilon} > C$ holds for any $\varepsilon > 0$. Barring the solutions for (iii) for $n$ even, $m$ odd, $r > 1$ odd, which are shown to be finite in number \cite{Bugeaud_1997}, solutions to all other diophantine equations above are enumerated in Appendix \S \ref{appendix_summary_of_solutions}. 

The generic case, where the triplet $\{A, B, C\}$ constitutes the product $ABC = 2^mx^ny^r$ with $x$ and $y$ as any two odd positive coprime integers, and $m$, $n$, $r$ as positive integers exponents satisfying any of the six diophantine combinations $\mbox{(a): } 2^{m} \pm 1 = x^{n}y^{r}, \mbox{ (b): } 2^{m} = x^{n} \pm y^{r}, \mbox{ and (c): }    2^{m}x^{n} \pm 1 = y^{r}$, still requires attention. \\

\noindent\rule{16cm}{1pt}  


\section*{Acknowledgment}{The author is sincerely indebted to Reese Scott, Robert Styer and Michel Waldschmidt for their generous guidance and support. Solutions to some diophantine equations (page \pageref{Maa_bobby153a}, item \ref{Maa_bobby153a}; page \pageref{Maa_17}, item \ref{Maa_17}; page \pageref{igmr24062024}, (footnote); and page \pageref{Maa_19}, item \ref{Maa_19}) are contributed by Reese Scott and Robert Styer. The author gratefully acknowledges discussions with K. B. Subramanium, K. V. Srikanth (IIT Guwahati), (Late) Arbind Lal, Alok Maloo. Special thanks to Shyam Sunder Nishad and Vitthal Khatik for patiently sitting through a series of presentations made by the author, and for helping him with corrections. Significant part of this work was carried out at RWTH Aachen University (many thanks to Burkhard Corves and Roger A. Sauer for generously hosting the author since 2010),  and IIT Guwahati.} \\
\noindent \rule{16cm}{1pt}

\noindent \rule{16cm}{1pt}
\begin{appendix}
\section{Solutions to (i) $2^{m}p+ \mu = p^{n}q^{r}$, (ii) $2^{m} = p^{n} + \mu q^{r}$ and (iii) $ 2^{m}p^{n} + \mu = q^{r}$}
\label{appendix_summary_of_solutions}
Cases corresponding to $2^{m} = p^{n} \pm q^{r}$, $n$ even, and $m$ and $r$ $(> 1)$ odd, are not listed. Finiteness of the number of solutions to this diophantine form is established in \cite{Bugeaud_1997}. $\varepsilon_o = \frac{log(C)}{log\left[rad(ABC)\right]} - 1$ is also computed so that $C > rad(ABC)^{1 + \varepsilon}$ is true for $0 < \varepsilon < \varepsilon_o$.

\begin{table}[h!]
 \label{tab:table1}
  \begin{center}
   \footnotesize
    \begin{tabular}{| c | c || c | c |}
     \hline  \hline
      $C = A + B $ & page & $C = A + B $ & page \\
     \hline
     \hline    
$3^3 \cdot 19 = 2^9 + 1$ & \pageref{concl1} &  $2^6 = 3^2 \cdot 7 + 1$ & \pageref{concl2} \\
$3^4 = 2^5 + 7^2$ & \pageref{concl3} &  $5^2  = 2^4 + 3^2$ & \pageref{concl4}, \pageref{concl5}, \pageref{concl8}\\ 
$3^3 = 2 + 5^2 $  & \pageref{concl5}, \pageref{concl9}, \pageref{concl16}, \pageref{concl18} & $3^2  = 2^2 + 5$  & \pageref{concl5} \\
$2^3 = 3 + 5 $   & \pageref{concl5} &  $2^5 = 3^3 + 5$  & \pageref{concl5} \\
$2^7 = 3 + 5^3$  & \pageref{concl5} &  $5 = 3 + 2$  & \pageref{concl5}, \pageref{concl10}  \\
$2^4 = 7 + 3^2$ & \pageref{concl6} &  $2^5 = 5^2 + 7$  & \pageref{concl7}\\
$7 = 3 + 2^2$ & \pageref{concl8pt5} &  $31 = 3^3 + 2^2$  & \pageref{concl8pt5} \\
$3^2 = 5 + 2^2$ & \pageref{concl11}, \pageref{concl13} & $3^4  = 2^6 + 17 $  & \pageref{concl12} \\
$(2^{y} + 1)^2 = 2^{y + 1} + (2^{2y} + 1)$ & \pageref{concl14} & $7^2 = 2^5 + 17$ & \pageref{concl15} \\
$17 = 2^3 + 3^2$ & \pageref{concl15a} & $17^{2} = 2^{5}3^{2} + 1$  & \pageref{concl19} \\
$3^4 = 2^4\cdot5 +1$  & \pageref{concl20} & $7^2 = 2^4 \cdot 3 + 1$ & \pageref{concl21}\\ 
$5^2 = 2^3 \cdot 3 + 1$ & \pageref{concl22} &   $2\cdot 5^1 = 3^2 + 1$ & \pageref{new_soln_eva1} \\
$2^2 \cdot 7 = 3^3 + 1$    & \pageref{concl23} & $ 5^3 =  2^2 \cdot 31 + 1$ & \pageref{concl24} \\
$7^1 + 2^1 = 3^2$  & \pageref{concl14pt5} &   &  \\

\hline
\end{tabular}
\end{center}
\end{table}

\noindent \rule{16cm}{1pt}

\section{GCD of factors of $g^{au} + \mu h^{av}, \mu = \pm 1, a > 1$ and odd (prime)}
\label{gcd_of_factors2}
\noindent Let $g$ and $h$ be any (positive or negative) integers. Let $u,v \geq 1$.\\
\noindent If $g$ and $h$ are coprime, $gcd(g^u + \mu h^v, \frac{g^{au} + \mu h^{av}}{g^u + \mu h^v}) = gcd(g^u + \mu h^v, a)$, as shown below\footnote{$gcd(x,y) = gcd(x,y + Dx)$ where $D$ is any integer.}.

\begin{align}
g^{au} + \mu h^{av} &= (g^u + \mu h^v)R \nonumber \\
& = (g^u + \mu h^v)\left[g^{u(a-1)} - \mu g^{u(a-2)}h^{v}  + ... - \mu g^{u}h^{v(a-2)} + h^{v(a-1)} \right] \nonumber \\
gcd(g^u + \mu h^v,R) &= \nonumber \\
& gcd[g^u + \mu h^v, g^{u(a-1)} - \mu g^{u(a-2)}h^{v} +  ... + h^{v(a-1)}  - (g^u + \mu h^v)g^{u(a-2)} ] \nonumber \\
&= gcd(g^u + \mu h^v, -2\mu g^{u(a-2)}h^{v} + g^{u(a-3)}h^{2v} - ... - \mu g^{u}h^{v(a-2)} + h^{v(a-1)} ) \nonumber \\
&= gcd(g^u + \mu h^v, -2\mu g^{u(a-2)}h^{v} + g^{u(a-3)}h^{2v} - ... - \mu g^{u}h^{v(a-2)} + h^{v(a-1)} + \nonumber \\
& \,\,  2 \mu h^v (g^u + \mu h^v)g^{u(a-3)}  ) \nonumber \\
&= gcd(g^u + \mu h^v, 3g^{u(a-3)}h^{2v} - \mu g^{u(a-4)}h^{3v} + ... - \mu g^{u}h^{v(a-2)} + h^{v(a-1)} ) \nonumber \\
&= ... \nonumber \\
&= gcd(g^u + \mu h^v, -(a-1)\mu g^u h^{v(a-2)} + h^{v(a-1)}) \nonumber \\
&= gcd(g^u + \mu h^v, -(a-1)\mu g^u h^{v(a-2)} + h^{v(a-1)} + (a-1)(g^u + \mu h^v)\mu h^{v(a-2)}) \nonumber \\
&= gcd(g^u + \mu h^v, h^{v(a-1)} + (a-1)h^{v(a-1)}) = gcd(g^u + \mu h^v,ah^{v(a-1)}) \nonumber
\end{align}

\noindent If $gcd(g,h) = 1$, $gcd(g^u,h^v) = gcd(g^u + 
\mu h^v,h^v) = gcd(g^u + \mu h^v, h^{v(a-1)}) = 1$ and so $gcd(g^u + \mu h^v,R) = gcd(g^u + \mu h^v, a)$. \\

\noindent {\bf Remark}: If $g, h \geq 2$ and $u, v \geq 1$ and $a \geq 3$, then $\frac{g^{au} + \mu h^{av}}{g^u + \mu h^v} > a$. 
\begin{enumerate}
\item One notes that $y^{x(a-1)} > a$, or, $xlog(y) > \frac{log(a)}{a-1}$ holds for integers $x \geq 1$, $y \geq 2$ and $a \geq 3$.
\item Consider $\mu = 1$. $\frac{g^{au} + h^{av}}{g^u + h^v} > a$, re-expressed as $g^u(g^{u(a-1)} - a) + h^v(h^{v(a-1)} - a) > 0$ holds per item 1 as individually, $g^u(g^{u(a-1)} - a) > 0$ and $h^v(h^{v(a-1)} - a) > 0$.

\item Consider $\mu = -1$. $\frac{g^{au} - h^{av}}{g^u - h^v} = g^{u(a-1)} + g^{u(a-2)}h^{v}  + ... + g^{u}h^{v(a-2)} + h^{v(a-1)}$ has all positive terms. $h^{v(a-1)} > a$ holds per item 1, and thus, $\frac{g^{au} - h^{av}}{g^u - h^v} > a$. 
\end{enumerate}
\noindent \rule{16cm}{1pt}

\section{$2^{m} + \mu, \mu = \pm 1,$ is not a perfect power}
\label{2_k_pm_mu}

A perfect power \cite{Waldschmidt_2009} is a positive integer $x^y$ with $x \geq 1$ and $y \geq 2$ as rational integers\footnote{Synonym for integers, used to distinguish from  cyclotomic, Eisenstein, Gaussian and Hamiltonian integers \cite{math_world}.}.

\begin{enumerate}
	\item Consider $d \geq 1$ an integer and $\mu = 1$. Let $2^{m} + 1$ (odd) be square so that \\ $ 2^{m} + 1 = (2d + 1)^2$, or, $2^{m-2} = d(d+1)$. $m > 2$ must hold. If $d > 1$, either $d$ or $d+1$ is odd and divides $2^{m-2}$, an impossibility. If $d = 1$, $m$ must be 3 so that $2^{3} + 1 = 3^2$, only exception to \S \ref{2_k_pm_mu}. 
	\item For $\mu = -1$, consider $2^{m} - 1$ (odd) as square, or, $ 2^{m} - 1 = (2d + 1)^2$ so that \\ $2^{m-1} = 2d(d+1) + 1$, an impossibility for $m > 1$ since the $LHS$ is even and $RHS$ odd.
	\item Let $2^{m} + \mu$ be an odd perfect power so that $2^m + \mu = (2d + \mu)^{r}, r > 1$ and odd. As $\mu^{j} = \mu$ for $j$ odd, $\mu^{k} = 1$ for $k$ even, and $^{r}C_{r-1}\ = r$,
	\begin{eqnarray}
		2^{m} = 2d\left[  (2d)^{r-1} + \comb[r]{1}(2d)^{r-2}\mu +  ... + \comb[r]{r-2}  (2d)\mu^{r-2} + r \right] = \left(2d \right) R \nonumber
	\end{eqnarray}
	Since $r $ is odd, $R$ is odd. $R \nmid 2^{m}$ unless $R = 1$, only possible when $r = 1$. 
\end{enumerate}

\subsection*{Corollary}
\label{conjecture_2p}

\noindent The ABC conjecture holds for the product $N = ABC = 2^m p^n$ composed of powers of only two distinct primes ($2$ and $p$) such that $A + B = C, gcd(A, B) = 1$ and $B > A \geq 1$.

\begin{enumerate}
	\item If $A > 1$, triple $\{A, B, C \}$, $gcd(A,B) = 1$, contributes at least three distinct primes to $rad(N = ABC)$. As $N$ must always be even, one among them is 2, and the other two are odd primes. 
	
	\item \label{Kaali_M4} Only possibility of $N$ being composed of powers of precisely two primes is when $A = 1$. As $N$ is divisible by only 2 and an odd prime $p$, following two possibilities arise:
	\subitem (i) $B = 2^m$ and $C = p^n$ so that $C = p^n = 2^m + 1$.
	\subitem (ii) $B = p^n$ and $C = 2^m$ implying $B = p^n = 2^m - 1$.
\end{enumerate}

\noindent Since $2^m + \mu$, $\mu = \pm 1$ is not a power (except when $m = 3 \text{ and } \mu = 1$), $n = 1$. \\ Inequalities $[2(2^m + 1)]^{1 + \varepsilon} > 2^m + 1 = C$ (for $m \geq 1$) and $[2(2^m - 1)]^{1 + \varepsilon} > 2^m = C$ (for $m \geq 2$) respectively, hold true for any $\varepsilon > 0$ suggesting that $rad(ABC)^{1 + \varepsilon} > C$ holds for any $\varepsilon > 0$ for any $p, m$ and $n$ satisfying $2^m + \mu = p^n$. \\
	

\end{appendix}


\begin{thebibliography}{9}

\bibitem{Masser85}
 Masser, D. W. (1985).
  \emph{"Open problems."}
  In Chen, W. W. L. Proceedings of the Symposium on Analytic Number Theory.
  London: Imperial College.
  
\bibitem{Oesterle88}
 Oesterl\'{e}, Joseph (1988).
  \emph{"Nouvelles approches du "th\'{e}or\`{e}me" de Fermat."}
  Ast\'{e}risque.
  S\'{e}minaire Bourbaki exp 694 (161): 165–186.

\bibitem{Beiler66}
 Beiler, A. H. (1966).
  \emph{"The Eternal Triangle." }
  Ch. 14. In Recreations in the Theory of Numbers: The Queen of Mathematics Entertains. NewYork. Dover, 1966.
  
\bibitem{Thms_Mersenne_primes}
  \emph{$https://en.wikipedia.org/wiki/Mersenne\_prime \# Theorems \_ about \_ Mersenne \_ numbers$}. accessed 22.05.2024.  
  
\bibitem{Bang1}
Bang. A. S. (1886). \emph{ "Taltheoretiske Undersøgelser"} .
Tidsskrift for Mathematik. 5. 4. Mathematica Scandinavica: 70–80.
   
\bibitem{Bang2}
   Bang. A. S. (1886). \emph{ "Taltheoretiske Undersøgelser "} (continued, see p. 80).
 Tidsskrift for Mathematik. 4: 130–137.
   
\bibitem{Lichtman_2025}
   Lichtman. J. D. (2025). \emph{ "The abc conjecture is almost always true"}.
  arXiv:2505.13991v1 [math.NT].

  \bibitem{Idowu_2025}
   Idowu. Michael. A. (2025). \emph{ "Symbolic Generation and Modular Embedding of High-Quality abc-Triples"}. arXiv:2506.10039v1 [cs.CR].

  \bibitem{Bright_2024}
 Bright. Curtis (2024). \emph{ "A New Lower Bound in the   Conjecture"}. 
Can. Math. Bull. 67 (2024) 369-378.

 \bibitem{Scoones_2022}
 Scoones Andrew(2022). \emph{ "On the   Conjecture in Algebraic Number Fields"}.   arXiv:2111.07791v2 [math.NT]




    

\bibitem{math_world}
Hardy. G. H. and Wright. E. M. (1979)
\emph{An Introduction to the Theory of Numbers}. 5th ed. Oxford, England: Clarendon Press,


\bibitem{Dorian_Goldfeld}
  Goldfeld Dorian (2002).
  \emph{Modular forms, elliptic curves, and the ABC-conjecture, A panorama in number theory: The view from Baker's garden. Based on a conference in honor of Alan Baker's 60th birthday, Cambridge University Press}
  Cambridge, 128–147.
    Mathematical Reviews (MathSciNet): MR1975449.
    Zentralblatt MATH: 1046.11035
    
\bibitem{wiki_abc}
Olivier Robert, Cameron L. Stewart and G\'{e}rald Tenenbaum. (2014)
\emph{A refinement of the abc conjecture}.
Bull. London Math. Soc. 46 (2014) 1156–1166.


    
\bibitem{Stewart_Tijdeman} 
Stewart, C. L. and Tijdeman, R. (1986).
\emph{On the Oesterl\'{e}-Masser conjecture.}
Monatshefte für Mathematik. 102 (3): 251–257. doi:10.1007/BF01294603.

\bibitem{Stewart_Yu1} 
Stewart, C. L. and Yu, Kunrui. (1991).
\emph{On the abc conjecture}. 
Mathematische Annalen. 291 (1): 225–230. doi:10.1007/BF01445201.

\bibitem{Stewart_Yu2}
Stewart, C. L. and Yu, Kunrui. (2001).
\emph{On the abc conjecture, II}.
Duke Mathematical Journal. 108 (1): 169–181. doi:10.1215/S0012-7094-01-10815-6.

\bibitem{Mochizuki1}
 Mochizuki, Shinichi. (2012a).
\emph{Inter-universal Teichm\"{u}ller Theory I: Construction of Hodge Theaters.}

\bibitem{Mochizuki2}
Mochizuki, Shinichi. (2012b).
\emph{Inter-universal Teichm\"{u}ller Theory II: Hodge–Arakelov-theoretic Evaluation.}

\bibitem{Mochizuki3}
 Mochizuki, Shinichi. (2012c).
\emph{Inter-universal Teichm\"{u}ller Theory III: Canonical Splittings of the Log-theta-lattice.}

\bibitem{Mochizuki4}
 Mochizuki, Shinichi. (2012d).
\emph{Inter-universal Teichm\"{u}ller Theory IV: Log-volume Computations and Set-theoretic Foundations.}

\bibitem{Mochizuki5}
Fesenko, Ivan. (2015).
\emph{Arithmetic deformation theory via arithmetic fundamental groups and nonarchimedean theta functions, notes on the work of Shinichi Mochizuki.}
Europ. J. Math., 1: 405–440.

\bibitem{Mochizuki6}
Shinichi Mochizuki, (2013, Dec.)
\emph{On the Verification of Inter-Universal Teichm\"{u}ller Theory: A Progress Report.}

\bibitem{Mochizuki7}
Shinichi Mochizuki, (2014, Dec.)
\emph{On the Verification of Inter-Universal Teichm\"{u}ller Theory: A Progress Report}

\bibitem{Wiles}
 Wiles, Andrew,  (1995).
\emph{Modular elliptic curves and Fermat's Last Theorem.}
Annals of Mathematics. 141 (3): 443–551. doi:10.2307/2118559. JSTOR 2118559. OCLC 37032255.

\bibitem{Ge_Yimin}
Ge, Yimin. (2008).
\emph{Elementary Properties of Cyclotomic Polynomials.}
Mathematical Reflections 2.

\bibitem{Birkhoff_Vandiver_1904}
Birkhoff, G. D. and Vandiver, H. S., (1904).
\emph{On the Integral Divisors of $a^n - b^n$.}
Annals of Mathematics. 5 (4):173-180.

\bibitem{Zsigmondy_1892}
Zsigmondy, K., (1892).
\emph{Zur Theorie der Potenzreste.}
Journal Monatshefte für Mathematik. 3 (1): 265–284.

\bibitem{Euler_maa}
Euler, L., (1758).
\emph{De numeris, qui sunt aggregata duorum quadratorum.}
Novi Commentarii academiae scientiarum Petropolitanae 4 (1758), pp. 3-40.
Translated by Bialek, Paul. R. $http://eulerarchive.maa.org/docs/translations/E228en.pdf$. accessed: 6th September, 2017, and 23rd May, 2024. 

\bibitem{Fermat_thm_sum_two_sqs1}
L. E. Dickson.
\emph{History of the Theory of Numbers}.
Vol II. Vol. 2. Chelsea Publishing Co., New York 1920;

\bibitem{Fermat_thm_sum_two_sqs2} \emph{Fermat's theorem on sums of two squares.} (2017).
https://en.wikipedia.org/wiki/Fermat's\textunderscore theorem\textunderscore on\textunderscore sums\textunderscore of\textunderscore two\textunderscore squares. accessed 7th September, 2017, and 23rd May, 2024.


\bibitem{Lucas_E} Lucas, Eduardo (1878) in \emph{Fermat Number}. $https://en.wikipedia.org/wiki/Fermat\textunderscore number$. accessed 29th September, 2017, and 24th May, 2024. 

\bibitem{Scott_Styer_2004}
Scott, R. and Styer, R. (2004).
\emph{On $p^x - q^y = c$ and related three term exponential Diophantine equations with prime bases.}
Journal of Number Theory 105 (2004) 212–234.

\bibitem{Scott_Styer_pc}
Scott, R. and Styer, R. (2018).
\emph{Solutions Contributed by Reese Scott and Robert Styer via personal communication with the author}. January 7th, 2018.

\bibitem{Scott_2018}
Scott, R. (2023).
\emph{Elementary treatment of $p^a \pm p^b + 1 = z^2$}. arXiv:math/0608796v1 [math.NT] 31 Aug 2006. revised 19th May, 2023. 

\bibitem{Scott_1993}
R. Scott. (1993).
\emph{On the Equations $p^x - b^y = c$ and $a^x + b^y = c^z$.}
Journal of Number Theory. 44 (2). 153–165.

\bibitem{Cao_1986}
Cao, Z. F. (1986).
\emph{On the equation $x^2 + 2^m = y^n$ and Hugh Edgar’s problem.}
 Kexue Tongbao. 31 (7). 555–556. (in Chinese).

\bibitem{Stroeker_and_Tijdeman_1982}
Stroeker, R.J. and Tijdeman, R. (1982).
\emph{Diophantine equations.}
Computational Methods in Number Theory. MC Track 155. Central Math Comp Sci, Amsterdam, 321–369.

\bibitem{Guy_1981}
Guy, R. K. (1981).
\emph{Unsolved Problems in Number Theory.}
Springer-Verlag. New York.

\bibitem{Waldschmidt_2009}
Waldschmidt, M. (2009).
\emph{Perfect Powers: Pillai's Works and their Developments.}
arXiv:0908.4031v1 [math.NT] 27 Aug 2009.

\bibitem{Pillai_SS_1945}
Pillai, S. S. (1945)
\emph{On the equation $2^{x} - 3^{y} = 2^{X} + 3^{Y}$.}
Bull. Calcutta Math. Soc., 37. 15–20. 

\bibitem{Catalan_Mihailescu}
Mih\u{a}ilescu, P. (2004).
\emph{Primary cyclotomic units and a proof of Catalan’s conjecture.}
J. Reine Angew. Math., 572. 167–195.

\bibitem{Bennett_Levin_2013}
Bennett, M. and Levin, A. (2015).
\emph{The Nagell–Ljunggren Equation via Runge's method.}
Monatsh Math. 177:15–31, DOI 10.1007/s00605-015-0748-1

\bibitem{Bender_Herzberg_1979}
Bender, E., Herzberg, N. (1979).
\emph{Some Diophantine equations related to the quadratic form $ax^2 + by^2$.} 
in: G.-C. Rota (Ed.), Studies in Algebra and Number Theory. Advances in Mathematics Supplementary Studies. Vol. 6. Academic Press, San Diego, 219-272.

\bibitem{Ljunggern}
Ljunggren, W. (1942).
\emph{Zur Theorie der Gleichung $x^2 + 1 = Dy^4$.}
Avh. Norske Vid. Akad. Oslo. I., 1942 (5): 27, MR 0016375.

\bibitem{Ljunggern1}
Steiner, R. Tzanakis, N. (1991).
\emph{Simplifying the solution of Ljunggren's equation $X^2 + 1 = 2 Y^4$}. 
Journal of Number Theory. 37 (2). 123–132. 

\bibitem{Ljunggern2}
Draziotis, K. A., (2007).
\emph{The Ljunggren equation revisited.}
Colloquium Mathematicum. 109 (1): 9–11. MR 2308822.

\bibitem{Ljunggern3}
Zhengjun Cao, Lihua Liu2. 2017
\emph{An Elementary Proof for Ljunggren Equation}.
arXiv:1705.03011v3 [math.NT]. 22nd May.

\bibitem{Stroeker_Tijdeman_1982}
Stroeker, R. J. and Tijdeman. R. (1982).
\emph{Diophantine equations.} 
Computational methods in number theory. Part II, vol. 155of Math. Centre Tracts, Math. Centrum, Amsterdam, pp.321–369.

\bibitem{Bennett_Skinner_2004}
Bennett, M. A. and Skinner, C. (2004).
\emph{Ternary Diophantine equations via Galois representations and modular forms.}
Canad. J. Math. 56. 23–54.

\bibitem{Rabinowitz}
Rabinowitz, S. (1977).
\emph{The solution of $y^2 \pm 2^n = x^3$}. 
Proc. Amer. Math. Soc. 62 (1). 1–6.

\bibitem{gp1}
Scott, R. Styer, R. (2011).
\emph{The generalized Pillai equation $\pm r a^x \pm s b^y=c$}. 
Journal of Number Theory. 131. 1037–1047.

\bibitem{gp2}
Le, M., (1992). 
\emph{A note on the Diophantine equation $a x^m - b y^n = k$}.
Indag. Math. (N.S.) 3 185–191.

\bibitem{gp3}
Shorey, T. N. (1986).
\emph{On the equation $ax^m - by^n=k$}.
Indag. Math. 48 (3) 353–358.

\bibitem{gp4}
M. Bennett, M., (2001). 
\emph{On some exponential equations of S.S. Pillai}.
Canad. J. Math. 53 (5) 897–922.

\bibitem{Nagell1}
Nagell, T., (1920).
\emph{Des \'{e}quations ind\'{e}termin\'{e}es $x^2 + x + 1 = y^n$ et $x^2 + x + 1 = 3y^n$}.
Norsk Mat. Forenings Skr. ser. 1 nr. 2, 14 pages.

\bibitem{Nagell2}
Nagell, T.
\emph{Note sur l’\'{e}quation ind\'{e}termin\'{e}e $\frac{x^n - 1}{x - 1} = y^q$}.
 Norsk Mat. Tidsskr. 2, 75–78.

\bibitem{Ljunggern4}
Ljunggren, W. (1943)
\emph{Noen Setninger om ubestemte likninger av formen $\frac{x^n - 1}{x - 1} = y^q$},
Norsk Mat. Tidsskrift. 25., 17–20.

\bibitem{Bugeaud_Mihailescu}
Bugeaud, Y. and Mih\u{a}ilescu, P. (2007)
\emph{On the Nagell-Ljunggern Equation $\frac{x^n - 1}{x - 1} = y^q$},
Math. Scand. 101. 177–183

\bibitem{Mahler}
Mahler, K. (1935) 
\emph{\"{U}ber den gr\"{o}ssten Primteiler spezieller Polynome zweiten Grades}
Archiv Math. og. Naturv. 41. 6. 3-26. 

\bibitem{da_Silva_et_al_2018}
N. da Silva, S. Raianu, and H. Salgado (2018) 
\emph{Differences of Harmonic Numbers and the abc-Conjecture}
The Pump Journal of Undergraduate Research. 1.  1-13.
 arXiv:1708.00620v1 [math.NT]

\bibitem{Lenstra_talk} 
H. Lenstra (1998).
\emph{Harmonic Numbers, Lecture at MSRI}.
available online at the URL: http://www.msri.org/realvideo/ln/msri/1998/mandm/lenstra/1/index.html

\bibitem{szalay} 
L. Szalay (2002).
\emph{The Equation $2^N \pm 2^M \pm 2^L = z^2$}.
Indag Math. N. S. 13 (1). 131-142.

\bibitem{Siksek_Stoll_2014} 
Samir Siksek and Michael Stoll (2014).
\emph{The generalised Fermat equation $x^2 + y^3 = z^{15}$}. 
Arch. Math. 102 (2014), 411-421.

\bibitem{Bugeaud_1997} 
Bugeaud Y. (1997).
\emph{On the diophantine equation $x^2 - 2^m = \pm y^{n}$}. 
Proc. Amer. Math. Soc. 125 (11) 3203-3208.

\bibitem{Bugeaud_Laurent_1996} 
Y. Bugeaud and M. Laurent. (1996).
\emph{Minoration effective de la distance p-adique entre puissances de nombres alg\'{e}briques}.
J. Number Th. 61. 311-342.

\bibitem{Laurent_Mignotte_Nesterenko_1995} 
M Laurent. M Mignotte. and Y. Nesterenko. (1995)
\emph{Formes lin\'{e}aires en deux logarithmes et d\'{e}terminants d'interpolation}.
J. Number Th. 55. 285-321.

\end{thebibliography}
\end{document}